\documentclass[a4paper,12pt]{book}

\usepackage{amssymb}
\usepackage{latexsym}
\usepackage{amsmath,amsthm,amsbsy}
\usepackage{xy}
\usepackage{xypic}

\newtheorem{thm}{Theorem}[section]
\newtheorem{cor}[thm]{Corollary}
\newtheorem{lem}[thm]{Lemma}
\newtheorem{prop}[thm]{Proposition}
\newtheorem{defn}[thm]{Definition}
\newtheorem{rem}[thm]{Remark}

\def\qed{\hfill$\Box$}

\def\CN{{\mathcal {N}}}
\def\CM{{\mathcal {M}}}
\def\CB{{\mathcal {B}}}
\def\CG{{\mathcal {G}}}
\def\CH{{\mathcal {H}}}
\def\CK{{\mathcal {K}}}
\def\CP{{\mathcal {P}}}
\def\CD{{\mathcal {D}}}
\def\CO{{\mathcal {O}}}
\def\CL{{\mathcal {L}}}
\def\CC{{\mathcal {C}}}
\def\CI{{\mathcal {I}}}

\def\bCK{{\mathbf {K}}}
\def\CE{{\mathcal {E}}}
\def\CF{{\mathcal {F}}}
\def\CR{{\mathcal {R}}}
\def\CA{{\mathcal {A}}}
\def\CG{{\mathcal {G}}}
\def\CT{{\mathcal {T}}}
\def\CU{{\mathcal {U}}}

\def\ua{{\underline{a}}}

\def\uz{{\underline{z}}}
\def\ut{{\underline{t}}}
\def\us{{\underline{s}}}
\def\uT{{\underline{T}}}
\def\uR{{\underline{R}}}
\def\uV{{\underline{V}}}
\def\uS{{\underline{S}}}
\def\uA{{\underline{A}}}
\def\uB{{\underline{B}}}
\def\uC{{\underline{C}}}

\def\uE{{\underline{E}}}

\def\uY{{\underline{Y}}}
\def\uZ{{\underline{Z}}}

\def\uL{{\underline{L}}}
\def\uR{{\underline{R}}}

\def\ua{{\underline{a}}}
\def\ub{{\underline{b}}}
\def\uc{{\underline{c}}}

\def\N{{\mathbb {N}}}
\def\C{{\mathbb {C}}}
\def\W{{\mathbb {W}}}
\def\F{{\mathbb {F}}}

\def\Z{{\mathbb {Z}}}
\def\D{{\mathbb {D}}}
\def\T{{\mathbb {T}}}

\newcommand{\eins}{\mathbf{1}}

\newcommand{\aso}{\stackrel{\circ}{a}^*}
\newcommand{\Ao}{\AA}
\newcommand{\Aso}{{\AA}^*}
\newcommand{\uAso}{\stackrel{\circ}{\underline{A}^*}}
\newcommand{\Ho}{\stackrel{\circ}{\CH}}

\newcommand{\uAo}{\stackrel{\circ}{\uA}}
\newcommand{\Do}{\stackrel{\circ}{D}}
\newcommand{\Dso}{\stackrel{\circ}{D}_*}
\newcommand{\CDo}{\stackrel{\circ}{\CD}}
\newcommand{\CDso}{\stackrel{\circ}{\CD}_*}
\newcommand{\Po}{\stackrel{\circ}{\CP}}
\newcommand{\Co}{\stackrel{\circ}{C}}

\newcommand{\oomega}{\overline{\omega}}
\newcommand{\uomega}{\underline{\omega}}

\newcommand{\Od}{\mathcal{O}_d}
\newcommand{\OH}{\Omega_{\CH}}
\newcommand{\OK}{\Omega_{\CK}}
\newcommand{\OP}{\Omega_{\CP}}
\newcommand{\DsA}{D_{*,A}}

\begin{document}

\thispagestyle{empty}

\begin{center} {\Large {\bf Characteristic Functions, Liftings and Modules

\vspace{1cm}

Habilitationsschrift}

\vspace{1cm}

 zur Erlangung des akademischen Grades\\
doctor rerum naturalium habilitatus\\

(Dr. rer. nat. habil.)\\
an der\\
Mathematisch-Naturwissenschaftlichen Fakult\"at\\
der Ernst-Moritz-Arndt-Universit\"at Greifswald

\vspace{2cm}

vorgelegt von Santanu Dey\\
geboren am 18. M\"arz 1976\\
in Bankura

\vspace{2cm}

Greifswald, 2009}

\end{center}

\include{erklaeren}

\frontmatter
\tableofcontents

\mainmatter
\chapter{Introduction}

This thesis specializes on certain topics in operator theory related
to dilation theory and to the theory of invariant subspaces. For any
contraction one associates  an isometry (or a unitary) unique up to
unitary equivalence called the minimal isometric (or unitary)
dilation.

Earlier pioneering works in dilation theory were by Sz.\,Nagy,
Stinespring and Kolmogorov. Since then it has been one of the
central techniques for working with interpolation, extension and
similarity problems, noncommutative probability, semigroups of
completely positive maps  and operator spaces. Surveys of such
applications can be found in \cite{FF90}, \cite{Ar69}, \cite{Pis01},
\cite{Pau03}, \cite{BP94}, \cite{Bh96} and \cite{ER00}.

Since late 90s there has been a lot of interest in dilation theory
of row contractions
(\cite{Ar98},\cite{DKS01},\cite{Po89a},\cite{Po02},\cite{BJKW00}). A
row contraction is a  tuple of operators on a common Hilbert space
which as a row operator is contractive. One instance when row
contractions appear is during Kraus decomposition of normal unital
completely positive maps. In place of the shift operator on $l^2$
(or $H^2$) appearing in the construction of dilation of a single
operator, the shift operators on the Fock space is needed in the
context of the dilation of a row contraction. Arveson used this type
of dilation to study some module structures induced by row
contractions.

In this thesis we investigate normal unital completely positive maps
on $B(\CH)$ for some Hilbert space $\CH$ with invariant vector states and
solve the classification problem of coisometric row contractions. There are many
recent works on non-coisometric case of this classification problem
and much of this case is well understood. We prove several
results about fixed points of completely positive maps on $B(\CH).$
The role of theories of fixed points and ergodicity is significant for von
Neumann algebras. A surprising observation of Arveson is that
certain normal unital completely positive maps on $B(\CH)$ produces geometric invariants.
Our above mentioned classification results are related to these invariants.
Finally we extend our theory to certain algebras on Hilbert modules which
includes analytic cross-products as special case.

Using any state or rather any positive linear functional on a given
C*-algebra $\CA,$ we get an embedding of $\CA$ as a subalgebra of
the algebra of bounded operators $B(\CH)$ for some  Hilbert space
$\CH$ through the GNS construction. Unital completely positive maps
are analogs of states when the codomain $\C$ is replaced by
 a  $C^*$-algebra. It was Stinespring who first realised
that  unital completely positive maps from C*-algebras to any
$B(\CH)$ can be  dilated to representations of those $C^*$-algebras.
Completely positive maps  have been used extensively in operator
algebras and its classification, and in mathematical physics. So the
problem of classifying completely positive maps is important.

Another significant result in the theory of invariant subspaces
 of interest to us here is Beurling's theorem. It states that
any invariant subspace of $H^2$ for the shift operator can be
realised as the range of an inner $H^\infty$ function. Sz.\,Nagy and
Foias classify certain big classes of contractions (namely the
completely noncounitary ones) up to unitary equivalence using  some
$H^\infty$ functions called the characteristic functions which are
identifiable with those from Beurling's theorem.

The notion of characteristic function and the related theory was
extended to row contractions by Popescu using his theory of dilation
of row contractions. He developed a noncommutative analogue of
$H^\infty$ for this purpose. Some analytic characteristic functions
in the sense of multivariable complex analysis for  commuting row
contractions were studied in \cite{BES05}, \cite{Po05}, \cite{BT07},
etc. A function intrinsically similar to this analytic function
already appeared in a preceding article (\cite{Ar98}) of Arveson.

The Cuntz algebras are primary example of simple nuclear purely
infinite $C^*$-algebras. Representations of Cuntz algebras are
obtained from minimal isometric dilations of coisometric row
contractions. These representations come very handy as a tool for
investigating endomorphisms. The commutant lifting theorem  of
Bratteli, Jorgensen, et al.\,and several other of their works
illustrates this. Further, the irreducible endomorphisms on $B(\CH)$
with unique invariant vector states via Kraus decomposition yield
 minimal isometric dilations of certain tuples of complex numbers
 corresponding to Cuntz states.

The weak Markov dilation (cf. \cite{BP94}) of an ergodic  completely
positive map (or discreet completely positive semigroup) on $B(\CH)$
for some $\CH$ with an invariant vector state is always an irreducible
endomorphism. While looking at the infinite tensor product picture
of weak Markov dilations, certain convergence results (quoted as
Appendix at the end of Chapter 3) had made us to look for (or
conjecture the existence of) multianalytic operators similar to
Sz.\,Nagy-Foias/Popescu's characteristic functions. We compute these
extended characteristic functions using techniques from
noncommutative probability in Section 3. They are complete unitary
invariants for a class of completely positive maps on $B(\CH).$ This
approach is motivated by the scattering theory for noncommutative
Markov chains of K\"ummerer and Maassen.

In chapter 4 we  reinterpret the above phenomenon of occurrence of
multianalytic invariants and thereby we notice that it has to do
with a much more general setup, namely that of subisometric
liftings. Such liftings for single contractions appeared first in a
work of Douglas and Foias \cite{DF84}. They established related
uniqueness and commutant lifting properties. In this context we
introduce characteristic functions of contractive liftings of row
contractions. They classify the certain liftings up to unitary
equivalence and provide a kind of functional model. The most general
setting is identified here which we call reduced liftings.  Finally
we discuss a factorization property of these multianalytic functions
and provide applications of our theory to completely positive maps.
Coisometric lifting is one of the most useful special case. It
relates via Kraus decomposition to the lifting of the corresponding normal unital
completely positive map. The commutant lifting theorem of Bratteli,
Jorgensen, Kishimoto and Werner gives an affine order isomorphism
between the fixed points of normal unital completely positive maps on some $B(\CH)$
and the commutants of the Cuntz algebra representations coming from
the associated Stinespring dilation. We apply our theory along with
this lifting theorem to prove a one-to-one correspondence between
the fixed point sets of any two normal unital completely positive
maps on algebras of bounded operators on Hilbert spaces
where one of the map is a coisometric lifting (or power
dilation) of the other by a $*$-stable row contraction.

Apart from the minimal isometric dilations there are other types of
dilations (see \cite{Po05}, \cite{Ar98}, \cite{BBD04}, \cite{BB02}
and \cite{BDZ06}) called  constrained dilations for row contractions
satisfying polynomial relations. An important class of the
constrained dilations is that of standard commuting dilation
initiated in the works of Drury \cite{Dr78} and Arveson \cite{Ar98}
and it comes with a multi-variable analytic functional calculus
model. Davidson \cite{DKS01} gave the complete structure of minimal
isometric dilations for row contractions on finite dimensional
spaces. Some of his techniques help us in chapter 5 to explore how
constrained dilations can be derived from minimal isometric
dilations in the context of row contractions on finite dimensional
Hilbert space. Our above theory of multi-analytic operators also extends
to constrained liftings.

A byproduct of the study of row contractions is the development of a
noncommutative multivariable complex analysis by Popescu. Many
results of classical theory like theorems of Cauchy, Liouville,
Schwartz, etc. have their noncommutative analogs. Here the open unit
disc in the complex plane is replaced by
$[B(\CH)^d]_1:=\{(X_1,\ldots,X_d) \in B(\CH)^d: \|X_1X^*_1  + \ldots
+ X_d X^*_d\| < 1\},$ for some Hilbert space $\CH.$

There after we  obtain that the two
invariants of Hilbert modules, namely, curvature invariant and Euler
characteristic are related to the characteristic function.
In the last section of chapter 6 we take a coordinate free approach
to our theory using Hilbert $C^*$-modules and their von Neumann
counterparts. Hilbert modules first appeared in the works of
Kaplansky while attempting to generalize Hilbert spaces by
structures where the inner products can take values in commutative
$C^*$-algebras. The motivation was to use the rich theory of vector
bundles in operator algebras. Later Paschke and Rieffel
independently extended this study to the case where the inner
product take values in arbitrary $C^*$-algebras. Today they are
extensively used in KK-theory, classification of $C^*$-algebras and
noncommutative dynamics.

Cuntz-Pimsner algebras are quotients of certain subalgebras (namely,
of Toeplitz algebras) of $C^*$-algebras of adjointable operators on
Hilbert bimodules. They contain crossed products by $\Z$ and
Cuntz-Krieger algebras as special cases. The von Neumann counterpart
of Toeplitz algebras are ``Hardy algebras". A notion of analytic
crossed product can be given using Hardy algebras. The Hardy
algebras were shown by Muhly and Solel to be generalizations of
$H^\infty$ in a  way similar to the noncommutative multivariable
complex analysis of Popescu mentioned above.

We have attempted to present the beautiful mathematical theory which
involves an interplay between the theory of invariant subspaces,
multivariable (noncommutative) complex analysis, the theory of
completely positive maps and dilation theory, and we discuss some of
its implications in operator algebras. The related theory has
developed significantly in recent years and we sincerely hope that
our monograph will be an invitation for the readers to pursue this
promising topic.

{\bf Acknowledgements:} I want to thank Rolf Gohm with whom I have
collaborated for two of the articles which are part of the present
thesis and for  permitting me to include them here. These two
articles are chapter 3 and 4 here. Discussions with him has been
very enriching and helped me to improve my mathematical skills. I am
indebted to M.\,Sch\"urmann for his encouraging me to work on this
monograph. I acknowledge the significant amount of Mathematics I
learned from the mathematical expertise of B.V.R.\,Bhat,
M.\,Sch\"urmann and U.\,Franz. I am grateful to
 G.\,Elliott, H.\,Osaka, J.\,Zacharias and many other with whom I
had discussions on these mathematical topics. Much of the work in
this monograph has been deeply influenced by the works of
W.\,Arveson, G.\,Popescu, R.\,Gohm, P.\,Muhly and B.\,Solel. I am
thankful to my colleagues at the Mathematics Department of the
University of Greifswald where most of the work was down. Research
work done during my stay at Fields Institute, Toronto in summer 2007
and at Ritsumeiken University Kyoto in summer 2008 has been very
useful for this thesis, and I gratefully acknowledge the financial
support I received for my visit at these Institutes.

\chapter[Preliminaries]{Dilations, Beurling's Theorem and Other Preliminaries}

\section{Minimal isometric dilations}

Isometries and unitaries have several preferred properties as
compared to  contractions. In dilation theory one utilizes the
existence of an isometric dilation of a given contraction. Consider
the simplest case of a contraction, namely a scalar operator $k$ on
$\C$ with $|k| \leq 1.$ The ``smallest" isometric dilation of $k$
will be the operator $\tilde{k} \in B(\C \oplus l^2(\N))$ given by
\[ \tilde{k} (h,h_1, h_2, \ldots): = (kh,(1 - |k|^2)^{\frac{1}{2}}h, h_1, h_2, \ldots) \]
where $h,h_i \in \C$ for  $i=1,2,\ldots$ and $\sum^d_{i=1}
|h_i|^2 < \infty.$ The isometric dilation is on a rather big space
as compared to the domain of the initial operator.

It is known that every contraction has isometric dilations and
uniqueness up to unitary equivalence of isometric dilations can be
ensured by making minimality assumption. Formally, we have:
\begin{defn}
Let $T$ be a contraction on a Hilbert space $\CH,$ i.e., $\|T\| \leq
1.$ An isometry $V$ on Hilbert space $\tilde{\CH}$ is an {\em
isometric dilation} of $T$  if
\[ \CH \subset \tilde{\CH} \mbox{~~and~~} P_{\CH} V^n |_{\CH} = T^n. \]
An isometric dilation is said to be minimal if
$\tilde{\CH}=\overline{\mbox{span}} \{ V^n h: n \in \N \cup \{0\}, h
\in \CH \}.$
\end{defn}
A minimal isometric dilation (mid for short) can be constructed in a
parallel way to the isometry constructed in the example above.

A tuple $\uT=(T_1,\ldots,T_d)$ of bounded operators on a Hilbert
space $\CH$ is said to be a row contraction if $T_1T^*_1 + \ldots +
T_d T^*_d \leq \eins.$ Treating a row contraction $\uT$ as a row
operator from $ \bigoplus^d_{i=1} \CH$ to $\CH,$ define $D_*:=
(\eins-\uT \uT^*)^\frac{1}{2}: \CH \rightarrow \CH$ and
$D:=(\eins-\uT^*\uT)^\frac{1}{2}: \bigoplus^d_{i=1} \CH \rightarrow
\bigoplus^d_{i=1} \CH$. This implies that
\begin{equation}
D_*= (\eins -\sum^d_{i=1} T_i T^*_i)^\frac{1}{2}, ~~~D=(\delta_{ij}
\eins -T_i^*T_j)^\frac{1}{2}_{d\times d}.
\end{equation}
Observe that $\uT D^2 = D^2_* \uT$ and hence $\uT D = D_* \uT.$ Let
$\CD:=\mbox{Range~} D$ and $\CD_*:=\mbox{Range~} D_*.$
\\

We  use the following {\em multi-index notation}. Let $\Lambda $
denote the set $\{ 1, 2, \ldots , d\}$ and
$\tilde{\Lambda}:=\cup_{n=0}^{\infty} \Lambda^n$, where $\Lambda
^0:=\{ 0\}$. If $\alpha \in \Lambda^n \subset \tilde{\Lambda}$ the
integer $n = |\alpha|$ is called its length. Now $T_\alpha$ with
$\alpha = (\alpha_1, \cdots ,\alpha_n) \in \Lambda^n$ means
$T_{\alpha_1} T_{\alpha_2} \ldots T_{\alpha_n}$.

The full Fock space over $\C^d$ ($d\geq 2$) denoted by $\Gamma
(\C^d)$ is the Hilbert space
$$\Gamma (\C^d)=\mathbb{C}\oplus \C^d \oplus (\C^d)^{\otimes ^2}\oplus \ldots
\oplus (\C^d) ^{\otimes ^m}\oplus \ldots.
$$
We will just write $\Gamma$ instead of $\Gamma (\C^d)$ at times. The
element $1\oplus 0\oplus \ldots$ of $\Gamma (\C^d)$ is called the
vacuum vector. Let $\{e_1, \ldots , e_d\}$ be the standard
orthonormal basis of $\C^d$. We include $d=\infty$ in which case
$\C^d$ stands for a complex separable Hilbert space of infinite
dimension. For $\alpha \in {\tilde \Lambda}$, $e_{\alpha }$ will
denote the vector $e_{\alpha _1}\otimes e_{\alpha _2}\otimes \ldots
\otimes e_{\alpha _m}$ in the full Fock space $\Gamma
(\mathbb{C}^d)$ and $e_0$ will denote the vacuum vector. Then the
(left) creation operators $L_i$'s on $\Gamma({\mathbb{C}^d})$ are
defined by
\[ L_i x = e_i \otimes x \]  for $1 \leq i \leq d$ and $x
\in \Gamma({\mathbb{C}}^d).$ The row contraction $\underline{L}=
(L_1, \ldots , L_d)$ consists of isometries with orthogonal ranges.
\\

Popescu in \cite{Po89a} gave the following explicit presentation of
the minimal isometric dilation of $\uT$ by $\uV$ on $\hat{\CH}=\CH
 \oplus (\Gamma (\C^d) \otimes \CD)$,
\begin{equation}
V_i(h \oplus \sum_{\alpha \in \tilde{\Lambda}}  e_\alpha \otimes
d_\alpha) = T_i h \oplus [e_0 \otimes D_i h + e_i \otimes
\sum_{\alpha \in \tilde{\Lambda}} e_\alpha \otimes d_\alpha]
\end{equation}
for $h \in \CH$ and $d_\alpha \in \CD.$ Here $D_i h := D
(0,\ldots,0,h,0,\ldots,0)$ and $h$ is embedded at the $i^{th}$
component. For $d=1$ this coincides with  Sch\"affer's construction.
\\

In other words, the $V_i$ are isometries with orthogonal ranges such
that $T^*_i = V^*_i |_{\CH}$ for $i=1,\ldots,d$ and the spaces
$V_\alpha \CH$ with $\alpha \in \tilde{\Lambda}$ together span the
Hilbert space on which the $V_i$ are defined. It is an important
fact, which we shall use repeatedly, that such minimal isometric
dilations are unique up to unitary equivalence (cf. \cite{Po89a}).
If in addition $\sum T_i T^*_i =\eins,$ then $\sum V_i V^*_i= \eins.$

For $d\geq 2$, the Cuntz algebra $\CO _d$ is the $C^*$-algebra
generated by $n$-isometries $\us= \{s_1, \ldots , s_d\}$, satisfying
Cuntz relations: $s_i^*s_j=\delta _{ij}\eins, 1\leq i, j\leq n,$ and
$\sum s_is_i^*=\eins.$ Therefore, $V_i$'s are representations of
Cuntz algebras, if $\sum T_i T^*_i =\eins.$

\section[Theorems of Stinespring and Beurling]
{Kraus decomposition, and theorems of Stinespring and Beurling}

Any normal unital completely positive map on $B(\CH)$ for any Hilbert space $\CH$ is associated to a coisometric
row contraction (may be of infinite length) on $\CH.$

\begin{thm} (Kraus) Let $\varphi: B(\CH) \to B(\CH)$ be a normal unital completely positive map.
Then there is a row contraction $\uT=(T_1,\ldots, T_d)$  such that
\[ \varphi (X) =\sum^d_{i=1} T_i X T^*_i, \mbox{~~for~} X \in B(\CH),\]
where $d$ may be infinite. In the $d=\infty$ case the convergence of
the above sum is to be taken in strong operator topology.
\end{thm}

The unital completely positive maps have proved to be the right
generalisations of states in the study of operator algebras. A role
parallel to the GNS construction is played by the following theorem
of Stinespring for completely positive maps:

\begin{thm}
If $\varphi$ is a unital completely positive map from a unital $C^*$-algebra $\CA$ to into $B(\CH),$ then there is a Hilbert space $\CK$
containing $\CH$ and $*$-representation $\sigma$ of $\CA$ on $\CK$ such that $\varphi(X)=P_{\CH} \sigma (X)|_{\CH}.$
\end{thm}

An elegant way of deriving Kraus decomposition from Stinespring
Theorem is described  in (\cite{Go04}, page 48).\\

Another result in operator theory crucial for us (cf. Section 4.1)
is the Beurling type theorem proved by Popescu (cf. \cite{Po89b}).
We give here a detailed proof as there is only a short discussion on
it in the existing literature. At first we quote the corresponding
classical theorem.  Let $\D$ denote the open unit ball of $\C.$ For a
Hilbert space $\CK$ we denote the set
\[ \{ u \hspace{.2cm} | \hspace{.2cm} u(z)=\sum^k_0 z^k a_k\mbox{~for ~} z \in \D, a_k \in \CK \mbox{~and~} \sup_{0 <r <1}
\int^{2 \pi}_0 \| u(r e^{it})\|^2 dt < \infty \}.\] by $H^2(\CK)$
and  define
\begin{align*}
H^\infty: = \{u: \overline{\D} \to \overline{\D} \hspace{.2cm} |
\hspace{.2cm} u \mbox{~is analytic
on~} \D, \mbox{~a.e. continuous on~} \partial \D,\\
\mbox{~and~}\sup_{z \in \D} |u(z)| < \infty \}.
\end{align*}
 Moreover, a function $u \in H^\infty$ is
called {\em inner} if $|u(e^{it})|=1$ a.e. for $t \in [0,2\pi).$

\begin{thm} (Beurling's Theorem)
The invariant subspaces $\CM$ for the multiplication operator $M_z$
on $H^2 (\CU),$ for some Hilbert space $\CU,$ are precisely those of
the form
\[ \CM = \theta H^2 (\CN)\]
where $\theta$ is an inner function with values in $B(\CN,\CU).$
\end{thm}

Let $\uV=(V_1,\ldots,V_d)$ be a row contraction consisting of
isometries on a Hilbert space $\CR.$ It follows that
$V_i,i=1,\ldots, d$ have mutually orthogonal ranges. A subspace
$\CK$ of $\CR$ is called {\em wandering} if
\[ V_\alpha \CK \bot V_\beta \CK \mbox{~~~~for distinct~} \alpha, \beta \in \tilde{\Lambda}.\]
Further, a row contraction $\uV$ of isometries is called a {\em
$d$-orthogonal shift} if there is a wandering subspace $\CK \subset
\CR,$ and
\[\CR= M(\CK):= \overline{\mbox{span}}\{V_\alpha \CK: \alpha \in \tilde{\Lambda}\}.\]
We have a unitary transformation $\Phi^\CK:M(\CK) \to \Gamma \otimes \CK$ defined by
\[\Phi^\CK (\sum_{\alpha \in \tilde{\Lambda}} V_\alpha k_\alpha)= \sum_{\alpha \in \tilde{\Lambda}} (L_\alpha \otimes \eins) (e_0 \otimes k_\alpha).\]

\begin{thm} (cf. \cite{Po89a}) (Wold decomposition)
Let $\uV=(V_1, \ldots,V_d)$ be a row contraction consisting of isometries on a common Hilbert space $\CR$. Then
$\CR$ decomposes as $\CR= \CR_0 \oplus \CR_1,$ such that $\CR_0$ and $\CR_1$ reduces $V_i, i=1,\ldots,d,$ and we have
\begin{itemize}
\item[(a)] $(\eins -\sum^d_{i=1} V_i V^*_i)|_{\CR_1}=0$ and $(V_1|_{\CR_0},\ldots,V_d|_{\CR_0})$ is a $d$-orthogonal shift for $\CR_0.$

\item[(b)] $\CR_1= \cap^\infty_{n=0} \overline{\mbox{span}}\{V_\alpha \CR: |\alpha| =n\}$ and $\CR_0= M(\CN),$ where $\CN:=\CR \ominus \overline{\mbox{span}}\{V_i \CR: i=1,\ldots,d\}.$
\end{itemize}
\end{thm}

For two Hilbert spaces $\CK$ and $\CK^\prime$ consider an operator $\theta: \CK \to \Gamma \otimes \CK^\prime.$  Take
$M_\theta : \Gamma \otimes \CK \to \Gamma \otimes \CK^\prime$ as the operator:
\[M_\theta (L_\alpha \otimes \eins) (e_0 \otimes k): = (L_\alpha \otimes \eins) \theta k, \hspace{.5cm} k \in \CK. \]
That $M_\theta (L_i \otimes \eins)= (L_i \otimes \eins) M_\theta$ is immediate. The operator $\theta$ is called {\em inner,} if
$M_\theta$ is an isometry.

\begin{lem}
Let $\uV$ and $\uV^\prime$ be two $d$-orthogonal shifts on Hilbert spaces $\CR$ and $\CR^\prime,$ with the wandering subspaces $\CK$
and $\CK^\prime$ respectively.   Let $Q$ be a contraction of $\CR$ into $\CR^\prime$ such that for $i =1,\ldots,d$
\[ Q V_i = V^\prime_i Q. \]
Then there is a contraction $\theta$ of $\CK$ into $\Gamma \otimes \CK^\prime$ such that
\begin{equation}
\Phi^{\CK^\prime} Q = M_\theta \Phi^\CK.
\end{equation}
If $Q$ is an isometry, then $\theta$ is inner.
\end{lem}

\begin{proof} Let $k$ be an arbitrary element in $\CK.$ Then from the
definition of $Q$ we get elements $k^\prime_\alpha \in \CK^\prime,$
such that
\[ Qk= \sum_{\alpha \in \tilde{\Lambda}} V^\prime_\alpha k^\prime_\alpha.   \]
It yields $\sum \|k^\prime_\alpha\| ^2 = \|Qk\|^2 \leq \|k\|^2.$ We set $k^\prime_\alpha = \theta_\alpha k$
for all $k \in \CK$ and thereby get $\theta_\alpha \in B(\CK,\CK^\prime).$ The function
\[\theta k := \sum_{\alpha \in \tilde{\Lambda}} e_\alpha \otimes \theta_\alpha k, \hspace{.5cm} k \in \CK \]
is  contractive. Now Equation (2.3) can be established in the following way: \\
We first observe that
\[ \Phi^{\CK'} Q k= \theta k, \hspace{.5cm} k \in \CK.\]
Therefore for $k_\alpha \in \CK$
\begin{eqnarray*}
\Phi^{\CK^\prime} Q \sum_{\alpha \in \tilde{\Lambda}} V_\alpha k_\alpha&=& \sum_{\alpha \in \tilde{\Lambda}}
\Phi^{\CK^\prime} V^\prime_\alpha Q k_\alpha
= \sum_{\alpha \in \tilde{\Lambda}} (L_\alpha \otimes \eins) \Phi^{\CK^\prime} Q k_\alpha \\
&=&  \sum_{\alpha \in \tilde{\Lambda}} (L_\alpha \otimes \eins) \theta k_\alpha =
M_\theta \sum_{\alpha \in \tilde{\Lambda}} (L_\alpha \otimes \eins) (e_0 \otimes k_\alpha)\\
&=& M_\theta \Phi^{\CK} \sum_{\alpha \in \tilde{\Lambda}} V_\alpha k_\alpha.
\end{eqnarray*}
As $\Phi^\CK$ and $\Phi^{\CK^\prime}$ are unitaries, we will have $M_\theta$ to be an isometry, if $Q$ is an isometry.
\end{proof}

\begin{thm} (Beurling type theorem) Let $\CU$ be a Hilbert space.
A subspace $\CM$ of $\Gamma \otimes \CU$ is invariant for $L_i
\otimes \eins, i=1, \ldots, d$ if and only if
\[\CM= M_\theta (\Gamma \otimes \CN)\]
for some Hilbert space $\CN$ and an inner function $\theta:\CN \to \Gamma \otimes \CU.$
\end{thm}

\begin{proof} The proof of the ``necessary" condition  is trivial. We prove  the ``sufficient" condition.
We consider the embedding of $\CU$ in $\Gamma \otimes \CU$ as $e_0
\otimes \CU,$ where $e_0$ is the vacuum vector in the Fock space
$\Gamma.$ We set $V_i= (L_i \otimes \eins)|_\CM, i=1,\ldots,d.$ The
row contraction $\uV=(V_1,\ldots,V_d)$ consists of isometries with
orthogonal ranges and
\[\cap^\infty_{n=0} \overline{\mbox{span}}\{V_\alpha \CM: |\alpha| =n\}=
\cap^\infty_{n=0} \overline{\mbox{span}}\{(L_\alpha \otimes \eins) (\Gamma \otimes \CU): |\alpha| =n\}=\{0\}.\]
The Wold decomposition gives us
\begin{equation}
\CM = \overline{\mbox{span}}\{V_\alpha \CN:  \alpha \in \tilde{\Lambda}\}=M(\CN)
\end{equation}
where $\CN:= \CM \ominus \overline{\mbox{span}}\{V_i \CM:
i=1,\ldots,d\}.$ We use Lemma 2.2.5 with
\[ \CR= \CM, \hspace{1cm} \CK= \CN,\]
\[ \CR^\prime=\Gamma \otimes \CU, \hspace{1cm} V^\prime_i = L_i \otimes \eins, \hspace{1cm} \CK^\prime = \CU \]
and $Q=id:\CM \to \Gamma \otimes \CU.$
Thereby an inner function $\theta: \CN \to \Gamma \otimes \CU$ is obtained with
\[\Phi^\CU Q = M_\theta \Phi^\CK. \]
The observation that $\Phi^\CU$ is identity leads us to
\[h=M_\theta \Phi^\CK h \mbox{~for~} h \in \CM.\]
Finally this implies
\[ \CM= \theta \Phi^\CK \CM= M_\theta (\Gamma \otimes \CN).\]
\end{proof}

\section{Characteristic function of Popescu}

 After the preparation done in the previous section we look at how Popescu's characteristic functions,
 which are unitary invariants for certain row contractions, are developed using Beurling type theorem. The other characteristic functions we  introduce  in later chapters bear many
 common features with the one discussed here.

Recall that for any contraction $T$ on a Hilbert space $\CH$ the mid
has the form
\[ V=
\left(\begin{array}{cc}
T & 0 \\
* & A \\
\end{array}
\right)   \]
acting on the dilation space $\tilde{\CH}= \CH \oplus \CH_A$ for some Hilbert space $\CH_A$ on which $A$ is defined.
Here $A$ is {\em $*$-stable,} i.e., $\lim_{n \to \infty} \|(A^*)^n h \| =0$ for all $h \in \CH_A.$
Infact the defect space $\CD$ (or $\CD_T$) is a wandering subspace for $\CH_A$ with respect to $A$.

At first consider the easy case, namely, of assuming $T$ to be $*$-stable.
Then the shift on $l^2 \otimes \CD_*$ is a mid of $T.$
Sch\"affer's construction (or Popescu's construction for $d=1$ case) described before gives another realisation
of the mid of $T.$  Here we observe that $\CH_A$ is embedded as an
invariant subspace for this shift. By Beurling's theorem an inner,
operator valued $H^\infty$ function $\theta$ exists
 such that
\[\CH_A = \theta (l^2 \otimes \CH')\]
for some Hilbert space $\CH'.$ (To be precise, we need to make an
identification of the shift on $l^2$ with the one on $H^2$.)
\\

A specific choice of such $\theta$ is the {\em characteristic
function} $\theta_T$ for $T.$ It is explicitly defined as

\begin{equation}
\theta_T (z): = (-T + z  D_* (\eins - z T^*)^{-1} D) |_{\CD}
\end{equation}
 where $z$ belongs to the open unit
ball of the complex plane. Here $\theta$ is bounded analytic
function on the open unit ball taking values in $B(\CD,\CD_*)$ and
\[ \| \theta_T (z)\|= 1 \mbox{~a.e. ~~~~for ~} \|z\|=1.\]
Obviously, we now have $ \CH_A = \theta_T (l^2 \otimes \CD).$
The identity (2.5) can be expanded in power series as
\[ \theta_T (z): = (-T + \sum^\infty_{n=1}z ^n D_* (T^*)^{n-1} D) |_{\CD}.\]

\begin{defn}
\begin{enumerate}
A row contraction $\uT=(T_1,\ldots, T_d)$ is called

\item $*$-stable
if $\lim_{n \to \infty} \sum_{|\alpha|=n} \|T^*_\alpha h \|^2 =0,$

\item  completely non-coisometric (c.n.c. in short) if
\[ \{ h \in \CH: \sum_{|\alpha|=n} \| T^*_\alpha h \|^2 =
\|h\|^2 \;\mbox{for all}\; n \in \N \}=\{0\}.\]

\end{enumerate}
\end{defn}

In operator matrix form Popescu's construction of mid on $\CH \oplus
(\Gamma \otimes \CD)$ is
\[ V_i= \left(\begin{array}{cc}
T_i & 0 \\
D (e_i \otimes . ) & L_i \otimes \eins \\
\end{array}
\right) \]

In case of an arbitrary row contraction $\uT,$ one way of
representing the mid is Popescu's construction. If $\uT$ is
$*$-stable, then the operator tuple $\uL \otimes \eins$ on $\Gamma
\otimes \CD$ is also a mid of $\uT$ and due to uniqueness of mids,
there exist a unitary
\[W:\CH \oplus (\Gamma \otimes \CD) \to \Gamma
\otimes \CD_* \]
 that intertwines between them. It is
immediate that this $W|_{\Gamma \otimes \CD}$ is a multi-analytic
inner operator, say $M_\theta: \Gamma \otimes \CD \to \Gamma \otimes
\CD_*.$ This is called the {\em characteristic function} of $\uT$
with symbol
\[ \theta_T:= M_{\theta}|_{e_0 \otimes \CD} \mbox{~where~} \theta_T
: \CD \to \Gamma \otimes \CD_*. \] Popescu introduced the above
defined characteristic functions in \cite{Po89b} and showed that for
c.n.c. row contractions $\uT$  they are complete invariant up to
unitary equivalence. With $P_j$ the orthogonal projection onto the
j-th component, $\theta_T$ in expanded form is
\[ \theta_T d = - e_0 \otimes \sum^d_{j=1} T_j P_j d + \sum^d_{j=1} e_j
\otimes \sum_{\alpha \in \tilde{\Lambda}} e_\alpha \otimes  D_* T^*_\alpha P_j D d, \hspace{.5cm} d \in \CD.\]

  Bhattacharyya,  Eschmeier and Sarkar defined in \cite{BES05} an analytic
characteristic function for any commuting row contraction $\uT$ as
\[ \theta_T (z): = (-T + z  D_* (\eins - z T^*)^{-1} D) |_{\CD}\]
 where $z$ belongs to the open unit ball of the $\C^d.$ These functions
 are also complete invariants for c.n.c. operator tuples.

Take $C: \CH \to \Gamma \otimes \CD_*$ to be the map
\begin{equation}
h \mapsto \sum_{\alpha \in \tilde{\Lambda}} e_\alpha \otimes D_* T_\alpha^* h,
\end{equation}
which is called the {\em Poisson kernel.} As $\CD \subset
\bigoplus^d_{i=1} \CH$ we take $P_i$ as the projection of $\CD$ onto
the $i^{th}$ component, $1 \leq i \leq d.$ From \cite{Po89b} we know
that the mid $\uV$ of a $*$-stable row contraction $\uT$ is
unitarily equivalent to $\uL \otimes \eins.$ In the following, using
an argument similar to \cite{FF90}, IX.6.4, this unitary is
constructed explicitly.

\begin{prop}
For a $*$-stable row contraction $\uT$ there exists a unitary
$$\hat{W}: \CH \oplus (\Gamma (\C^d) \otimes \CD) \to
\Gamma (\C^d) \otimes \CD_*$$ such that $(L_i \otimes \eins) \hat{W}
= \hat{W} V_i$. Together with this intertwining relation it is
determined by
\[
\hat{W} |_{\CH} = C,\quad \hat{W} |_{e_0 \otimes \CD} = \theta_T.
\]
\end{prop}
\begin{proof}
It is easy to see that when $\uT$ is $*$-stable, $C$ is an isometry
(see \cite{Po89b}) and $(L^*_i \otimes \eins) C = C T^*_i.$ Let
$\CH_*:= C \CH$ and $T_{i*}^*:=(L^*_i \otimes \eins)|_{\CH_*}.$ So
$\uL \otimes \eins$ is the mid of $\uT_*.$
Moreover $\uT_*$ is unitarily equivalent to $\uT$ as $T_{i*} C_*=C_*
T_i$ where $C_*$ is the unitary given by $C$ as a map from $\CH$ to
$\CH_*$.

Because mids are unique up to unitary
equivalence (cf. \cite{Po89a}) this implies that there is a unitary
$\hat{W}: \CH \oplus (\Gamma (\C^d) \otimes \CD) \to \Gamma (\C^d)
\otimes \CD_*$ such that $(L_i \otimes \eins) \hat{W} = \hat{W} V_i$
and $\hat{W}|_\CH=C.$ For $h_i \in \CH$ we have
\begin{eqnarray*}
&& \hat{W} \big( e_0 \otimes D(h_1,\ldots,h_d) \big) =
\hat{W}\sum_i(V_i-T_i )h_i\\
&=& \sum_i (L_i \otimes \eins) C h_i - \sum_j C T_j h_j\\
&=&\sum_i e_i \otimes \sum_{\alpha \in \tilde{\Lambda}} e_\alpha \otimes D_* T_\alpha^* h_i -\sum_j \sum_{\beta \in \tilde{\Lambda}} e_\beta \otimes D_*T_\beta^*T_jh_j\\
&=& -e_0 \otimes \sum_i D_* T_i h_i + \sum_j e_j \otimes \sum_{\alpha \in \tilde{\Lambda}} e_\alpha \otimes D_* T_\alpha^*P_j D^2(h_1,\ldots,h_d)\\
&=& -e_0 \otimes \sum_i T_i P_i D(h_1,\ldots,h_d)\\
& +& \sum_j e_j \otimes \sum_{\alpha \in \tilde{\Lambda}} e_\alpha
\otimes  D_* T_\alpha^*P_j D(D(h_1,\ldots,h_d))\\
&=& \theta_T D(h_1,\ldots,h_d).
\end{eqnarray*}
Hence $\hat{W} |_{e_0 \otimes \CD} = \theta_T$. From the explicit
form of Popescu's dilation we see that the restriction of $V_i$ to
$\Gamma (\C^d) \otimes \CD$ coincides with $L_i \otimes \eins$. This
shows that $\hat{W}$ is determined by $C$ and $\theta_T$ (together
with $(L_i \otimes \eins) \hat{W} = \hat{W} V_i$).
\end{proof}

It is observed from this that $\theta_T$ is an isometry when $\uT$
is $*$-stable (also realized in Remark 3.2 of \cite{Po89b}).
\\

The noncommutative Beurling type theorem (Theorem 2.2.6) and the
characteristic functions of row contractions indicate that there
should be a parallel theory of noncommutative complex analysis. Such
a theory was given by Popescu (\cite{Po91},\cite{Po06}). There one
replaces open unit disc of the complex plane by
$[B(\CH)^d]_1:=\{(X_1,\ldots,X_d) \in B(\CH)^d: \|X_1X^*_1 + \ldots
+ X_d X^*_d\| < 1\}.$ Let $Z_1,\ldots, Z_d$ denote $d$ noncommuting
variables and
\[Hol (\D^N_d):= \{ F=\sum_{\alpha \in \tilde{\Lambda}} a_k Z_\alpha: \limsup_{k \rightarrow \infty} (\sum_{|\alpha|=k}
|a_\alpha|^2)^{\frac{1}{2k}} \leq 1 \}.\]
Then the noncommutative analogue of $H^\infty$ is given by
\[H^\infty (\D^N_d)=\{F \in Hol (\D^N_d) : \hspace{1cm} \|F\|_\infty := \sup \|F(X_1,\ldots,X_d)\| < \infty \]
\[\mbox{~where supremum is taken over $(X_1,\ldots,X_d)$ in ~} [B(\CH)^d]_1\},\]
Note that $\D^N_d$ alone has no explicit meaning  and the notation
$Hol (\D^N_d)$ is used to denote the set which is the counterpart of
the set of classical analytic functions on the open unit ball of the
complex plane.

\section{Hilbert C*-module}

Next we review some selected topics about Hilbert $C^*$-modules
which we need in Section 6.3.

For a given $C^*$-algebra $\CA,$ a pre-Hilbert $C^*$-module $\CG$ on
$\CA$ is a right $\CA$-module with a sesquilinear inner product:

\[
<.,.> : \CG \times \CG \to \CA,
\]
 for $\chi,\eta \in \CG,a \in \CA$ such that:
\begin{itemize}
\item[(a)] $<\chi, \eta a> =<\chi,\eta> a,$
\item[(b)] $<\chi,\chi> \geq 0$ and
\item[(c)] $<\chi,\chi>=0 \Leftrightarrow \chi =0$
\end{itemize}
A norm can be defined on $\CG$ by $\|.\|:= \|
<.,.>\|^{\frac{1}{2}}_{\CA}$ where $\|.\|_{\CA}$ is the
$C^*$-algebra norm of $\CA.$ Such a module $\CG$ is called {\em
Hilbert $C^*$-module} on $\CA$ if it is complete with respect to the
above norm. By $\CL(\CG)$ we denote the set of all adjointable maps
on $\CG.$ The set  $\CL(\CG)$ turns out to be a $C^*$-algebra with
respect to the operator norm induced by $(\CG,\|\hspace{2mm} \|).$

Let $\varphi: \CA \to \CL(\CG)$ be a $*$-homomorphism and  write
$\varphi (a) \eta :=a \eta$ for all $a \in \CA, \eta \in \CE.$ The
Hilbert $C^*$-module $\CG$ together with such a $*$-homomorphism is
called a {\em Hilbert bimodule} on $\CA.$ There is an associated
tensor product structure $\CG \otimes \CG$ for any Hilbert bimodule
$\CG$ given by
\[  <\eta_1 \otimes \chi_1, \eta_2 \otimes \chi_2>=<\chi_1,<\eta_1, \eta_2> \chi_2>
\mbox{~~for~} \eta_1,\eta_2,\chi_1,\chi_2 \in \CG.\]

\begin{defn} A {\em $W^*$-correspondence} $\CE$ over a von Neumann algebra
$\CM$ is a self-dual Hilbert $\CM$-bimodule where the corresponding
left action $\varphi$ of $\CM$ on $\CE$ is normal.
\end{defn}

\begin{defn}
A pair $(T,\sigma)$ is called a covariant representation of
$W^*$-correspondence $\CE$ over $\CM$ on a Hilbert space $\CH,$ if:
\begin{itemize}
\item[(a)] $T: \CE \to B(\CH)$ is a linear map that is continuous w.r.t. the $\sigma$-topology of \cite{BDH88} on $\CE$
and the ultraweak topology on $B(\CH).$
\item[(b)] $\sigma: \CM \to B(\CH)$ is a normal homomorphism.
\item[(c)] $T(a \xi) = \sigma (a) T(\xi), T(\xi a)=T(\xi) \sigma (a)$ for all $\xi \in \CE, a \in \CM.$
\end{itemize}
If a completely contractive covariant representation satisfies
\[ T(\xi)^* T(\eta) = \sigma (\langle \xi,  \eta\rangle)\]
in addition,  it is called isometric.
\end{defn}

One of our prime objective is to use dilation theory of covariant representations
of $W^*$-correspondences to extend our theory of characteristic functions and then apply it to classify
these covariant representations. First we remark that for any $W^*$-correspondence $\CE$ over a von Neumann algebra $\CM$
and a normal representation $\sigma: \CM \to B(\CH)$ there is an induced
tensor product $\CE \otimes_{\sigma} \CH,$ which is a Hilbert space,
with the defining identities:
\[\xi_1 a \otimes \eta_1 = \xi_1 \otimes \sigma (a) \eta_1 \mbox{~~and} \]
\[ \langle \xi_1 \otimes \eta_1, \xi_2 \otimes \eta_2 \rangle = \langle \eta_1, \sigma
(\langle \xi_1, \xi_2\rangle) \eta_2 \rangle\] for $\xi_1,\xi_2 \in
\CE$ and $\eta_1,\eta_2 \in \CH.$ In addition if $\sigma$ is also faithful, define
\[\CE^\sigma := \{\mu \in B(\CH, \CE \otimes_{\sigma} \CH): \mu \sigma (a) = (\varphi(a) \otimes
\eins) \mu~~~~ \forall a \in \CM \}.\] $\CE^\sigma$ is called the
{\em $\sigma$-dual} of $\CE.$ Its open unit ball is denoted by
$\D(\CE^\sigma).$
\\

For every covariant representation $(T, \sigma)$ on $\CH,$ we can
get an operator $\tilde{T}$ in $B(\CE \otimes \CH, \CH)$ given by
\[ \tilde{T} (\eta \otimes h): = T (\eta) h.\]
It is handy to work with the operator $\tilde{T}$ instead of $(T,\sigma)$ in dilation theory.
The next lemma provides a dictionary to translate some important notions for $(T,\sigma)$ in terms of that for
$\tilde{T}$ and vice versa.

\begin{lem} \cite{MS98}
\begin{itemize}
\item[(a)] $T$ is completely bounded if and only if $\tilde{T}$ is bounded and
\[ \| T\|_{cb} = \|\tilde{T} \| \]

\item[(b)]
$T$ is completely contractive $\Leftrightarrow \|\tilde{T}\| \leq 1.$
\item[(c)] $(T,\sigma)$ is isometric if and only if $\tilde{T}$ is isometric.

\end{itemize}
\end{lem}

\begin{proof} We prove only  the part (b) here. Start with the assumption that
$\|\tilde{T}\| \leq 1.$ Let $(\eta_{ij})_{k \times k}$  be an
element in $M_k (\CE).$ For $h=(h_1,\ldots,h_k),h_i \in \CH$ we have
\begin{eqnarray*}
\|(T(\eta_{ij}))_{k \times k} h\|^2 &=& \sum^k_{i=1} \|\sum^k_{j=1}  T(\eta_{ij}) h_j\|^2 = \sum^k_{i=1} \|\tilde{T}(\sum^k_{j=1} \eta_{ij} \otimes h_j)\|^2\\
&\leq& \sum^k_{i=1} \|(\sum^k_{j=1} \eta_{ij} \otimes h_j)\|^2 = \sum_{i,j,l} \langle \eta_{ij} \otimes h_j, \eta_{il} \otimes h_l \rangle\\
&=& \sum_{i,j,l} \langle h_j, \sigma (\langle \eta_{ij}, \eta_{il}\rangle) h_l\rangle .
\end{eqnarray*}
Therefore
\[\|(T(\eta_{ij}))_{k \times k} h\|^2 \leq \| (\sum^k_{i=1} \sigma (\langle \eta_{ij}, \eta_{il}\rangle))_{k \times k} \| \|h \|^2. \]
 $\sigma$ is completely contractive because it is a $*$-representation. This implies
\[ \|(T(\eta_{ij}))_{k \times k} h\|^2 \leq \|(\eta_{ij})_{k \times k} \|^2 \|h\|^2,\]
Hence $\|T\|_{cb} \leq 1.$

For the converse first we claim that: If $T$ is completely contractive, then for every $\zeta_1,\zeta_2,\ldots, \zeta_k \in \CE$
\begin{equation}
(T(\zeta_i)^*T(\zeta_j))_{k \times k} \leq (\sigma (\langle \zeta_i, \zeta_j \rangle))_{k \times k}.
\end{equation}

Using the claim, the converse of part (b) can be proved in the following way:\\
Let $\zeta_1, \ldots, \zeta_k \in \CE$ and $h=(h_1,\ldots,h_k), h_i
\in \CH.$  Assuming equation (2.7) holds, we observe that
\[\|\sum^k_{i=1} \zeta_i \otimes h_i  \|^2_{\CE \otimes \CH} = \sum_{i,j} \langle h_i, \sigma (\langle \zeta_i,\zeta_j\rangle) h_j \rangle
= \langle h, (\sigma (\langle \zeta_i,\zeta_j\rangle))_{k \times k}
h \rangle \] and on the other hand
\begin{eqnarray*}
\|\tilde{T}(\sum^k_{i=1} \zeta_i \otimes h_i)  \|^2 &=& \| \sum^k_{i=1} T(\zeta_i) h_i \|^2\\
&=& \sum_{i,j} \langle h_i, T ( \zeta_i)^* T (\zeta_j) h_j \rangle\\
&=& \langle h, (T ( \zeta_i)^* T (\zeta_j))_{k \times k} h \rangle .
\end{eqnarray*}
From the above two calculations we conclude that $\|\tilde{T}\| \leq
1.$

Proof of the Claim: We start with some $\zeta_1,\ldots, \zeta_k \in
\CE.$ Because we will be studying how $T$ operates on
$\zeta_i,i=1,\ldots,d,$ we assume without loss of generality that
$\CE$ is generated by the $\zeta_i$'s. A corollary of Kasparov's
stabilization theorem (cf. \cite{La95} or Appendix of chapter 6)
allows us to pick vectors $\{\eta_i\}^\infty_{i=1}$ in $\CE$ such
that $\{ \eta_i \otimes \eta^*_i\}^\infty_{i=1}$ is a contractive
approximate identity of the algebra $K(\CE)$ of compacts.  The
projection $Q:=(\langle \eta_i, \eta_j \rangle)_{\infty \times
\infty}$ is an element of the multiplier algebra of $\CM \otimes
\mathfrak{K}$ (where $\mathfrak{K}$ stands for the $C^*$-algebra of
compact operators on any separable infinite dimensional Hilbert
space). Denote $(\sigma (\langle \eta_i, \eta_j \rangle))_{\infty
\times \infty}$ (in $B(\CH^{(\infty)})$) and $(\sigma (\langle
\zeta_i, \eta_j \rangle))_{\infty \times k}$ (in
$B(\CH^{(\infty)},\CH^{(n)}))$ by $\tilde{\sigma}(Q)$ and $R$
respectively. We have
\[ R \tilde{\sigma}(Q)R^* =  (\sigma (\langle \zeta_i, \zeta_j \rangle))_{k \times k}, \hspace{.5cm} R(T(\eta_i)^* T(\eta_j))_{\infty \times \infty} R^*
=(T(\zeta_i)^* T(\zeta_j))_{k \times k}. \]
Therefore we just need to show
\[ (T(\eta_i)^* T(\eta_j))_{\infty \times \infty} \leq \tilde{\sigma}(Q)\]
and this is an easy exercise on observing that
\[ (T(\eta_1),T(\eta_2),\ldots) \tilde{\sigma}(Q)=  (T(\eta_1),T(\eta_2),\ldots)\]
and employing the fact that $T$ is completely contractive.
\end{proof}

\chapter
 {Characteristic Functions for Ergodic Tuples}

\vspace{.5cm}

{\bf Abstract:} {\em Motivated by a result on weak Markov dilations,
we define a notion of characteristic function for ergodic and
coisometric row contractions with a one-dimensional invariant
subspace for the adjoints. This extends a definition given by
G.\,Popescu. We prove that our characteristic function is a complete
unitary invariant for such tuples and show how it can be computed. }

\vspace{4cm}

Joint work with Rolf Gohm. Published in Integral Equations and
Operator Theory 58 (2007), 43-63.

\newpage

\setcounter{section}{-1}
\section{Introduction}

If $Z = \sum^d_{i=1} A_i \cdot A^*_i$ is a normal, unital, ergodic,
completely positive map on $B(\CH)$, the bounded linear operators on
a complex separable Hilbert space, and if there is a (necessarily
unique) invariant vector state for $Z$, then we also say that $\uA =
(A_1,\ldots,A_d)$ is a coisometric, ergodic row contraction with a
one-dimensional invariant subspace for the adjoints. Precise
definitions are given below. This is the main setting to be
investigated in this paper.

In Section 3.1 we give a concise review of a result  on the
dilations of $Z$ obtained by R.\,Gohm in \cite{Go04} in a chapter
called `Cocycles and Coboundaries'. There exists a conjugacy between
a homomorphic dilation of $Z$ and a tensor shift, and we emphasize
an explicit infinite product formula that can be obtained for the
intertwining unitary. \cite{Go04} may also be consulted for
connections of this topic to a scattering theory for noncommutative
Markov chains by B.\,K\"ummerer and H.\,Maassen (cf. \cite{KM00})
and more general for the relevance of this setting in applications.

In this work we are concerned with its relevance in operator theory
and correspondingly in Section 3.2 we shift our attention to the row
contraction $\uA = (A_1,\ldots,A_d)$. Our starting point has been
the observation that the intertwining unitary mentioned above has
many similarities with the notion of characteristic function
occurring in the theory of functional models of contractions, as
initiated by B.\,Sz.-Nagy and C.\,Foias (cf. \cite{NF70,FF90}). In
fact, the center of our work is the commuting diagram 3.5 in Section
3.3, which connects the results in \cite{Go04} mentioned above with
the theory of minimal isometric dilations of row contractions by
G.\,Popescu (cf. \cite{Po89a}) and shows that the intertwining
unitary determines a multi-analytic inner function, in the sense
introduced by G.\,Popescu in \cite{Po89c,Po95}. We call this inner
function the {\it extended characteristic function} of the tuple
$\uA$, see Definition 3.3.3.

Section 3.4 is concerned with an explicit computation of this inner
function. In Section 3.5 we show that it is an extension of the
characteristic function of the $*$-stable part $\uAo$ of $\uA$, the
latter in the sense of Popescu's generalization of the
Sz.-Nagy-Foias theory to row contractions (cf. \cite{Po89b}). This
explains why we call our inner function an {\it extended}
characteristic function. The row contraction $\uA$ is a
one-dimensional extension of the $*$-stable row contraction $\uAo$,
and in our analysis we separate the new part of the characteristic
function from the part already given by Popescu.

G.\,Popescu has shown in \cite{Po89b} that for completely
non-coisometric tuples, in particular for $*$-stable ones, his
characteristic function is a complete invariant for unitary
equivalence. In Section 3.6 we prove that our extended
characteristic function does the same for the tuples $\uA$ described
above. In this sense it is {\it characteristic}. This is remarkable
because the strength of Popescu's definition lies in the completely
non-coisometric situation while we always deal with a coisometric
tuple $\uA$. The extended characteristic function also does not
depend on the choice of the decomposition $\sum^d_{i=1} A_i \cdot
A^*_i$ of the completely positive map $Z$ and hence also
characterizes $Z$ up to conjugacy. We think that together with its
nice properties established earlier this clearly indicates that the
extended characteristic function is a valuable tool for classifying
and investigating such tuples respectively such completely positive
maps.

Section 3.7 contains a worked example for the constructions in this
paper.

\section{Weak Markov dilations and conjugacy}

In this section we give a brief and condensed review of results in
\cite{Go04}, Chapter 2, which will be used in the following and
which, as described in the introduction, motivated the
investigations documented in this paper. We also introduce notation.

A theory of {\em weak Markov dilations} has been developed in
\cite{BP94}. For a (single) normal unital completely positive map
$Z: B(\CH) \rightarrow B(\CH)$, where $B(\CH)$ consists of the
bounded linear operators on a (complex, separable) Hilbert space, it
asks for a normal unital $^*-$endomorphism $\hat{J}: B(\hat{\CH})
\rightarrow B(\hat{\CH})$, where $\hat{\CH}$ is a Hilbert space
containing $\CH$, such that for all $n\in\N$ and all $x\in B(\CH)$
\[
Z^n(x) = p_{\CH}\, \hat{J}^n(x p_{\CH})\; |_{\CH}.
\]
Here $p_{\CH}$ is the orthogonal projection onto $\CH$. There are
many ways to construct $\hat{J}$. In \cite{Go04}, 3.2.3, we gave a
construction analogous to the idea of `coupling to a shift' used in
\cite{Ku85} for describing quantum Markov processes. This gives rise
to a number of interesting problems which remain hidden in other
constructions.
\\

We proceed in two steps. First note that there is a Kraus
decomposition $Z(x) = \sum^d_{i=1} a_i\, x\, a^*_i$ with
$(a_i)^d_{i=1} \subset B(\CH)$. Here $d=\infty$ is allowed in which
case the sum should be interpreted as a limit in the strong operator
topology. Let $\CP$ be a $d$-dimensional Hilbert space with
orthonormal basis $\{\epsilon_1,\ldots,\epsilon_d\}$, further $\CK$
another Hilbert space with a distinguished unit vector $\OK \in
\CK$. We identify $\CH$ with $\CH \otimes \OK \subset \CH \otimes
\CK$ and again denote by $p_{\CH}$ the orthogonal projection onto
$\CH$. For $\CK$ large enough there exists an isometry
\[
     u: \CH \otimes \CP \rightarrow \CH \otimes \CK
\quad \mbox{s.t.}\quad p_{\CH}\, u (h \otimes \epsilon_i) = a_i(h),
\]
for all $h\in\CH,\; i=1,\ldots,d$, or equivalently,
\[
u^* (h \otimes \OK) = \sum^d_{i=1} a^*_i(h) \otimes \epsilon_i.
\]
Explicitly, one may take $\CK = \C^{d+1}$ (resp.
infinite-dimensional) and identify
\[
\CH \otimes \CK \;\simeq\; (\CH \otimes \OK) \oplus \bigoplus^d_1
\CH \;\simeq\; \CH \oplus \bigoplus^d_1 \CH.
\]

Then, using isometries $u_1,\ldots,u_d: \CH \rightarrow \CH \oplus
\bigoplus^d_1 \CH$ with orthogonal ranges and such that $a_i =
p_{\CH} u_i$ for all $i$ (for example, such isometries are
explicitly constructed in Popescu's formula for isometric dilations,
cf. \cite{Po89a} or equation 3.4 in Section 3), we can define
\[
u (h \otimes \epsilon_i) := u_i(h)
\]
for all $h\in\CH,\; i=1,\ldots,d$ and check that $u$ has the desired
properties. Now we define a $^*-$homomorphism
\begin{eqnarray*}
J: B(\CH) &\rightarrow& B(\CH \otimes \CK), \\
     x &\mapsto& u \, (x\otimes \eins_{\CP})\, u^*.
\end{eqnarray*}
It satisfies
\[
 p_{\CH}\, J(x) (h \otimes \OK) = p_{\CH}\,u \, (x\otimes \eins) u^*
(h \otimes \OK)
\]
\[
= p_{\CH}\,u \, (x\otimes \eins) \big( \sum^d_{i=1} a^*_i(h) \otimes
\epsilon_i \big) = \sum^d_{i=1} a_i\, x\, a^*_i (h) = Z(x)(h),
\]
which means that $J$ is a kind of first order dilation for $Z$.
\\

For the second step we write $\tilde{\CK} := \bigotimes^\infty_1
\CK$ for an infinite tensor product of Hilbert spaces along the
sequence $(\OK)$ of unit vectors in the copies of $\CK$. We have a
distinguished unit vector $\Omega_{\tilde{\CK}}$ and a (kind of)
tensor shift
\[
R: B(\tilde{\CK}) \rightarrow B(\CP \otimes \tilde{\CK}), \quad
\tilde{y} \mapsto \eins_{\CP} \otimes \tilde{y}.
\]
Finally $\tilde{\CH} := \CH \otimes \tilde{\CK}$ and we define a
normal $^*-$endomorphism
\begin{eqnarray*}
\tilde{J}: B(\tilde{\CH}) &\rightarrow& B(\tilde{\CH}), \\
B(\CH) \otimes B(\tilde{\CK}) \ni x \otimes \tilde{y} &\mapsto& J(x)
\otimes \tilde{y} \in B(\CH \otimes \CK) \otimes B(\tilde{\CK}).
\end{eqnarray*}
Here we used von Neumann tensor products and (on the right hand
side) a shift identification $\CK \otimes \tilde{\CK} \simeq
\tilde{\CK}$. We can also write $\tilde{J}$ in the form
\[
\tilde{J}(\cdot) = u\, (Id_{\CH} \otimes R)(\cdot) \, u^*,
\]
where $u$ is identified with $u \otimes \eins_{\tilde{\CK}}$. The
natural embedding $\CH \simeq \CH \otimes \Omega_{\tilde{\CK}}
\subset \tilde{\CH}$ leads to the restriction $\hat{J} := \tilde{J}
|_{\hat{\CH}}$ with $\hat{\CH} := \overline{\mbox{span}}_{n\geq 0}
\tilde{J}^n(p_{\CH})(\tilde{\CH})$, which can be checked to be a
normal unital $^*$-endomorphism satisfying all the properties of a
weak Markov dilation for $Z$ described above. See \cite{Go04}, 2.3.

A Kraus decomposition of $\hat{J}$ can be written as
\[
\hat{J}(x) = \sum^d_{i=1} t_i\, x\, t^*_i,
\]
where $t_i \in B(\hat{\CH})$ is obtained by linear extension of $\CH
\otimes \tilde{\CK} \ni h \otimes \tilde{k} \mapsto u_i(h) \otimes
\tilde{k} = u (h \otimes \epsilon_i) \otimes \tilde{k} \in (\CH
\otimes \CK)\otimes \tilde{\CK} \simeq \CH \otimes \tilde{\CK}$.
Because $\hat{J}$ is a normal unital $^*-$endomorphism the
$(t_i)^d_{i=1}$ generate a representation of the Cuntz algebra $\Od$
on $\hat{\CH}$ which we called a {\em coupling representation} in
\cite{Go04}, 2.4. Note that the tuple $(t_1,\ldots,t_d)$ is an
isometric dilation of the tuple $(a_1,\ldots,a_d)$, i.e., the $t_i$
are isometries with orthogonal ranges and $p_{\CH} t^n_i |_{\CH} =
a^n_i$ for all $i=1,\ldots,d$ and $n\in\N$.
\\

The following {\em multi-index notation} will be used frequently in
this work. Let $\Lambda $ denote the set $\{ 1, 2, \ldots , d\}.$
For  operator tuples $(a_1, \ldots , a_d),$ given $\alpha =(\alpha
_1, \ldots , \alpha _m)$ in $\Lambda ^m$,\, $a_{\alpha }$ will stand
for the operator $a_{\alpha _1}a_{\alpha _2}\ldots a_{\alpha _m}$,
\; $|\alpha| := m$. Further $\tilde{\Lambda}:=\cup_{n=0}^{\infty}
\Lambda^n$, where $\Lambda ^0:=\{ 0\}$ and $a_0$ is the identity
operator. If we write $a^*_{\alpha}$ this always means
$(a_{\alpha})^* = a^*_{\alpha _m}\ldots a^*_{\alpha _1}$.
\\

Back to our isometric dilation, it can be checked that
\[
\overline{\mbox{span}}\{t_{\alpha} h: h \in \CH, \alpha \in
\tilde{\Lambda} \} =\hat{\CH},
\]
which means that we have a {\em minimal isometric dilation}, cf.
\cite{Po89a} or the beginning of Section 3.3. For more details on
the construction above see \cite{Go04}, 2.3 and 2.4.

Assume now that there is an invariant vector state for $Z: B(\CH)
\rightarrow B(\CH)$ given by a unit vector $\OH \in \CH$.
Equivalent: There is a unit vector $\OP = \sum^d_{i=1}
\overline{\omega}_i \epsilon_i \in \CP$ such that $u (\OH\otimes\OP)
= \OH \otimes \OK$. Also equivalent: For $i=1,\ldots,d$ we have
$a^*_i\, \OH = \overline{\omega}_i\, \OH$. Here $\omega_i \in \C$
with $\sum^d_{i=1} |\omega_i|^2 = 1$ and we used complex conjugation
to get nice formulas later. See \cite{Go04}, A.5.1, for a proof of
the equivalences.

On $\tilde{\CP} := \bigotimes^\infty_1 \CP$ along the unit vectors
$(\OP)$ in the copies of $\CP$ we have a tensor shift
\[
S: B(\tilde{\CP}) \rightarrow B(\tilde{\CP}), \quad \tilde{y}
\mapsto \eins_{\CP} \otimes \tilde{y}.
\]
Its Kraus decomposition is $S(\tilde{y}) = \sum^d_{i=1}
s_i\,\tilde{y}\,s^*_i$ with $s_i \in B(\tilde{\CP})$ and
$s_i(\tilde{k}) = \epsilon_i \otimes \tilde{k}$ for
$\tilde{k}\in\tilde{\CP}$ and $i=1,\ldots,d$. In \cite{Go04}, 2.5,
we obtained an interesting description of the situation when the
dilation $\hat{J}$ is conjugate to the shift endomorphism $S$. This
result will be further analyzed in this paper. We give a version
suitable for our present needs but the reader should have no
problems to obtain a proof of the following from \cite{Go04}, 2.5.

\begin{thm}
Let $Z: B(\CH) \rightarrow B(\CH)$ be a normal unital completely
positive map with an invariant vector state $\langle \OH, \cdot\,
\OH \rangle$. Notation as introduced above, $d\geq 2$. The following
assertions are equivalent:
\begin{itemize}
\item[(a)]
$Z$ is ergodic, i.e., the fixed point space of $Z$ consists of
multiples of the identity.
\item[(b)]
The vector state $\langle \OH, \cdot\, \OH \rangle$ is absorbing for
$Z$, i.e., if $n\to\infty$ then $\phi(Z^n(x)) \to \langle \OH, x \OH
\rangle\;$ for all normal states $\phi$ and all $x\in B(\CH)$. (In
particular, the invariant vector state is unique.)
\item[(c)]
$\hat{J}$ and $S$ are conjugate, i.e., there exists a unitary ${\bf
w}: \hat{\CH} \rightarrow \tilde{\CP}$ such that
\[
\hat{J}(\hat{x}) = {\bf w}^*\, S({\bf w}\, \hat{x}\, {\bf w}^*)\,
{\bf w}.
\]
\item[(d)]
The $\Od-$representations corresponding to $\hat{J}$ and $S$ are
unitarily equi\-valent, i.e.,
\[
{\bf w}\, t_i = s_i\, {\bf w} \quad \mbox{for}\; i=1,\ldots,d.
\]
\end{itemize}
An explicit formula can be given for an intertwining unitary as
occurring in (c) and (d). If any of the assertions above is valid
then the
 following limit exists strongly,
\[
\tilde{{\bf w}} = \lim_{n\to\infty} u^*_{0n} \ldots u^*_{01}: \; \CH
\otimes \tilde{\CK} \rightarrow \CH \otimes \tilde{\CP},
\]
where we used a leg notation, i.e., $u_{0n} = (Id_{\CH} \otimes
R)^{n-1}(u)$. In other words $u_{0n}$ is $u$ acting on $\CH$ and on
the $n-$th copy of $\CP$. Further $\tilde{{\bf w}}$ is a partial
isometry with initial space $\hat{\CH}$ and final space $\tilde{\CP}
\simeq \OH \otimes \tilde{\CP} \subset \CH \otimes \tilde{\CP}$ and
we can define ${\bf w}$ as the corresponding restriction of
$\tilde{{\bf w}}$.
\end{thm}

To illustrate the product formula for ${\bf w}$, which will be our
main interest in this work, we use it to derive (d).
\[
{\bf w}\,t_i (h \otimes \tilde{k}) = {\bf w} \, \big[ u (h \otimes
\epsilon_i) \otimes \tilde{k} \big] = \lim_{n\to\infty} u^*_{0n}
\ldots u^*_{01} u_{01} (h \otimes \epsilon_i \otimes \tilde{k})
\]
\[
= \lim_{n\to\infty} u^*_{0n} \ldots u^*_{02} (h \otimes \epsilon_i
\otimes \tilde{k}) = s_i\, {\bf w} (h \otimes \tilde{k}).
\]
Let us finally note that Theorem 3.1.1 is related to the conjugacy
results in \cite{Pow88} and \cite{BJP96}. Compare also Proposition
3.2.4.

\section{Ergodic coisometric row contractions}

In the previous section we considered a map $Z: B(\CH) \rightarrow
B(\CH)$ given by $Z(x) = \sum^d_{i=1} A_i\, x\, A^*_i$, where
$(A_i)^d_{i=1} \subset B(\CH)$. We can think of $(A_i)^d_{i=1}$ as a
$d$-tuple $\underline{A} = (A_1,\ldots,A_d)$ or (with the same
notation) as a linear map
\[
\underline{A} = (A_1,\ldots,A_d):\; \bigoplus^d_{i=1} \CH
\rightarrow \CH.
\]
(Concentrating now on the tuple we have changed to capital letters
$A$. We will sometimes return to lower case letters $a$ when we want
to emphasize that we are in the (tensor product) setting of Section
3.1.) We have the following dictionary.

\begin{eqnarray*}
Z(\eins) \leq \eins &\Leftrightarrow&
\sum^d_{i=1} A_i\,A^*_i \leq \eins \\
&\Leftrightarrow&
\underline{A} \;\mbox{is a contraction} \\
\\
Z(\eins) = \eins & \Leftrightarrow&
\sum^d_{i=1} A_i\,A^*_i = \eins \\
\big( Z \;\mbox{is called unital}\big) & &
\big( \underline{A} \;\mbox{is called coisometric} \big) \\
\\
\langle \OH, \cdot \OH \rangle = \langle \OH, Z(\cdot) \OH \rangle
&\Leftrightarrow& A^*_i\, \OH = \overline{\omega}_i\, \OH, \;
\omega_i \in \C, \; \sum^d_{i=1} |\omega_i|^2 = 1 \\
\big(\mbox{ invariant vector state} \big) & &
\big(\mbox{ common eigenvector for adjoints} \big) \\
\\
Z \;\mbox{ergodic} &\Rightarrow&
\{A_i, A^*_i\}^\prime = \C\,\eins \\
\big(\mbox{trivial fixed point space} \big) & & \big(\;\mbox{trivial
commutant}\; \big)
\end{eqnarray*}

The converse of the implication at the end of the dictionary is not
valid. This is related to the fact that the fixed point space of a
completely positive map is not always an algebra. Compare the
detailed discussion of this phenomenon in \cite{BJKW00}.

By a slight abuse of language we call the tuple (or row contraction)
$\underline{A} = (A_1,\ldots,A_d)$ {\em ergodic} if the
corresponding map $Z$ is ergodic. With this terminology we can
interpret Theorem 3.1.1 as a result about ergodic coisometric row
contractions $\underline{A}$ with a common eigenvector $\OH$ for the
adjoints $A^*_i$. This will be examined starting with Section 3.3.
To represent these objects more explicitly let us write $\Ho := \CH
\ominus \C\,\OH$. With respect to the decomposition $\CH = \C\,\OH
\oplus \Ho$ we get $2 \times 2 -$ block matrices
\begin{equation}
A_i = \left( \begin{array}{cc}
        \omega_i & 0 \\
        |\ell_i\rangle & \Ao_i \\
        \end{array}
\right), \quad A^*_i = \left( \begin{array}{cc}
        \overline{\omega}_i & \langle \ell_i| \\
        0 & \Aso_i \\
        \end{array}
\right).
\end{equation}

Here $\Ao_i \in B(\Ho)$ and $\ell_i\in\Ho$. For the off-diagonal
terms we used a Dirac notation that should be clear without further
comments.

Note that the case $d=1$ is rather uninteresting in this setting
because if $A$ is a coisometry with block matrix $ \left(
\begin{array}{cc}
        \omega & 0 \\
        |\ell\rangle & \Ao \\
        \end{array}
\right) $ then because
\[
\left( \begin{array}{cc}
        1 & 0 \\
        0 & \eins \\
        \end{array}
\right) = A\,A^* =\left( \begin{array}{cc}
        |\omega|^2 & \omega\, \langle\ell| \\
        \overline{\omega}\, |\ell\rangle & |\ell\rangle\langle\ell|
          + \Ao\,\Aso\\
        \end{array}
\right)
\]
we always have $\ell = 0$. But for $d \geq 2$ there are many
interesting examples arising from unital ergodic completely positive
maps with invariant vector states. See Section 3.1 and also Section
3.7 for an explicit example. We always assume $d \geq 2$.

\begin{prop}
A coisometric row contraction $\underline{A} = (A_1,\ldots,A_d)$ is
ergodic with common eigenvector $\OH$ for the adjoints
$A^*_1,\ldots,A^*_d$ if and only if $\Ho$ is invariant for
$A_1,\ldots,A_d$ and the restricted row contraction
$(\Ao_1,\ldots,\Ao_d)$ on $\Ho$ is $*$-stable, i.e., for all
$h\in\Ho$
\[
\lim_{n\to\infty} \sum_{|\alpha|=n} \| \Aso_\alpha h\|^2 = 0\;.
\]
\end{prop}

Here we used the multi-index notation introduced in Section 3.1.
Note that $*$-stable tuples are also called pure, we prefer the
terminology from \cite{FF90}.
\\

\begin{proof}
It is clear that $\OH$ is a common eigenvector for the adjoints if
and only if $\Ho$ is invariant for $A_1,\ldots,A_d$. Let $Z(\cdot) =
\sum^d_{i=1} A_i \cdot A^*_i$ be the associated completely positive
map. With $q \, := \eins - | \OH \rangle \langle \OH |$, the
orthogonal projection onto $\Ho$, and by using $ q\, A_i\,q = A_i\,
q \simeq\; \Ao_i$ for all $i$, we get
\[
Z^n(q) = \sum_{|\alpha|=n} A_\alpha\, q\, A^*_\alpha =
\sum_{|\alpha|=n} \Ao_\alpha\, \Aso_\alpha
\]
and thus for all $h\in\Ho$
\[
\sum_{|\alpha|=n} \| \Aso_\alpha h\|^2 = \langle h, Z^n(q)\, h
\rangle.
\]
Now it is well known that ergodicity of $Z$ is equivalent to $Z^n(q)
\to 0$ for $n\to\infty$ in the weak operator topology. See
\cite{GKL06}, Prop. 3.3.2. This completes the proof.
\end{proof}

\begin{rem}
Given a coisometric row contraction $\underline{a} =
(a_1,\ldots,a_d)$ we also have the isometry $u: \CH \otimes \CP
\rightarrow \CH \otimes \CK$ from Section 3.1. We introduce the
linear map $a: \CP \rightarrow B(\CH),\; k \mapsto a_k$ defined by
\[
a^*_k(h) \otimes k := (\eins_{\CH} \otimes | k \rangle \langle k
|)\, u^*(h \otimes \OK).
\]
Compare \cite{Go04}, A.3.3. In particular $a_i = a_{\epsilon_i}$ for
$i=1,\ldots,d$, where $\{\epsilon_1,\ldots,\epsilon_d\}$ is the
orthonormal basis of $\CP$ used in the definition of $u$. Arveson's
metric operator spaces, cf. \cite{Ar03}, give a conceptual
foundation for basis transformations in the operator space linearly
spanned by the $a_i$. Similarly, in our formalism a unitary in
$B(\CP)$ transforms $\underline{a} = (a_1,\ldots,a_d)$ into another
tuple $\underline{a}^\prime = (a^\prime_1,\ldots,a^\prime_d)$. If
$\OH$ is a common eigenvector for the adjoints $a^*_i$ then $\OH$ is
also a common eigenvector for the adjoints $(a^\prime_i)^*$ but of
course the eigenvalues are transformed to another tuple
$\uomega^\prime = (\omega^\prime_1,\ldots,\omega^\prime_d)$. We
should consider the tuples $\ua$ and $\underline{a}^\prime$ to be
essentially the same. This also means that the complex numbers
$\omega_i$ are not particularly important and they should not play a
role in classification. They just reflect a certain choice of
orthonormal basis in the relevant metric operator space. Independent
of basis transformations is the vector $\OP = \sum^d_{i=1}
\overline{\omega}_i\, \epsilon_i \in \CP$ satisfying $u
(\OH\otimes\OP) = \OH \otimes \OK$ (see Section 3.1) and the
operator $a_{\OP} = \sum^d_{i=1} \overline{\omega}_i\, a_i$.
\end{rem}
For later use we show

\begin{prop}
Let $\underline{A} = (A_1,\ldots,A_d)$ be an ergodic coisometric row
contraction such that $A^*_i\, \OH = \overline{\omega}_i\, \OH$ for
all $i$, further $A_{\OP} := \sum^d_{i=1} \overline{\omega}_i\,
A_i$. Then for $n\to\infty$ in the strong operator topology
\[
(A^*_{\OP})^n \to | \OH \rangle \langle \OH |.
\]
\end{prop}

\begin{proof}
We use the setting of Section 3.1 to be able to apply Theorem 3.1.1.
From $u^* (h \otimes \OK) = \sum^d_{i=1} a^*_i(h) \otimes
\epsilon_i$ we obtain
\[
u^* (h \otimes \OK) = a^*_{\OP}(h) \otimes \OP\; \oplus h^\prime\]
with $h^\prime \in \CH \otimes \OP^\perp$. Assume that $h \in \Ho$.
Because $u^*$ is isometric on $\CH \otimes \OK$ we conclude that
\begin{equation}
u^*(\OH\otimes\OK) = \OH \otimes \OP \perp u^*(h \otimes\OK)
\end{equation}
and thus also $a^*_{\OP}(h) \in \Ho$. In other words,
\[
a^*_{\OP}(\Ho) \subset\, \Ho .
\]
Let $q_n$ be the orthogonal projection from $\CH \otimes
\bigotimes^n_1 \CP$ onto $\OH \otimes \bigotimes^n_1 \CP$. From
Theorem 3.1.1 it follows that
\[
(\eins - q_n) u^*_{0n} \ldots u^*_{01} (h \otimes \bigotimes^n_1
\OK) \to 0 \quad (n\to\infty).
\]
On the other hand, by iterating the formula from the beginning,
\[
u^*_{0n} \ldots u^*_{01} (h \otimes \bigotimes^n_1 \OK) =
\big((a^*_{\OP})^n(h) \otimes \bigotimes^n_1 \OP \big) \oplus
h^\prime
\]
with $h^\prime \in \CH \otimes (\bigotimes^n_1 \OP)^\perp$. It
follows that also
\[
(\eins - q_n) \big( (a^*_{\OP})^n(h) \otimes \bigotimes^n_1 \OP
\big) \to 0.
\]
But from $a^*_{\OP}(\Ho) \subset\, \Ho $ we have $q_n \big(
(a^*_{\OP})^n(h) \otimes \bigotimes^n_1 \OP \big) = 0$ for all $n$.
We conclude that $(a^*_{\OP})^n(h) \to 0$ for $n\to\infty$. Further
\[
a^*_{\OP} \OH = \sum^d_{i=1} \omega_i\, a^*_i\, \OH = \sum^d_{i=1}
\omega_i\, \overline{\omega}_i\, \OH = \OH,
\]
and the proposition is proved.
\end{proof}

The following proposition summarizes some well known properties of
minimal isometric dilations and associated Cuntz algebra
representations.

\begin{prop}
Suppose $\uA$ is a coisometric tuple on $\CH$ and $\uV$ is its
minimal isometric dilation. Assume $\Omega_\CH$ is a distinguished
unit vector in $\CH$ and $\uomega = (\omega_1,\ldots,\omega_d) \in
\C^d,\;\sum_i |\omega_i|^2=1$. Then the following are equivalent.
\begin{enumerate}
\item $\uA$ is ergodic and $A^*_i\, \OH = \overline{\omega}_i\, \OH$ for all $i$.

\item $\uV$ is ergodic and
$V^*_i\, \OH = \overline{\omega}_i\, \OH$ for all $i$.

\item $V^*_i\, \OH = \overline{\omega}_i\, \OH$ and $\uV$
generates the GNS-representation of the Cuntz algebra $\CO_d =
C^*\{g_1,\cdots,g_d \}$ ($g_i$ its abstract generators) with respect
to the Cuntz state which maps
\[ g_\alpha\, g^*_\beta \mapsto
\omega_\alpha\, \oomega_\beta, ~~\forall \alpha, \beta \in
\tilde{\Lambda}.
\]
\end{enumerate}
Cuntz states are pure and the corresponding GNS-representations are
irreducible.
\end{prop}

This Proposition clearly follows from Theorem 5.1 of \cite{BJKW00},
Theorem 3.3 and Theorem 4.1 of \cite{BJP96}. Note that in Theorem
3.1.1(d) we already saw a concrete version of the corresponding
Cuntz algebra representation.

\section{A new characteristic function}

First we recall some more details of the theory of minimal isometric
dilations for row contractions (cf. \cite{Po89a}) and introduce
further notation.

The full Fock space over $\C^d$ ($d\geq 2$) denoted by $\Gamma
(\C^d)$ is
$$\Gamma (\C^d) := \mathbb{C}\oplus \C^d \oplus (\C^d)^{\otimes ^2}\oplus \cdots
\oplus (\C^d) ^{\otimes ^m}\oplus \cdots.
$$
$1\oplus 0\oplus \cdots$ is called the vacuum vector. Let $\{e_1,
\ldots , e_d\}$ be the standard orthonormal basis of $\C^d$. Recall
that we include $d=\infty$ in which case $\C^d$ stands for a complex
separable Hilbert space of infinite dimension. For $\alpha \in
{\tilde \Lambda}$, $e_{\alpha }$ will denote the vector $e_{\alpha
_1}\otimes e_{\alpha _2}\otimes \cdots \otimes e_{\alpha _m}$ in the
full Fock space $\Gamma (\mathbb{C}^d)$ and $e_0$ will denote the
vacuum vector. Then the (left) creation operators $L_i$ on
$\Gamma({\mathbb{C}^d})$ are defined by
\[ L_i x = e_i \otimes x \]  for $1 \leq i \leq d$ and $x
\in \Gamma({\mathbb{C}}^d).$ The row contraction $\underline{L}=
(L_1, \ldots , L_d)$ consists of isometries with orthogonal ranges.
\\

Let $\uT=(T_1,\cdots,T_d)$ be a row contraction on a Hilbert space
$\CH$. Treating $\uT$ as a row operator from $ \bigoplus^d_{i=1}
\CH$ to $\CH,$ define $D_*:= (\eins-\uT \uT^*)^\frac{1}{2}: \CH
\rightarrow \CH$ and $D:=(\eins-\uT^*\uT)^\frac{1}{2}:
\bigoplus^d_{i=1} \CH \rightarrow \bigoplus^d_{i=1} \CH$. This
implies that
\begin{equation}
D_*= (\eins -\sum^d_{i=1} T_i T^*_i)^\frac{1}{2}, ~~~D=(\delta_{ij}
\eins -T_i^*T_j)^\frac{1}{2}_{d\times d}.
\end{equation}
Observe that $\uT D^2 = D^2_* \uT$ and hence $\uT D = D_* \uT.$ Let
$\CD:=\mbox{Range~} D$ and $\CD_*:=\mbox{Range~} D_*.$ Popescu in
\cite{Po89a} gave the following explicit presentation of the minimal
isometric dilation of $\uT$ by $\uV$ on $\CH  \oplus (\Gamma (\C^d)
\otimes \CD)$,
\begin{equation}
V_i(h \oplus \sum_{\alpha \in \tilde{\Lambda}}  e_\alpha \otimes
d_\alpha) = T_i h \oplus [e_0 \otimes D_i h + e_i \otimes
\sum_{\alpha \in \tilde{\Lambda}} e_\alpha \otimes d_\alpha]
\end{equation}
for $h \in \CH$ and $d_\alpha \in \CD.$ Here $D_i h := D
(0,\ldots,0,h,0,\ldots,0)$ and $h$ is embedded at the $i^{th}$
component.

In other words, the $V_i$ are isometries with orthogonal ranges such
that $T^*_i = V^*_i |_{\CH}$ for $i=1,\ldots,d$ and the spaces
$V_\alpha \CH$ with $\alpha \in \tilde{\Lambda}$ together span the
Hilbert space on which the $V_i$ are defined. It is an important
fact, which we shall use repeatedly, that such minimal isometric
dilations are unique up to unitary equivalence (cf. \cite{Po89a}).
\\

Now, as in Section 3.2, let $\uA =(A_1,\cdots,A_d),\, A_i \in
B(\CH)$, be an ergodic coisometric tuple with $A^*_i \Omega_\CH =
\oomega_i \Omega_\CH$ for some unit vector $\Omega_\CH \in \CH$ and
some $\uomega \in \C^d,\;\sum_i |\omega_i|^2=1$. Let $\uV =
(V_1,\cdots,V_d)$ be the minimal isometric dilation of $\uA$ given
by Popescu's construction (see equation 3.4) on $\CH \oplus
\big(\Gamma (\C^d) \otimes \CD_A \big)$. Because $A^*_i = V^*_i
|_\CH$ we also have $V^*_i \Omega_\CH = \oomega_i \Omega_\CH$ and
because $\uV$ generates an irreducible $\Od-$representation
(Proposition 3.2.4), we see that $\uV$ is also a minimal isometric
dilation of $\uomega: \C^d \rightarrow \C$. In fact, we can think of
$\uomega$ as the most elementary example of a tuple with all the
properties stated for $\uA$. Let $\tilde{\uV} =
(\tilde{V}_1,\cdots,\tilde{V}_d)$ be the minimal isometric dilation
of $\uomega$ given by Popescu's construction on $\C \oplus (\Gamma
(\C^d) \otimes \CD_\omega)$.
\\

Because $\uA$ is coisometric it follows from equation 3.3 that $D$
is in fact a projection and hence
$D=(\delta_{ij}\eins-A_i^*A_j)_{d\times d}.$ We infer that $D
(A^*_1, \cdots,A^*_d)^T =0$, where $T$ stands for transpose. Applied
to $\uomega$ instead of $\uA$ this shows that $D_{\omega} = (\eins -
|\overline{\uomega} \rangle \langle \overline{\uomega} |)$ and
$$\CD_\omega \oplus \C (\oomega_1, \cdots, \oomega_d)^T = \C^d,$$
where $\underline{\oomega} = (\oomega_1, \cdots, \oomega_d)$.
\\

\begin{rem}
Because $\OH$ is cyclic for $\{ V_\alpha,\, \alpha \in
\tilde{\Lambda} \}$ we have
$$\overline{\mbox{span}}\{A_\alpha \Omega_{\CH}:\alpha \in
\tilde{\Lambda}\}=\overline{\mbox{span}}\{p_\CH\,V_\alpha
\Omega_{\CH}: \alpha \in \tilde{\Lambda} \}=\; \CH.$$ Using the
notation from equation 3.1 this further implies that
$$\overline{\mbox{span}}\{ \Ao_\alpha\, l_i: \alpha \in \tilde{\Lambda},
1 \leq i \leq d\} =\; \Ho .$$
\end{rem}
\vspace{0.4cm}

As minimal isometric dilations of the tuple $\uomega$ are unique up
to unitary equi\-valence, there exists a unitary
$$W: \CH \oplus (\Gamma (\C^d) \otimes \CD_A) \to
\C \oplus (\Gamma (\C^d) \otimes \CD_\omega),$$ such that $W V_i =
\tilde{V}_i W$ for all $i.$

After showing the existence of $W$ we now proceed to compute $W$
explicitly. For $\uA$, by using Popescu's construction, we have its
minimal isometric dilation $\uV$ on $\CH \oplus (\Gamma (\C^d)
\otimes \CD_A).$ Another way of constructing a minimal isometric
dilation $\ut$ of $\ua$ was demonstrated in Section 3.1 on the space
$\hat{\CH}$ (obtained by restricting to the minimal subspace of $\CH
\otimes \tilde{\CK}$ with respect to $\ut$). Identifying $\uA$ and
$\ua$ on the Hilbert space $\CH$ there is a unitary $\Gamma_A:
\hat{\CH} \rightarrow \CH \oplus (\Gamma (\C^d) \otimes \CD_A)$
which is the identity on $\CH$ and satisfies $V_i \Gamma_A =
\Gamma_A t_i$.
\\

By Theorem 3.1.1(d) the tuple $\us$ on $\tilde{\CP}$ arising from
the tensor shift is unitarily equivalent to $\ut$ (resp. $\uV$),
explicitly ${\bf w}\,t_i = s_i\,{\bf w}$ for all $i$. An alternative
viewpoint on the existence of ${\bf w}$ is to note that $\us$ is a
minimal isometric dilation of $\uomega.$ In fact, $s^*_i\,
\Omega_{\tilde{\CP}} = \langle \epsilon_i, \OP \rangle
\Omega_{\tilde{\CP}} = \oomega_i\, \Omega_{\tilde{\CP}}$ for all
$i$. Hence there is also a unitary $\Gamma_\omega: \tilde{\CP}
\rightarrow \C \oplus (\Gamma (\C^d) \otimes \CD_\omega)$ with
$\Gamma_\omega \Omega_{\tilde{\CP}} = 1 \in \C$ which satisfies
$\tilde{V}_i \Gamma_\omega = \Gamma_\omega s_i$.
\\

\begin{rem}
It is possible to describe $\Gamma_\omega$ in an explicit way and in
doing so to construct an interesting and natural (unitary)
identification of $\bigotimes^\infty_1 \C^d$ and $\C \oplus (\Gamma
(\C^d) \otimes \C^{d-1})$. In fact, recall (from Section 3.1) that
$\tilde{\CP} = \bigotimes^\infty_1 \CP$ and the space $\CP$ is
nothing but a $d$-dimensional Hilbert space. Hence we can identify
\[
\C^d \simeq \CP =\, \Po \oplus\, \C \OP \simeq \CD_\omega \oplus
\C\, \underline{\oomega}^T \simeq \C^{d-1} \oplus \C
\]
In this identification the orthonormal basis $(\epsilon_i)^d_{i=1}$
of $\CP$ goes to the canonical basis $(e_i)^d_{i=1}$ of $\C^d$, in
particular the vector $\OP = \sum_i \oomega_i\, \epsilon_i$ goes to
$\underline{\oomega}^T = (\oomega_1, \cdots, \oomega_d)^T$ and we
have $\Po\, \simeq \CD_\omega$. Then we can write
\begin{eqnarray*}
\Gamma_\omega:\qquad \Omega_{\tilde{\CP}} &\mapsto& 1 \in \C, \\
k \otimes \Omega_{\tilde{\CP}} &\mapsto& e_0 \otimes k \\
\epsilon_\alpha \otimes k \otimes \Omega_{\tilde{\CP}} &\mapsto&
e_\alpha \otimes k,
\end{eqnarray*}
where $k \in \Po,\, \alpha \in \tilde{\Lambda},\, \epsilon_\alpha =
\epsilon_{\alpha_1} \otimes \ldots \epsilon_{\alpha_n} \in
\bigotimes^n_1 \CP$ (the first $n$ copies of $\CP$ in the infinite
tensor product $\tilde{\CP}$), $e_\alpha = e_{\alpha_1} \otimes
\ldots e_{\alpha_n} \in \Gamma(\C^d)$ as usual. It is easily checked
that $\Gamma_\omega$ given in this way indeed satisfies the equation
$\tilde{V}_i \Gamma_\omega = \Gamma_\omega s_i$ (for all $i$), which
may thus be seen as the abstract characterization of this unitary
map (together with $\Gamma_\omega \Omega_{\tilde{\CP}} = 1$).
\end{rem}
\vspace{0.4cm}

Summarizing, for $i=1,\ldots,d$
\begin{eqnarray*}
V_i\, \Gamma_A = \Gamma_A\, t_i,\quad {\bf w}\, t_i = s_i\, {\bf
w},\quad \tilde{V}_i\, \Gamma_\omega = \Gamma_\omega\, s_i
\end{eqnarray*}
and we have the commuting diagram
\begin{equation}
\xymatrix{ \hat{\CH} \ar[r]^{{\bf w}} \ar[d]_{\Gamma_A} &
\tilde{\CP} \ar[d]^{\Gamma_\omega}
\\
\CH \oplus (\Gamma (\C^d) \otimes \CD_A) \ar[r]^{W} & \C \oplus
(\Gamma (\C^d) \otimes \CD_\omega). }
\end{equation}
From the diagram we get
$$W=\Gamma_\omega {\bf w} \Gamma^{-1}_A.$$
Combined with the equations above this yields $W V_i =
\tilde{V}_i\,W$ and we see that $W$ is nothing but the
dilations-intertwining map which we have already introduced earlier.
Hence ${\bf w}$ and $W$ are essentially the same thing and for the
study of certain problems it may be helpful to switch from one
picture to the other.

In the following we analyze $W$ to arrive at an interpretation as a
new kind of characteristic function. First we have an isometric
embedding
\begin{equation}
\hat{C} := W |_\CH : \CH \rightarrow \C \oplus (\Gamma (\C^d)
\otimes \CD_\omega).
\end{equation}
Note that $\hat{C}\,\OH = W\,\OH = 1 \in \C$. The remaining part is
an isometry
\begin{equation}
M_{\hat{\Theta}} := W |_{\Gamma (\C^d) \otimes \CD_A} : \Gamma
(\C^d) \otimes \CD_A \rightarrow \Gamma (\C^d) \otimes \CD_\omega .
\end{equation}
From equation 3.4 we get for all $i$
\[
V_i\, |_{\Gamma (\C^d) \otimes \CD_A} = (L_i \otimes \eins_{\CD_A}),
\]
\[
\tilde{V}_i\, |_{\Gamma (\C^d) \otimes \CD_\omega} = (L_i \otimes
\eins_{\CD_\omega}),
\]
and we conclude that
\begin{equation}
M_{\hat{\Theta}} (L_i \otimes \eins_{\CD_A}) = (L_i \otimes
\eins_{\CD_\omega}) M_{\hat{\Theta}},~~\forall 1 \leq i \leq d.
\end{equation}
In other words, $M_{\hat{\Theta}}$ is a multi-analytic inner
function in the sense of \cite{Po89c,Po95}. It is determined by its
symbol
\begin{equation}
\hat{\theta} := W|_{e_0 \otimes \CD_A} : \CD_A \rightarrow \Gamma
(\C^d) \otimes \CD_\omega,
\end{equation}
where we have identified $e_0 \otimes \CD_A$ and $\CD_A$. In other
words, we think of the symbol $\hat{\theta}$ as an isometric
embedding of $\CD_A$ into $\Gamma (\C^d) \otimes \CD_\omega$.

\begin{defn}
We call $M_{\hat{\Theta}}$ (or $\hat{\theta}$) the extended
characteristic function of the row contraction $\uA$,
\end{defn}

See Sections 3.5 and 3.6 for more explanation and justification of
this terminology.

\section[Explicit computation]{Explicit computation of the extended characteristic function}

To express the extended characteristic function more explicitly in
terms of the tuple $\uA$ we start by defining
\begin{eqnarray}
\hat{D}_*:\; \Ho = \CH \ominus \C \OH &\rightarrow& \Po\, = \CP
\ominus \C \OP \simeq \CD_\omega,
\end{eqnarray}
\[
h \mapsto \big( \langle \OH | \otimes \eins_{\CP} \big) \, u^*(h
\otimes \OK),
\]
where $u: \CH \otimes \CP \rightarrow \CH \otimes \CK$ is the
isometry introduced in Section 3.1. That indeed the range of
$\hat{D}_*$ is contained in $\Po$ follows from equation 3.2, i.e.,
$u^*(h \otimes\OK) \perp \OH \otimes \OP$ for $h \in \Ho$. With
notations from equation 3.1 we can get a more concrete formula.

\begin{lem}
For all $h \in \Ho$ we have $\hat{D}_*(h) = \sum^d_{i=1} \langle
\ell_i, h \rangle \epsilon_i$.
\end{lem}

\begin{proof}
$ \big( \langle \OH | \otimes \eins_{\CP} \big) \, u^*(h \otimes
\OK) = \sum^d_{i=1} \langle \OH, a^*_i h \rangle \otimes \epsilon_i
= \sum^d_{i=1} \langle \ell_i, h \rangle \epsilon_i. $
\end{proof}
\vspace{0.2cm}

\begin{prop}
The map $\hat{C}:\CH \to \C \oplus (\Gamma (\C^d) \otimes
\CD_\omega)$ from equation 3.7 is given explicitly by
$\hat{C}\Omega_\CH = 1$ and for $h \in \Ho$ by
$$\hat{C} h = \sum_{\alpha\in\tilde{\Lambda}} e_\alpha \otimes \hat{D}_* \Aso_\alpha h.$$
\end{prop}
\begin{proof}
As $W \Omega_\CH = 1$ also $\hat{C} \Omega_\CH =1$. Assume $h \in
\Ho$. Then
\begin{eqnarray*}
u_{01}(h \otimes \Omega_{\tilde{\CK}}) &=&\sum_i a^*_i h \otimes
\epsilon_i
\otimes \Omega_{\tilde{\CK}}\\
& = & \sum_i \langle \ell_i, h \rangle \Omega_\CH \otimes \epsilon_i
\otimes \Omega_{\tilde{\CK}} + \sum_i \aso_i\!h \otimes \epsilon_i
\otimes \Omega_{\tilde{\CK}}.
\end{eqnarray*}
Because $u^*(\OH \otimes \OK) = \OH \otimes \OP$ we obtain (with
Lemma 3.4.1) for the first part
\begin{eqnarray*}
&&\lim_{n \to \infty} u^*_{0n} \cdots u^*_{02} (\sum_i \langle
\ell_i, h \rangle \Omega_\CH \otimes \epsilon_i
\otimes \Omega_{\tilde{\CK}})\\
&=&\sum_i \langle \ell_i, h \rangle \Omega_\CH \otimes \epsilon_i
\otimes \Omega_{\tilde{\CP}} = \OH \otimes \hat{D}_* h \otimes
\Omega_{\tilde{\CP}} \simeq \hat{D}_* h \otimes \Omega_{\tilde{\CP}}
\in \tilde{\CP}.
\end{eqnarray*}
Using the product formula from Theorem 3.1.1 and iterating the
argument above we get
\[
\hat{C}(h) \,=\, Wh \,=\, \Gamma_\omega {\bf w} \Gamma^{-1}_A(h)
\]
\[
= \Gamma_\omega(\hat{D}_* h \otimes \Omega_{\tilde{\CP}}) \,+\,
\Gamma_\omega  \lim_{n \to \infty} u^*_{0n} \cdots u^*_{02} \sum_i
\aso_i\!h \otimes \epsilon_i \otimes \Omega_{\tilde{\CK}}
\]
\[
= e_0 \otimes \hat{D}_* h \,+\, \Gamma_\omega  \lim_{n \to \infty}
u^*_{0n} \cdots u^*_{03} \sum_{j,i} \big(\langle \ell_j, \aso_i\!h
\rangle \OH \,+\, \aso_j \aso_i\!h \big) \otimes \epsilon_i \otimes
\epsilon_j \otimes \Omega_{\tilde{\CK}}
\]
\[
= e_0 \otimes \hat{D}_* h + \sum^d_{i=1} e_i \otimes \hat{D}_*
\aso_i\!h \,+\, \Gamma_\omega  \lim_{n \to \infty} u^*_{0n} \cdots
u^*_{03} \sum_{j,i} \aso_j \aso_i\!h \otimes \epsilon_i \otimes
\epsilon_j \otimes \Omega_{\tilde{K}}
\]
\[
= \ldots
\]
\[
= \sum_{|\alpha|<m} e_\alpha \otimes \hat{D}_* \aso_\alpha h \,+\,
\Gamma_\omega  \lim_{n \to \infty} u^*_{0n} \cdots u^*_{0,m+1}
\sum_{|\alpha|=m} \aso_\alpha h \otimes \epsilon_\alpha \otimes
\Omega_{\tilde{\CK}}.
\]
From Proposition 3.2.1 we have $\sum_{|\alpha|=m} \| \aso_\alpha\!h
\|^2 \to 0$ for $m \to \infty$ and we conclude that the last term
converges to $0$. It follows that the series converges and this
proves Proposition 3.4.2.
\end{proof}
\begin{rem}
Another way to prove Proposition 3.4.2 for $h \in \Ho$ consists in
repeatedly applying the formula
\[
u^* (h\otimes \OK) = a^*_{\OP} h \otimes \OP + h^\prime, \quad
h^\prime \in \CH \otimes \Po
\]
to the $u^*_{0n}(h\otimes \OK)$ and then using $(a^*_{\OP})^n h \to
0$, see Proposition 3.2.3. This gives some insight how the infinite
product in Theorem 3.1.1 transforms into the infinite sum in
Proposition 3.4.2.
\end{rem}

Now we present an explicit computation of the extended
characteristic function. One way of writing $\CD_A$ is
$$\CD_A= \overline{\mbox{span}}\{ (V_i - A_i) h: i \in \Lambda, h \in \CH\}.$$
Let $d^i_h:=(V_i -A_i)h.$ Then
$$\hat{\theta}\, d^i_h =W(V_i-A_i) h = \tilde{V}_i \hat{C} h - \hat{C} A_i h.$$

{\bf Case I:} Take $h = \Omega_\CH.$
$$\tilde{V}_i \hat{C} \Omega_\CH = \tilde{V}_i 1 = \omega_i \oplus
[e_0 \otimes (\eins -  |\underline{\oomega} \rangle \langle
\underline{\oomega} |) \epsilon_i],$$
$$\hat{C}A_i\, \Omega_\CH = \omega_i \oplus \sum_\alpha e_\alpha \otimes \hat{D}_* \Aso_\alpha l_i$$
and thus
\begin{eqnarray*}
\hat{\theta}\; d^i_{\Omega_\CH} &=& e_0 \otimes [(\eins -
|\underline{\oomega} \rangle \langle \underline{\oomega} |)
\epsilon_i- \hat{D}_*l_i] - \sum_{|\alpha|\geq 1} e_\alpha \otimes
\hat{D}_*
\Aso_\alpha l_i\\
&=& e_0 \otimes [\epsilon_i - \sum_j \oomega_j \omega_i \epsilon_j -
\sum_j \langle l_j, l_i\rangle \epsilon_j] - \sum_{|\alpha|\geq 1}
e_\alpha \otimes \sum_j \langle l_j,
\Aso_\alpha l_i \rangle \epsilon_j\\
&=&e_0 \otimes [\epsilon_i - \sum_j (\oomega_j \omega_i + \langle
l_j, l_i\rangle) \epsilon_j] - \sum_{|\alpha|\geq 1} e_\alpha
\otimes \sum_j \langle \Ao_\alpha l_j, l_i \rangle \epsilon_j
\end{eqnarray*}
\begin{equation}
= e_0 \otimes [\epsilon_i - \sum_j \langle A_j\,\OH, A_i\,\OH
\rangle\, \epsilon_j] - \sum_{|\alpha|\geq 1} e_\alpha \otimes
\sum_j \langle \Ao_\alpha l_j, l_i \rangle \epsilon_j.
\end{equation}

{\bf Case II:} Now let $h \in\, \Ho$. With $i \in \Lambda$
\[
\tilde{V}_i\, \hat{C} h = (L_i \otimes \eins) \hat{C} h =
\sum_\alpha e_i \otimes e_\alpha \otimes \hat{D}_* \Aso_\alpha h,$$
$$\hat{C}A_i h= \sum_\beta e_\beta \otimes \hat{D}_* \Aso_\beta \Ao_i h.
\]
Finally
\begin{equation*}
\hat{\theta}\; d^i_h =  \sum_\alpha e_i \otimes e_\alpha \otimes
\hat{D}_* \Aso_\alpha h- \sum_\beta e_\beta \otimes \hat{D}_*
\Aso_\beta \Ao_i h
\end{equation*}
\begin{equation*}
=\; - e_0 \otimes \hat{D}_* \Ao _i h +  e_i \otimes \sum_\alpha
e_\alpha \otimes \hat{D}_* \Aso_\alpha (\eins - \Aso_i \Ao_i) h +
\sum_{j \neq i} e_j \otimes \sum_\alpha e_\alpha \otimes \hat{D}_*
\Aso_\alpha (-\Aso_j \Ao_i) h
\end{equation*}
\begin{equation}
=\; - e_0 \otimes \hat{D}_* \Ao _i h + \sum^d_{j=1} e_j \otimes
\sum_{\alpha \in \tilde{\Lambda}} e_\alpha \otimes \hat{D}_*
\Aso_\alpha (\delta_{ji} \eins - \Aso_j \Ao_i) h .
\end{equation}

\section{Case II is Popescu's characteristic function}

In this section we show that case II in the previous section can be
identified with the characteristic function of the $*$-stable tuple
$\uAo$, in the sense introduced by Popescu in \cite{Po89b}. This is
the reason why we have called $\hat{\theta}$ an {\it extended}
characteristic function. All information about $\uA$ beyond $\uAo$
must be contained in case I.

First recall the theory of characteristic functions for row
contractions, as developed by G.\,Popescu in \cite{Po89b},
generalizing the theory of B.\,Sz.-Nagy and C.\,Foias (cf.
\cite{NF70}) for single contractions. We only need the results about
a $*$-stable tuple $\uAo = (\Ao_1,\ldots,\Ao_d)$ on $\Ho$. In this
case, with $\Dso \;=\; (\eins- \uAo \uAso)^{\frac{1}{2}}:\; \Ho
\rightarrow \Ho$ and $\CDso$ its range, the map
\begin{equation}
\Co:\; \Ho \to \Gamma (\C^d) \otimes \CDso
\end{equation}
\[
h \mapsto \sum_{\alpha \in \tilde{\Lambda}} e_\alpha \otimes \Dso
\Ao_\alpha^* h
\]
is an isometry (Popescu's Poisson kernel). If, as usual, $\Do \;=\;
(\eins - \uAso \uAo)^{\frac{1}{2}}: \bigoplus^d_1 \Ho \rightarrow
\bigoplus^d_1 \Ho$, with $\CDo$ its range, and if $P_j$ is the
projection onto the $j$-th component, then the characteristic
function $\theta_{\Ao}$ of $\uAo$ can be defined as
\begin{equation}
\theta_{\Ao}:\; \CDo \to \Gamma (\C^d) \otimes \CDso
\end{equation}
\[
f \mapsto -e_0 \otimes \sum^d_{j=1} \Ao_j P_j f + \sum^d_{j=1} e_j
\otimes \sum_{\alpha \in \tilde{\Lambda}} e_\alpha  \otimes \Dso
\Ao_\alpha^* P_j \Do f.
\]
See \cite{Po89b} for details, in particular for the important result
that $\theta_{\Ao}$ characterizes the $*$-stable tuple $\uAo$ up to
unitary equivalence.
\\

Now consider again the tuple $\uA$ of the previous section, with
extended characteristic function $\hat{\theta}$. From equation 3.1
\[
A_i = \left( \begin{array}{cc}
        \omega_i & 0 \\
        |\ell_i\rangle & \Ao_i \\
        \end{array}
\right), \quad A^*_i =\left( \begin{array}{cc}
        \overline{\omega}_i & \langle \ell_i| \\
        0 & \Aso_i \\
        \end{array}
\right)
\]
and hence
\[ A_i A^*_i = \left( \begin{array}{cc} |\oomega_i|^2 & \langle \oomega_i l_i |\\
|\oomega_i l_i \rangle & |l_i \rangle \langle l_i| + \Ao_i \Aso_i
\end{array} \right).
\]
Recall that $D^2_*= \eins - \sum_i A_i A^*_i$ which is $0$ as $\uA$
is coisometric. Thus $\sum_i \oomega_i l_i =0$ and $\eins-\sum_i
\Ao_i \Aso_i= \sum_i |l_i \rangle \langle l_i|.$ The first equation
means that $A^*_{\OP}(\Ho) \subset \Ho$ and that
\[
\langle \hat{\CD}_* h, \OP \rangle = \langle \sum_i \langle \ell_i,
h \rangle \epsilon_i, \sum_j \oomega_j \epsilon_j \rangle = \langle
\sum_i \oomega_i \ell_i, h \rangle = 0,
\]
which we already know (see 3.10). \\
The second equation yields
\[
\Dso^2 = \eins-\sum_i \Ao_i \Aso_i= \sum_i |l_i \rangle \langle
l_i|.
\]

\begin{lem}
There exists an isometry $\gamma:\, \CDso\, \rightarrow\, \Po\,
\simeq \CD_{\omega}$ defined for $h \in \Ho$ as
\[
\Dso h \mapsto \sum_i \langle l_i,h \rangle \epsilon_i = \hat{D}_*
h.
\]
\end{lem}
\begin{proof}
Take $h \in \Ho.$ By Lemma 3.4.1 we have $\hat{D}_*(h) =
\sum^d_{i=1} \langle \ell_i, h \rangle \epsilon_i$. Now we can
compute
\[
\|\hat{D}_* h\|^2 = \langle \sum_i \langle l_i,h \rangle \epsilon_i,
\sum_j \langle l_j,h \rangle \epsilon_j\rangle = \sum_i \langle
h,l_i \rangle \langle l_i,h \rangle = \langle h, \Dso^2
h\rangle=\|\Dso h\|^2.
\]
Hence $\gamma:\; \Dso h \mapsto \hat{D}_* h$ is isometric.
\end{proof}

\begin{thm}
Let $\uA =(A_1,\cdots,A_d),\, A_i \in B(\CH)$, be an ergodic
coisometric tuple with $A^*_i \Omega_\CH = \oomega_i \Omega_\CH$ for
some unit vector $\Omega_\CH \in \CH$ and some $\uomega \in
\C^d,\;\sum_i |\omega_i|^2=1$. Let $\hat{\theta}$ be the extended
characteristic function of $\uA$ and let $\theta_{\Ao}$ be the
characteristic function of the ($*$-stable) tuple $\underline{\Ao}$.
For $h \in \Ho$
\begin{eqnarray*}
\gamma \Dso h &=& \hat{D}_*h,\\
(\eins \otimes \gamma)\Co h &=& \hat{C}h,
\end{eqnarray*}
\[
(\eins \otimes \gamma)\, \theta_{\Ao}\, d^i_h = \hat{\theta}\,
d^i_h~~ \mbox{for}~ i \in \Lambda.
\]
In other words, the part of $\hat{\theta}$ described by case II in
the previous section is equivalent to $\theta_{\Ao}$.
\end{thm}

\begin{proof}
We only have to use Lemma 3.5.1 and compare Proposition 3.4.2 and
equation 3.13 as well as equations 3.12 and 3.14. For the latter
note that $d^i_h = \Do(0,\ldots,0,h,0\ldots,0)$, where $h$ is
embedded at the $i$-th position. Hence
\[
\gamma \sum_j \Ao_j P_j d^i_h = \gamma \sum_j \Ao_j P_j
\Do(0,\ldots,0,h,0\ldots,0) = \gamma \uAo \;
\Do(0,\ldots,0,h,0\ldots,0)
\]
\[
= \gamma \Dso \uAo(0,\ldots,0,h,0\ldots,0) = \hat{D}_* \, \Ao_i h
\]
and also
\[
P_j \Do d^i_h = P_j \Do^2(0,\ldots,0,h,0\ldots,0) = (\delta_{ji}
\eins - \Aso_j \Ao_i) h .
\]
\end{proof}

Of course, Theorem 3.5.2 explains why we have called $\hat{\theta}$
an {\it extended} characteristic function.

\section[A complete unitary invariant]
{The extended characteristic function is a complete unitary invariant}

In this section we prove that the extended characteristic function
is a complete invariant with respect to unitary equivalence for the
row contractions investigated in this paper. Suppose that $\uA =
(A_1,\ldots,A_d)$ and $\uB = (B_1,\ldots,B_d)$ are ergodic and
coisometric row contractions on Hilbert spaces $\CH_A$ and $\CH_B$
such that $A^*_i \Omega_A = \oomega_i \Omega_A$ and $B^*_i \Omega_B
= \oomega_i \Omega_B$ for $i=1,\ldots,d$, where $\Omega_A \in \CH_A$
and $\Omega_B \in \CH_B$ are unit vectors and $\uomega =
(\omega_1,\ldots,\omega_d)$ is a tuple of complex numbers. Recall
from Remark 3.2.2 that it is no serious restriction of generality to
assume that it is the same tuple of complex numbers in both cases
because this can always be achieved by a transformation with a
unitary $d\times d-$matrix (with scalar entries). We will use all
the notations introduced earlier with subscripts $A$ or $B$.

Let us say that the extended characteristic functions
$\hat{\theta}_A$ and $\hat{\theta}_B$ are {\em equivalent} if there
exists a unitary $V: \CD_A \rightarrow \CD_B$ such that
$\hat{\theta}_A = \hat{\theta}_B\,V$. Note that the ranges of
$\hat{\theta}_A$ and $\hat{\theta}_B$ are both contained in $\Gamma
(\C^d) \otimes \CD_\omega$ and thus this definition makes sense. Let
us further say that $\uA$ and $\uB$ are {\em unitarily equivalent}
if there exists a unitary $U: \CH_A \rightarrow \CH_B$ such that $
U\,A_i = B_i\,U$ for $i=1,\ldots,d$. By ergodicity the unit
eigenvector $\Omega_A$ (resp. $\Omega_B$) is determined up to an
unimodular constant (see Theorem 3.1.1(b)) and hence in the case of
unitary equivalence we can always modify $U$ to satisfy additionally
$U\,\Omega_A = \Omega_B$.

\begin{thm}
The extended characteristic functions $\hat{\theta}_A$ and
$\hat{\theta}_B$ are equivalent if and only if $\uA$ and $\uB$ are
unitarily equivalent.
\end{thm}
\begin{proof}
If $\uA$ and $\uB$ are unitarily equivalent then all constructions
differ only by naming and it follows that $\hat{\theta}_A$ and
$\hat{\theta}_B$ are equivalent. Conversely, assume that there is a
unitary $V: \CD_A \rightarrow \CD_B$ such that $\hat{\theta}_A =
\hat{\theta}_B\,V$. Now from the commuting diagram 3.5 and the
definitions following it
\begin{eqnarray*}
W_B \CH_B &=& \C \oplus \big( \Gamma (\C^d) \otimes \CD_\omega \big)
\ominus M_{\hat{\Theta}_B} \big( \Gamma (\C^d) \otimes \CD_B \big) \\
&=& \C \oplus \big( \Gamma (\C^d) \otimes \CD_\omega \big)
\ominus M_{\hat{\Theta}_B} \big( \Gamma (\C^d) \otimes V\,\CD_A \big) \\
&=& \C \oplus \big( \Gamma (\C^d) \otimes \CD_\omega \big)
\ominus M_{\hat{\Theta}_A} \big( \Gamma (\C^d) \otimes \CD_A \big) \\
&=& W_A \CH_A,
\end{eqnarray*}
where we used equation 3.8, i.e., $M_{\hat{\Theta}} (L_i \otimes
\eins_{\CD}) = (L_i \otimes \eins_{\CD_\omega})
M_{\hat{\Theta}},~~\forall 1 \leq i \leq d$, to deduce
$M_{\hat{\Theta}_A} = M_{\hat{\Theta}_B} (\eins \otimes V)$ from
$\hat{\theta}_A = \hat{\theta}_B\,V$. Now we define the unitary $U$
by
\[
U := W^{-1}_B\,W_A |_{\CH_A}: \CH_A \rightarrow \CH_B.
\]
Because $W_A \Omega_A = 1 = W_B \Omega_B$ we have $U\,\Omega_A =
\Omega_B$. Further for all  $i=1,\ldots,d$ and $h\in \CH_A$,
\[
U A_i\, h = W^{-1}_B\,W_A\, A_i\,h = W^{-1}_B\,W_A P_{\CH_A}\,V^A_i
h = P_{\CH_B} W^{-1}_B\,W_A\,V^A_i h
\]
\[
= P_{\CH_B} W^{-1}_B\,\tilde{V}_i\, W_A h = P_{\CH_B}
V^B_i\,W^{-1}_B\,W_A h = B_i\,U h,
\]
i.e., $\uA$ and $\uB$ are unitarily equivalent.
\end{proof}
\begin{rem}
An analogous result for completely non-coisometric tuples has been
shown by G.\,Popescu in \cite{Po89b}, Theorem 3.5.4.
\end{rem}
Note further that if we change  $\uA = (A_1,\ldots,A_d)$ into
$\uA^\prime = (A^\prime_1,\ldots,A^\prime_d)$ by applying a unitary
$d\times d-$matrix with scalar entries (as described in Remark
3.2.2), then $\hat{\theta}_A = \hat{\theta}_{A^\prime}$. In fact,
this follows immediately from the definition of $W$ as an
intertwiner in Section 3.3, from which it is evident that $W$ does
not change if we take the same linear combinations on the left and
on the right. This does not contradict Theorem 3.6.1 because
$\uomega$ and $\uomega\prime$ are now different tuples of
eigenvalues and Theorem 3.6.1 is only applicable when the same tuple
of eigenvalues is used for $\uA$ and $\uB$.

For another interpretation, let $Z$ be a normal, unital, ergodic,
completely positive map with an invariant vector state $\langle
\Omega_A,\cdot\, \Omega_A \rangle$. If we consider two minimal Kraus
decompositions of $Z$, i.e.,
\[
Z = \sum^d_{i=1} A_i \cdot A^*_i = \sum^d_{i=1} A^\prime_i \cdot
(A^\prime_i)^*,
\]
with $d$ minimal, then the tuples $\uA = (A_1,\ldots,A_d)$ into
$\uA^\prime = (A^\prime_1,\ldots,A^\prime_d)$ are related in the way
considered above (see for example \cite{Go04}, A.2). It follows that
$\hat{\theta}_A = \hat{\theta}_{A^\prime}$ does not depend on the
decomposition but can be associated to $Z$ itself. Hence we have the
following reformulation of Theorem 3.6.1.

\begin{cor}
Let $Z_1,\,Z_2$ be normal, unital, ergodic, completely positive maps
on $B(\CH_1),\,B(\CH_2)$ with invariant vector states $\langle
\Omega_1,\cdot\, \Omega_1 \rangle$ and $\langle \Omega_2,\cdot\,
\Omega_2 \rangle$. Then the associated extended characteristic
functions $\hat{\theta}_1$ and $\hat{\theta}_2$ are equivalent if
and only if $Z_1$ and $Z_2$ are conjugate, i.e., there exists a
unitary $U: \CH_1 \rightarrow \CH_2$ such that
\[
Z_1(x) = U^* Z_2( U x U^* ) U \quad \mbox{for all}\; x \in B(\CH_1).
\]
\end{cor}

\section{Example}

The following example illustrates some of the constructions in this
paper.
\\

Consider $\CH=\C^3$ and
$$A_1=\frac{1}{\sqrt{2}}
\left( \begin{array}{ccc} 0 & 0 & 0 \\ 1 & 0 & 0 \\ 0 & 1 & 1
\end{array} \right), A_2= \frac{1}{\sqrt{2}}
\left( \begin{array}{ccc} 1 & 1 & 0 \\ 0 & 0 & 1 \\ 0 & 0 & 0
\end{array} \right).$$ Then $\sum^2_{i=1} A_iA^*_i=\eins.$ Take the
unital completely positive map $Z:M_3 \to M_3$ by $Z(x)=\sum^2_{i=1}
A_i x A^*_i.$ It is shown in Section 3.5 of \cite{GKL06} (and not
difficult to verify directly) that this map is ergodic. We will use
the same notations here as in previous sections. Observe that the
vector $\Omega_\CH:=\frac{1}{\sqrt{3}} (1,1,1)^T$ gives an invariant
vector state for $Z$ as
$$\langle \Omega_\CH,Z(x)\Omega_\CH \rangle =
\langle \Omega_\CH ,x \Omega_\CH \rangle = \frac{1}{3}\sum^3_{i,j=1}
x_{ij}.$$ $A^*_i \Omega_\CH = \frac{1}{\sqrt{2}}  \Omega_\CH$ and
hence $\uomega=\frac{1}{\sqrt{2}}(1,1).$ The orthogonal complement
$\Ho$ of $\C \Omega_\CH$ in $\C^3$ and the orthogonal projection $Q$
onto $\Ho$ are given by
$$\Ho=\{ \left(\begin{array}{c} k_1 \\ k_2 \\ -(k_1 + k_2)\end{array} \right):
k_1,k_2 \in \C\}, \qquad Q = \frac{1}{3} \left( \begin{array}{ccc} 2
& -1 & -1 \\ -1 & 2 & -1 \\ -1 & -1 & 2 \end{array} \right).
$$
From this we get for $\Ao_i = Q A_i Q = A_i Q$
$$\Ao_1 = \frac{1}{3 \sqrt{2}} \left( \begin{array}{ccc}
0 & 0 & 0 \\ 2 & -1 & -1 \\ -2 & 1 & 1 \end{array} \right), ~~\Ao_2
= \frac{1}{3 \sqrt{2}} \left( \begin{array}{ccc} 1 & 1 & -2 \\ -1 &
-1 & 2 \\ 0 & 0 & 0 \end{array} \right).$$ We notice that the tuple
$\uAo=(\Ao_1,\Ao_2)$ is $*$-stable as (by induction)
$$\sum_{|\alpha|=n} \Ao_\alpha \Aso_\alpha
= \frac{1}{3\times 2^{n-1}} \left( \begin{array}{ccc} 1 & -1 & 0 \\
-1 & 2 & -1 \\ 0 & -1 & 1 \end{array} \right) \to 0 \quad (n \to
\infty).
$$
Here $\CP=\C^2$ and $\Po:=\CP \ominus \C \Omega_\CP$ with $\OP =
\frac{1}{\sqrt{2}} (1,1)^T$. Easy calculation shows that
$\hat{D}_*:\; \Ho \to \Po$ is given by
$$\left( \begin{array}{c} k_1 \\ k_2 \\ -(k_1 + k_2)\end{array} \right)
\mapsto \frac{1}{\sqrt{6}} (2k_1 + k_2)\left( \begin{array}{c} -1 \\
1
\end{array} \right).$$
Moreover $\Dso=\frac{1}{\sqrt{6}} \left( \begin{array}{ccc} 1 & 0 &
-1 \\ 0 & 0 & 0 \\ -1 & 0 & 1 \end{array} \right).$ There exists an
isometry $\gamma:\,\CDso \to \Po$ such that $\left( \begin{array}{c}
1 \\ 0 \\ -1 \end{array} \right) \mapsto \left( \begin{array}{c} -1
\\ 1
\end{array} \right)$ and $\gamma(\Dso h) = \hat{D}_* h$ for $h \in \Ho$.

The map $\hat{C}: \CH \to \Gamma (\C^d) \otimes \CD_{\omega}$ is
given by $\hat{C}(\Omega_\CH)=1$ and for $h \in \Ho$ by
\begin{eqnarray*}
&&\hat{C} \left( \begin{array}{c} k_1 \\ k_2 \\ -(k_1 +
k_2)\end{array} \right)\\
 &=& e_0 \otimes \frac{(2k_1 +
k_2)}{\sqrt{6}} \left( \begin{array}{c} -1 \\ 1 \end{array}
\right) + \sum_{\alpha,\alpha_1=1} e_\alpha \otimes (\frac{1}{\sqrt{2}})^{|\alpha|} \frac{(k_1 + 2k_2)}{\sqrt{6}} \\
&& \times \left( \begin{array}{c} -1 \\ 1 \end{array} \right) +
\sum_{\alpha,\alpha_1=2} e_\alpha \otimes
(\frac{1}{\sqrt{2}})^{|\alpha|} \frac{(k_1 - k_2)}{\sqrt{6}} \left(
\begin{array}{c} -1 \\ 1 \end{array} \right)
\end{eqnarray*}
where the summations are taken over all $0 \not=\alpha \in
\tilde{\Lambda}$ such that $\alpha_i \neq \alpha_{i+1}$ for all $1
\leq i \leq |\alpha|$ and fixing $\alpha_1$ to $1$ or $2$ as
indicated. This simplification occurs because $\Ao^2_i = 0$ for
$i=1,2$. All the summations below in this section are also of the
same kind.

Now using the equations 3.11 and 3.12 for $\hat{\theta}_A :\CD_A \to
\Gamma (\C^d) \otimes \CD_{\omega}$ and simplifying we get
\begin{eqnarray*}
\hat{\theta}_A\, d^1_{\Omega_\CH}  &=& -e_0 \otimes \frac{1}{6}
\left( \begin{array}{c} -1 \\ 1 \end{array} \right) +
\sum_{\alpha,\alpha_1=1} e_\alpha \otimes
(\frac{1}{\sqrt{2}})^{|\alpha|} \frac{1}{6}
\left( \begin{array}{c} -1 \\ 1 \end{array} \right)\\
&&+ \sum_{\alpha,\alpha_1=2} e_\alpha \otimes
(\frac{1}{\sqrt{2}})^{|\alpha|} \frac{1}{6} \left( \begin{array}{c}
-1 \\ 1 \end{array} \right) = - \hat{\theta}_A\, d^2_{\Omega_\CH},
\end{eqnarray*}
and for $h \in \Ho,$
\begin{eqnarray*}
\hat{\theta}_A d^1_h
&=& -e_0 \otimes \frac{k_1}{2 \sqrt{3}} \left( \begin{array}{c} -1 \\
1 \end{array} \right)+ e_1 \otimes \frac{(k_1 + k_2)}{\sqrt{6}}
\left( \begin{array}{c} -1 \\ 1 \end{array} \right)
+ \sum_{\alpha,\alpha_1=1} e_1 \otimes e_\alpha \\
&&\otimes (\frac{1}{\sqrt{2}})^{|\alpha|} \frac{(k_1 +
2k_2)}{\sqrt{6}} \left( \begin{array}{c} -1 \\ 1 \end{array} \right)
- \sum_{\alpha,\alpha_1=2} e_1 \otimes e_\alpha \otimes
(\frac{1}{\sqrt{2}})^{|\alpha|} \frac{ k_2}{\sqrt{6}}
\left( \begin{array}{c} -1 \\ 1 \end{array} \right)\\
&&+ \sum_{\alpha,\alpha_1=2} e_\alpha \otimes
(\frac{1}{\sqrt{2}})^{|\alpha|} \frac{k_1}{2 \sqrt{3}} \left(
\begin{array}{c} -1 \\ 1 \end{array} \right),
\end{eqnarray*}

\begin{eqnarray*}
\hat{\theta}_A d^2_h&=& -e_0 \otimes \frac{(k_1 + k_2)}{2 \sqrt{3}}
\left( \begin{array}{c} -1 \\ 1 \end{array} \right) +
\sum_{\alpha,\alpha_1=1} e_\alpha \otimes
(\frac{1}{\sqrt{2}})^{|\alpha|} \frac{(k_1 + k_2)}{2 \sqrt{3}}
\left( \begin{array}{c} -1 \\ 1 \end{array}
\right)\\
&& + e_2 \otimes \frac{k_1}{\sqrt{6}} \left( \begin{array}{c} -1 \\
1
\end{array} \right) + \sum_{\alpha,\alpha_1=1} e_2 \otimes e_\alpha
\otimes (\frac{1}{\sqrt{2}})^{|\alpha|}\frac{k_2}{\sqrt{6}} \left(
\begin{array}{c} -1 \\ 1 \end{array}
\right)\\
&& + \sum_{\alpha,\alpha_1=2} e_2 \otimes e_\alpha \otimes
(\frac{1}{\sqrt{2}})^{|\alpha|} \frac{(k_1 - k_2)}{\sqrt{6}} \left(
\begin{array}{c} -1 \\ 1 \end{array} \right).
\end{eqnarray*}

Form this we can easily obtain $\Co$ and $\theta_{\Ao}$ for $h \in
\Ho$ by using the following relations from Theorem 3.5.2,
$$(\eins \otimes \gamma) \Co h = \hat{C} h,$$
$$(\eins \otimes \gamma) \theta_{\Ao} d^i_h = \hat{\theta}_A d^i_h.$$
Further
$$l_1= A_1 \Omega_\CH - \frac{1}{\sqrt{2}} \Omega_\CH=  \frac{1}{\sqrt{6}}
\left( \begin{array}{c} -1 \\ 0 \\ 1 \end{array} \right),~~ l_2= A_2
\Omega_\CH - \frac{1}{\sqrt{2}} \Omega_\CH=  \frac{1}{\sqrt{6}}
\left( \begin{array}{c} 1 \\ 0 \\ -1 \end{array} \right),$$ $\Ao_1
l_1=  \frac{1}{2 \sqrt{3}}  \left( \begin{array}{c} 0 \\ -1 \\ 1
\end{array} \right)$ and  clearly
$\Ho=\overline{\mbox{span}} \{ \Ao_\alpha l_i: i=1,2 \mbox{~and~}
\alpha \in \tilde{\Lambda}\},$ as already observed in Remark 3.1.

\section{Appendix}

Here for $*$-stable $\uT$ we generalize the
computation of \cite{Go06} to tuples. The following result
(together with the results in Section 3.1) made us to ask and
investigate the relation between Popescu's characteristic function
and ergodic tuples. Let
$$R_k:\CH \oplus (\oplus_{\alpha \in \tilde{\Lambda}} e_\alpha \otimes
\tilde{\CD}_\alpha ) \to \CH \oplus (\oplus_{\alpha \in
\tilde{\Lambda}} e_\alpha \otimes \CD ) = \CH \oplus (\Gamma (\C^d)
\otimes \CD) $$ where $ \tilde{\CD}_\alpha = \left\{
\begin{array}{cc} \CD_* & \mbox{if~} |\alpha|= k \\ \CD & \mbox{if~}
|\alpha|\neq k\end{array} \right.$ be given by

\begin{eqnarray*}
h \oplus \big(\sum_{\alpha \in \tilde{\Lambda}} e_\alpha \otimes
\tilde{d}_\alpha \big) &\mapsto& \{\sum_i T_ih + D_*
(\sum_{|\alpha|=k}\tilde{d}_\alpha)\} \oplus
[\sum_{|\alpha|<k} e_\alpha \otimes \tilde{d}_\alpha + \sum_{|\alpha|=k} e_\alpha \\
&& \otimes\{ D(h,\ldots,h)-(T^*_1 \tilde{d}_\alpha, \ldots,T^*_d
\tilde{d}_\alpha)\} + \sum_{|\alpha|>k} e_\alpha \otimes
\tilde{d}_\alpha].
\end{eqnarray*}

Let us use the presentation of the minimal isometric dilation given
by Popescu on $\CH  \oplus (\Gamma (\C^d) \otimes \CD)$. First
consider the isometry $U:=\frac{1}{\sqrt{d}}\sum_{i=1}^d V_i.$
\begin{eqnarray*}
&& d^{\frac{k}{2}}\, U^k(h\oplus \sum_{\delta \in \tilde{\Lambda}}
e_\delta \otimes d_\delta) \\
&=& \sum_{|\alpha|=k} T_\alpha h \oplus \sum^{k-1}_{i=0}
\sum_{|\beta|=i} e_\beta \otimes D(\sum_{|\gamma|=k-1-i} T_\gamma
h,\ldots,
\sum_{|\gamma|=k-1-i} T_\gamma h) \\
&& + \sum_{|\epsilon|=k} e_\epsilon \otimes \sum_{\delta \in
\tilde{\Lambda}} e_\delta \otimes d_\delta = R_0 \ldots R_{k-1}
(h\oplus \sum_{|\epsilon|=k} e_\epsilon \otimes \sum_{\delta \in
\tilde{\Lambda}} e_\delta \otimes d_\delta).
\end{eqnarray*}
We conclude that $U^k = R_0 \ldots R_{k-1} \big( \eins_{\CH} \oplus
(L \otimes \eins) \big)^k$ with $L := \frac{1}{\sqrt{d}}\sum_{i=1}^d
L_i$, i.e., the product of $R_i$'s is a kind of cocycle relating the
isometries $U$ and $L$.
\\

On the other hand we can use the product of adjoints to factorize
the unitary $\hat{W}$ corresponding to the characteristic function.
Note that
\begin{eqnarray*}
R^*_k(h \oplus \sum_{\alpha \in \tilde{\Lambda}} e_\alpha \otimes
d_\alpha) & = & (\sum_i T^*_ih + \sum_i \sum_{|\alpha|=k}
 P_i D d_\alpha) \oplus (\sum_{|\alpha|<k} e_\alpha \otimes d_\alpha \\
&& + \sum_{|\alpha|=k} e_\alpha \otimes(D_* h- \sum_i T_i P_i
d_\alpha) + \sum_{|\alpha|>k} e_\alpha \otimes d_\alpha),
\end{eqnarray*}
where $P_i$ is the projection onto the $i^{th}$ component of $\CD$.

\begin{eqnarray*}
&& R^*_{k-1} \ldots R^*_0(h \oplus \sum_{\alpha \in \tilde{\Lambda}}
e_\alpha \otimes d_\alpha) \\
&=& \{ \sum_{|\beta|=k} T_\beta^* h + \sum_{|\beta|=k-1} \sum_i
T_\beta^*  P_i D d_0 + \ldots + \sum_i \sum_{|\beta|=k-1}
P_i D d_\beta \} \\
&& \oplus\; [ (D_*h - \sum_i T_i P_i d_0 ) \\
&& + \sum_{|\beta|=1} e_\beta \otimes \{\sum_{|\gamma|=1}
D_*T_\gamma^* h + \sum_i D_* P_i D d_0 - \sum_{|\gamma|=1} \sum_i
T_i
P_i d_\gamma\}\\
&& +\ldots  + \sum_{|\beta|=k-1} e_\beta \otimes \{\sum_{|\gamma|=k-1} D_*T_\gamma^* h + \sum_{|\gamma|=k-2} D_*\sum_i T_\gamma^* P_i D d_0 + \ldots \\
&& - \sum_{|\gamma|=d-1} \sum_i T_i P_i d_\gamma\}  +
\sum_{|\alpha|>k-1} e_\alpha \otimes d_\alpha].
\end{eqnarray*}
Now the first bracket $\{\cdot\}$ converges to $0$ for $k \to
\infty$, and a comparison of the second bracket $[\cdot]$ with
Proposition 3.3.1 shows that
$$\hat{W}=\lim_{k\to \infty} R^*_{k-1} \ldots R^*_0,$$
which is analogous to the product formula for $\tilde{{\bf w}}$ in
Theorem 3.1.1.

\chapter{Characteristic Functions of Liftings}


{\bf Abstract:} {\em We introduce characteristic functions for
certain contractive liftings of row contractions. These are
multi-analytic operators which classify the liftings up to unitary
equivalence and provide a kind of functional model. The most
important cases are subisometric and coisometric liftings. We also
identify the most general setting which we call reduced liftings. We
derive properties of these new characteristic functions and discuss
the relation to Popescu's definition of the characteristic function
for completely non-coisometric row contractions. Finally we apply
our theory to completely positive maps and prove a one-to-one
correspondence between the fixed point sets of completely positive
maps related to each other by a subisometric lifting.}

\vspace{2cm}

Joint work with Rolf Gohm, accepted in Journal of Operator Theory.

\newpage

\section*{Introduction}

Let $C$ be a contraction on a Hilbert space $\CH_C$. Then a
contraction $E$ on a Hilbert space $\CH_E \supset \CH_C$ is called a
contractive lifting of $C$ if $P E = C P$, where $P$ is the
orthogonal projection from $\CH_E$ onto $\CH_C$. In other words, we
have an operator matrix
\begin{align} \label{0.1}
E = \left(
\begin{array}{cc}
C & 0 \\
B & A \\
\end{array}
\right).
\end{align}
See Chapter 5 of \cite{FF90}. In this book C.\,Foias and
A.E.\,Frazho amply demonstrate the importance of understanding the
structure of contractive liftings, in particular in connection with
the commutant lifting theorem and its applications.

The minimal isometric dilation (mid for short) of $C$ is the most
prominent example of a contractive lifting. In \cite{DF84}
R.G.\,Douglas and C.\,Foias introduced subisometric dilations (see
also Chapter 8.3 of \cite{Ber88} for a discussion closer to our
point of view). These are contractive liftings with the property
that the mid of $E$ is also minimal as an isometric dilation of $C$.
In this context Douglas and Foias were especially interested in
problems of uniqueness and of commutant lifting. We arrived at the
subisometric property in a completely different way and ask
different questions about it. Let us briefly describe the most
relevant aspects of this development.

Many results of the Sz.-Nagy/Foias-theory for contractions
\cite{NF70} can be generalized to row contractions $\underline{C} =
(C_1,\ldots,C_d)$, i.e. tuples of operators such that $\sum^d_{i=1}
C_i C^*_i \leq \eins$. This has been done very systematically by
G.\,Popescu starting with \cite{Po89a} and many people contributed
to this development, an incomplete list of work related to our
interests is \cite{Ar98,BBD04,BDZ06,BES05,DKS01,Po89b,Po03,Po05}. In
particular in \cite{Po89b} G.\,Popescu described a class of
multi-analytic operators which classify completely non-coisometric
(c.n.c.) row contractions up to unitary equivalence and called them
characteristic functions, in analogy to a similar concept in the
Sz.-Nagy/Foias-theory. In \cite{DG07a} S.\,Dey and R.\,Gohm started
from some seemingly unrelated questions in noncommutative
probability theory arising in \cite{Go04,GKL06} and established a
class of multi-analytic operators which are associated to certain
rather special coisometric row contractions (i.e., $\sum^d_{i=1} C_i
C^*_i = \eins$). Investigating their properties we came to the
conclusion that there are good reasons to think of them as of
characteristic functions for these tuples. This is not covered by
Popescu's theory.

In this paper we will show that it is the property of being a
subisometric lifting which makes this analysis possible. This is a
vast generalization of the setting of \cite{DG07a} and it clarifies
the mechanism behind it. It is straightforward to define liftings
for row contractions. Let $\uE = (E_1,\ldots,E_d)$ be a row
contraction on a Hilbert space $\CH_E \supset \CH_C$. If for all
$i=1,\ldots,d$ (with $d$ countable) we have an operator matrix
\begin{align} \label{0.2}
E_i =\left(
\begin{array}{cc}
C_i  & 0 \\
B_i  & A_i \\
\end{array}
\right)
\end{align}
with respect to $\CH_C \oplus \CH_C^\perp$ then we say that $\uE$ is
a lifting of $\uC=(C_1,\ldots,C_d)$ by $\uA=(A_1,\ldots,A_d)$ (or
that $\uE$ is an extension of $\uA$ by $\uC$). The subisometric
property in the form given here also makes sense for row
contractions, using Popescu's theory of mid for row contractions
\cite{Po89a}. This is worked out in Section 4.1 below. It then turns
out that there is a Beurling-type classification of subisometric
liftings, involving a correspondence to certain multi-analytic inner
operators (Theorem 4.1.6). They classify subisometric liftings up to
unitary equivalence, so we call them characteristic functions of
(subisometric) liftings.

In Section 4.2 we focus on coisometric liftings, i.e. $\sum^d_{i=1}
E_i E^*_i = \eins$, emphasizing another type of classification which
uses an isometry $\gamma$ mapping the defect space $\CD_{*,A}$ of
$\uA$ into the defect space $\CD_C$ of $\uC$ (Theorem 4.2.1). The
connection to Section 4.1 lies in the fact that coisometric liftings
by $*-$stable $\uA$ are subisometric (Proposition 4.2.3). But this is
only a special case and we have to generalize further.

This is done in Section 4.3. We get a hint from a result about
contractive liftings for single contractions. Lemma 2.1 in Chap.IV
of \cite{FF90} states that $E = \left(
\begin{array}{cc}
C & 0 \\
B & A \\
\end{array}
\right)$ is a contraction if and only if $C$ and $A$ are
contractions and there exists a contraction $\gamma: \CD_{*,A}
\rightarrow \CD_C$ such that
\begin{align} \label{0.3}
B = D_{*,A}\, \gamma^* \, D_C,
\end{align}
where $D_{*,A}$ and $D_C$ are the defect operators of $A^*$ and $C$.
We establish an analogous result for row contractions (Proposition
4.3.1). This shows that the isometry $\gamma$ occurring for
coisometric liftings in Section 4.2 has to be replaced in a more
general setting by a contraction.

The most general situation where we can establish a satisfactory
theory of characteristic functions for liftings is identified in
Section 4.3 and we call such liftings reduced. The technical tool here
is to use the Wold decomposition for the mid's. For $\gamma$ we
isolate the special property needed and call it resolving. Reduced
liftings include subisometric liftings as well as coisometric
liftings by c.n.c. row contractions. We define characteristic
functions for reduced liftings (Definition 4.3.6) and we argue that
this is the most general setting which is natural for that. These
characteristic functions are multi-analytic operators (not inner in
general) and they characterize reduced liftings up to unitary
equivalence. They also provide a kind of functional model for the
lifting which is useful for a closer investigation of the structure
of the lifting in the same sense as the characteristic functions of
Sz.-Nagy/Foias and of Popescu are useful in their context.

In Section 4.4 we study some further properties of these
characteristic functions. In particular we clarify the connection to
Popescu's characteristic functions and we investigate iterated
liftings, showing a factorization result for our characteristic
functions (Theorem 4.4.1). This is another indication that our
definition leads to a promising theory.

We believe that in particular the theory of subisometric liftings
may be even more interesting for row contractions than it is for
single contractions. There is a straightforward way to transfer
results from a row contraction $\uC = (C_1,\ldots,C_d)$ to the
completely positive map $\Phi_C: X \mapsto \sum^d_{i=1} C_i X
C^*_i$. This topic is taken up in Section 4.5. We define
characteristic functions for liftings of completely positive maps
and show in which way they are characteristic in this case
(Corollary 4.5.2). We investigate what subisometric lifting means in
this context and prove a one-to-one correspondence between the fixed
point sets (Theorem 4.5.4). In particular we consider the situation
where a normal invariant state is restricted to its support
(Corollary 4.5.6). From our point of view these applications give a
strong motivation for further developing the theory of liftings for
row contractions.

In an Appendix we reprove a commutant lifting theorem by
O.\,Bratteli, P.\,Jorgensen, A.\,Kishimoto and R.F.\,Werner
\cite{BJKW00}, used in Section 4.5, in a way that helps to understand
its role in our theory.

\section{Subisometric liftings}

In this section we define subisometric liftings in the setting of
row contractions and show that there is a nice Beurling-type
classification for them.
\\

We recall the notion of a minimal isometric dilation for a row
contraction, cf. \cite{Po89a}. Let $\uT=(T_1,\cdots,T_d)$ be a row
contraction on a Hilbert space $\CH$. Treating $\uT$ as an operator
from $ \bigoplus^d_{i=1} \CH$ to $\CH$, define $D_*:= (\eins-\uT
\uT^*)^\frac{1}{2}: \CH \rightarrow \CH$ and
$D:=(\eins-\uT^*\uT)^\frac{1}{2}: \bigoplus^d_{i=1} \CH \rightarrow
\bigoplus^d_{i=1} \CH$. This implies that
\begin{align}\label{1.1}
D_*= (\eins -\sum^d_{i=1} T_i T^*_i)^\frac{1}{2}, ~~~D=(\delta_{ij}
\eins -T_i^*T_j)^\frac{1}{2}_{d\times d}
\end{align}
Let $\CD_*:= \overline{\mbox{range~} D_*}$ and
$\CD:=\overline{\mbox{range~} D}$.

We  use the following {\em multi-index notation}. Let $\Lambda $
denote the set $\{ 1, 2, \ldots , d\}$ and
$\tilde{\Lambda}:=\cup_{n=0}^{\infty} \Lambda^n$, where $\Lambda
^0:=\{ 0\}$. If $\alpha \in \Lambda^n \subset \tilde{\Lambda}$ the
integer $n = |\alpha|$ is called its length. Now $T_\alpha$ with
$\alpha = (\alpha_1, \cdots ,\alpha_n) \in \Lambda^n$ means
$T_{\alpha_1} T_{\alpha_2} \ldots T_{\alpha_n}$.

The full Fock space over $\C^d$ ($d\geq 2$) denoted by $\Gamma
(\C^d)$ is
\begin{align}\label{1.2}
\Gamma (\C^d) := \mathbb{C}\oplus \C^d \oplus (\C^d)^{\otimes
^2}\oplus \cdots \oplus (\C^d) ^{\otimes ^m}\oplus \cdots.
\end{align}
To simplify notation we shall often only write $\Gamma$ instead of
$\Gamma (\C^d)$. The vector $e_0 := 1\oplus 0\oplus \cdots$ is
called the vacuum vector. Let $e_1, \ldots , e_d$ be the standard
orthonormal basis of $\C^d$. We include $d=\infty$ in which case
$\C^d$ stands for a complex separable Hilbert space of infinite
dimension. For $\alpha \in \Lambda^n$, $e_{\alpha}$ will denote the
vector $e_{\alpha _1}\otimes e_{\alpha _2}\otimes \cdots \otimes
e_{\alpha _n}$ in the full Fock space $\Gamma$. Then $e_{\alpha}$
over all $\alpha \in \tilde{\Lambda}$ forms an orthonormal basis of
the full Fock space. The (left) creation operators $L_i$ on
$\Gamma({\mathbb{C}^d})$ are defined by $ L_i x = e_i \otimes x $
for $1 \leq i \leq d$ and $x \in \Gamma({\mathbb{C}}^d).$ Then
$\underline{L}= (L_1, \ldots , L_d)$ is a row isometry, i.e., the
$L_i$ are isometries with orthogonal ranges.
\\

Using the definition of lifting in the introduction a minimal
isometric dilation (mid for short) can be described as an isometric
lifting $\uV$ of $\uT$ such that the spaces $V_\alpha \CH$ with
$\alpha \in \tilde{\Lambda}$ together span the Hilbert space on
which the $V_i$ are defined. It is an important fact, which we shall
use repeatedly, that such minimal isometric dilations are unique up
to unitary equivalence (cf. \cite{Po89a}). A useful model for the
mid is given by a version of the Sch\"{a}ffer construction, given in
\cite{Po89a}. Namely, we can realize a mid $\uV$ of $\uT$ on the
Hilbert space $\hat{\CH} := \CH  \oplus (\Gamma \otimes \CD)$,
\begin{align}\label{1.3}
V_i(h \oplus \sum_{\alpha \in \tilde{\Lambda}}  e_\alpha \otimes
d_\alpha) = T_i h \oplus [e_0 \otimes D_i h + e_i \otimes
\sum_{\alpha \in \tilde{\Lambda}} e_\alpha \otimes d_\alpha]
\end{align}
for $h \in \CH$ and $d_\alpha \in \CD.$ Here $D_i h := D
(0,\ldots,0,h,0,\ldots,0)$ and $h$ is embedded at the $i^{th}$
component.

If we have more than one row contraction at the same time then we
shall use the above notations with superscripts or subscripts, as
convenient. We are now ready for the basic definition in this
section.
\begin{defn}
Let $\uC=(C_1,\cdots,C_d)$ be a row contraction on a Hilbert space
$\CH_C$. A lifting $\underline{E}$ of $\underline{C}$ on $\CH_E
\supset \CH_C$ is called subisometric if the corresponding mids
$\underline{V}^E$ (on the Hilbert space $\hat{\CH}_E$) and
$\underline{V}^C$ (on the Hilbert space $\hat{\CH}_C$) are unitarily
equivalent, in the sense that there exists a unitary $W: \hat{\CH}_E
\rightarrow \hat{\CH}_C$ such that $W |_{\CH_C} = \eins |_{\CH_C}$
and $W V_i^E = V_i^C W$.
\end{defn}
For $d=1$ this is consistent with the definition of subisometric
dilation in \cite{DF84}, see the discussion in the introduction.
Note that the mid $\underline{V}^C$ is an example of a subisometric
lifting in this sense. Another (trivial) example is $\uC$ itself
(considered as a lifting of $\uC$). Further note that, given the
mids $\underline{V}^E$ and $\underline{V}^C$, the unitary $W$ is
uniquely determined by its properties (use the minimality of
$\underline{V}^C$).
\\

We want to make the structure of subisometric liftings more
explicit. Let $\underline{E} = (E_1,\ldots,E_d)$ be a subisometric
lifting of $\uC=(C_1,\cdots,C_d)$ on ${\CH_E} = {\CH}_C \oplus
{\CH}_A$ as in Definition 4.1.1, so that for all $i=1,\ldots,d$ we
have block matrices
\begin{align}\label{1.4}
E_i =\left(
\begin{array}{cc}
C_i  & 0 \\
B_i  & A_i \\
\end{array}
\right)
\end{align}
Let $\underline{V}^C$ be the mid of $\underline{C}$, realized as in
\eqref{1.3} on the space $\hat{\CH}_C = \CH_C  \oplus (\Gamma
\otimes \CD_C)$. Because ${\CH_E} = {\CH}_C \oplus {\CH}_A \subset
\hat{\CH}_E$ we can use the unitary $W$ from the subisometric
lifting property to obtain a subspace $\CH_{A*} := W \CH_A \subset
\Gamma \otimes \CD_C$. Further $\CH_{E*} := {\CH}_C \oplus
{\CH}_{A*} \subset \hat{\CH}_C$, and $\underline{V}^C$ is also a mid
of the row contraction $\uE_*$ which is transferred by $W$ from the
unitarily equivalent original $\uE$. We can write
\begin{align}\label{1.5}
E_{i*} =\left(
\begin{array}{cc}
C_i  & 0 \\
B_{i*} & A_{i*} \\
\end{array}
\right)
\end{align}
so $\uE_*$ is also a lifting of $\uC$.

Because $\underline{V}^C$ is a mid of $\uE_*$ it follows that
$\CH_{E*}$ is coinvariant for $\uV^C$ (by which we mean that it is
invariant for all $(V^C_i)^*$, $i=1,\ldots,d$). Note that
\begin{align}\label{1.6}
V_i^C \;|_{\Gamma \otimes \CD_C} = L_i \otimes \eins.
\end{align}
Hence $\underline{L \otimes \eins}$ is an isometric lifting of
$\uA_*$, in particular $\CH_{A*}$ is coinvariant for $\underline{L
\otimes \eins}$. An isometric lifting always contains the mid. In
particular the mid of $\uA_*$ lives on the space
$\overline{span}\{(L_\alpha \otimes \eins) \CH_{A*}, \alpha \in
\tilde{\Lambda} \}$. This subspace is reducing for the $L_i \otimes
\eins$ for all $i=1,\ldots,d$ and hence has the form $\Gamma \otimes
{\cal E}$ for a subspace ${\cal E}$ of $\CD_C$, see for example
Cor.1.7 of \cite{Po05}, where it is done in a more general setting.
In this reference the space ${\cal E}$ is described as the closure
of the image of $\CH_{A*}$ under the orthogonal projection onto $e_0
\otimes \CD_C$.

We can obtain a more concrete formula for ${\cal E}$ by comparing
this result with another way of writing the mid. First note that, as
a compression of $\underline{L \otimes \eins}$, the row contraction
$\uA_*$ (and hence also $\uA$) is $*$-stable, i.e., for all $h \in
\CH_A$
\begin{align}\label{1.7}
\lim_{n\to\infty} \sum_{|\alpha|=n} \| A^*_\alpha h\|^2 = 0\;,
\end{align}
cf. \cite{Po89a}, Prop.2.3 (where it is called pure). In this case,
with $D_{*,A} \;=\; (\eins- \uA \uA^*)^{\frac{1}{2}}:\; \CH_A
\rightarrow \CH_A$ and $\CD_{*,A}$ its closed range, the map
\begin{align}\label{1.8}
\CH_A \to \Gamma \otimes \CD_{*,A}
\\
h \mapsto \sum_{\alpha \in \tilde{\Lambda}} e_\alpha \otimes D_{*,A}
A^*_\alpha h \nonumber
\end{align}
is isometric (Popescu's Poisson kernel, cf. \cite{Po03}). With this
embedding of $\CH_A$ it can be checked that now $\underline{L
\otimes \eins}$ on $\Gamma \otimes \CD_{*,A}$ is a mid of $\uA$.

Because mids are unique up to unitary equivalence we have a unitary
$u: \Gamma \otimes \CD_{*,A} \rightarrow \Gamma \otimes {\cal E}$
such that $u \CH_A = \CH_{A*}$ and $u (L_i \otimes \eins) = (L_i
\otimes \eins) u$ for all $i=1,\ldots,d$. The commutation relation
implies that $u$ is of the form $\eins \otimes u'$, where $u'$ is a
unitary from $\CD_{*,A}$ onto ${\cal E}$ (you may use the fact that
$e_0 \otimes \CD_{*,A}$ respectively $e_0 \otimes {\cal E}$ are the
uniquely determined wandering subspaces). Thinking of $u'$ as an
isometry from $\CD_{*,A}$ into $\CD_C$ we call it $\gamma$. So
$\gamma: \CD_{*,A} \rightarrow \CD_C$ has ${\cal E}$ as its range
and it is canonically associated to a subisometric lifting in the
way shown above.

Using $\gamma$ we see that the embedding of $\CH_{A}$ into $\Gamma
\otimes {\CD_C}$ is  automatically of Poisson kernel type
\eqref{1.8}, namely
\begin{align}\label{1.9}
\CH_A \ni h \mapsto \sum_{\alpha \in \tilde{\Lambda}} e_\alpha
\otimes \gamma D_{*,A} A^*_\alpha h \in \Gamma \otimes {\CD_C}
\end{align}
which is an explicit formula for the embedding $W |_{\CH_A}: \CH_A
\rightarrow \CH_{A*} \subset \Gamma \otimes \CD_C$.

Note also that the isometry $\gamma$ is closely related to the
$\uB$-part of the lifting $\uE$. In fact, because $E^*_{i*} =
(V^C_i)^* |_{\CH_{E*}}$ we obtain $B^*_{i*} = p_C (V^C_i)^* p_{A*}$,
where $p_C, p_{A*}$ are the orthogonal projections onto $\CH_C,
\CH_{A*}$. Combining this with \eqref{1.3} and \eqref{1.9} yields
$B^*_i = D^*_{i,C}\, \gamma\, D_{*,A}: \CH_A \rightarrow \CH_C,\;
i=1,\ldots,d$. Or in a more compact form
\begin{align}\label{1.10}
\uB^* = D^*_C \, \gamma \, D_{*,A}.
\end{align}

\begin{prop}
A lifting $\uE$ of a row contraction $\uC$ with
\[
E_i =\left(
\begin{array}{cc}
C_i  & 0 \\
B_i  & A_i \\
\end{array}
\right), \qquad i=1,\ldots,d,
\]
is subisometric if and only if $\uA$ is $*$-stable and $\uB =
D_{*,A} \gamma^* D_C$ with an isometry $\gamma: \CD_{*,A}
\rightarrow \CD_C$.
\end{prop}

\begin{proof}
We have already seen above that if $\uE$ is subisometric then the
conditions are satisfied. Conversely, if $\uA$ is $*$-stable then
use the isometry $\gamma$ to embed $\uA$ (as $\uA_*$) and its mid
into $\Gamma \otimes \CD_C$ as in \eqref{1.9}. Then the formula for
$\uB$ (or \eqref{1.10}) combined with \eqref{1.3} for $\uC$ shows
that $\uV^C$ is a mid for $\uE_*$ which is unitarily equivalent to
$\uE$. (Clearly $\uV^C$ is minimal for $\uE_*$ because it is already
minimal for $\uC$.) Hence $\uE$ is subisometric.
\end{proof}

\begin{rem}
This is consistent with the results for $d=1$ in \cite{DF84} which
we mentioned in the introduction. $\gamma$ unitary corresponds to
what Douglas and Foias call a minimal subisometric dilation. We have
no reason for imposing this condition and continue to consider
general subisometric liftings. Compare also Chapter 8.3 of
\cite{Ber88}.
\end{rem}

Classifying subisometric liftings becomes especially transparent by
focusing on the invariant subspace associated to it.

\begin{defn}
Let $\uE$ on ${\CH_E} = {\CH}_C \oplus {\CH}_A$ be a subisometric
lifting of $\uC$ on ${\CH}_C$, notation as in Definition 4.1.1. Then
we call
\begin{align}\label{1.11}
{\cal N} := (\Gamma \otimes \CD_C) \ominus W \CH_A
\end{align}
the invariant subspace associated to the subisometric lifting.
Clearly ${\cal N}$ is invariant for $L_i \otimes \eins$,
$i=1,\ldots,d$.
\end{defn}

We can go the way back. Let $\uC$ on ${\cal H}_C$ be a row
contraction. If ${\cal N} \subset \Gamma \otimes \CD_C$ is a
subspace which is invariant for all $L_i \otimes \eins, \;
i=1,\ldots,d$ then we can define
\begin{align}\label{1.12,1.13}
\CH_{A*} := (\Gamma \otimes \CD_C) \ominus \CN
\\
\CH_* := \CH_C \oplus \CH_{A*}
\end{align}
On $\CH_C \oplus (\Gamma \otimes \CD_C)$ we have the mid $\uV^C$ of
$\uC$, as in \eqref{1.3}, so we can further define
\begin{align}\label{1.14}
\uE_* = (E_{1*},\ldots,E_{d*}),\quad E_{i*} := P_{\CH_*} V_i^C
|_{\CH_{E*}}: \CH_{E*} \rightarrow \CH_{E*}
\end{align}
Then $\uE_*$ is a row contraction and
\begin{align}\label{1.15}
E_{i*} =\left(
\begin{array}{cc}
C_i  & 0 \\
B_{i*}  & A_{i*} \\
\end{array}
\right)
\end{align}
with respect to the decomposition $\CH_{E*} := \CH_C \oplus
\CH_{A*}$, i.e., $\uE_*$ is a lifting of $\uC$. Then $\uV^C$ is a
mid of $\uE_*$ (minimal because it is already minimal for $\uC$).
Hence we have constructed a subisometric lifting. We are back in the
setting of Proposition 4.1.2.
\\

These considerations suggest a classification of subisometric
liftings along a Beurling type theorem for the associated invariant
subspaces. It is instructive to introduce the generalized inner
functions occurring here directly from the definition of
subisometric lifting.

So let $\uE$ be a subisometric lifting of $\uC$. Then the mids
$\uV^E$ of $\uE$ and $\uV^C$ of $\uC$ are connected by the unitary
\begin{align}\label{1.16}
W: \hat{\CH}_E = \CH_E \oplus (\Gamma \otimes \CD_E) \rightarrow
\hat{\CH}_C = \CH_C \oplus (\Gamma \otimes \CD_C)
\end{align}
such that $W |_{\CH_C} = \eins |_{\CH_C}$ and $W V^E_i = V_i^C W$
for $i=1,\ldots,d$. If we define the isometry
\begin{align}\label{1.17}
M_{C,E} := W |_{\Gamma \otimes \CD_E}
\end{align}
then from \eqref{1.3} and \eqref{1.16} we obtain
\begin{align}\label{1.18}
M_{C,E} (L_i \otimes \eins_E) = (L_i \otimes \eins_C) M_{C,E}
\end{align}
which means that $M_{C,E}: \Gamma \otimes \CD_E \rightarrow \Gamma
\otimes \CD_C$ is a {\em multi-analytic inner operator} determined
by its {\em symbol}
\begin{align}\label{1.19}
\Theta_{C,E}: \CD_E \rightarrow \Gamma \otimes \CD_C, \quad
\Theta_{C,E} = W |_{e_0 \otimes \CD_E}.
\end{align}
according to the terminology introduced in \cite{Po89b}. Obviously
this is nothing but the multi-analytic inner operator corresponding
to the invariant subspace $\CN$, in fact it is easy to check that
\begin{align}\label{1.20}
\CN = M_{C,E} (\Gamma \otimes \CD_E),
\end{align}
compare the Beurling type theorem in \cite{Po89b}. Our new insight
is that it is connected to the subisometric lifting $\uE$ of $\uC$.

\begin{defn}
We call $M_{C,E}$ (or $\Theta_{C,E}$) the characteristic function of
the subisometric lifting $\uE$ of $\uC$.
\end{defn}

It is not difficult to check that two multi-analytic inner operators
$M: \Gamma \otimes \CD \rightarrow \Gamma \otimes {\cal E}$ and
$M^\prime: \Gamma \otimes \CD^\prime \rightarrow \Gamma \otimes
{\cal E}$ with symbols $\Theta, \Theta^\prime$ describe the same
invariant subspace if and only if there exists a unitary $v: \CD
\rightarrow \CD^\prime$ such that $\Theta = \Theta^\prime v$. Let us
call multi-analytic functions {\em equivalent} if they are related
in this way. We are ready for our classification result.

\begin{thm}
Let $\uC = (C_1,\ldots,C_d)$ be a row contraction on a Hilbert space
$\CH_C$. Then there is a one-to-one correspondence between
\begin{itemize}
\item[(a)]
unitary equivalence classes of subisometric liftings $\uE$ of $\uC$,
\item[(b)]
$\underline{L \otimes \eins}$-invariant subspaces $\CN$ of $\Gamma
\otimes \CD_C$,
\item[(c)]
multi-analytic inner operators $M$ with symbols $\Theta: \CD
\rightarrow \Gamma \otimes \CD_C$ up to equivalence.
\end{itemize}
The correspondence is described above. In particular if $\uE$ is the
lifting then $\CD = \CD_E$, $M = M_{C,E}$ with symbol $\Theta =
\Theta_{C,E}$ and $(b) \leftrightarrow (c)$ is Beurling's theorem.
\end{thm}

Theorem 4.1.6 shows that the characteristic function of a subisometric
lifting characterizes the lifting up to unitary equivalence,
justifying to call it characteristic.

\begin{proof}
$(b) \leftrightarrow (c)$ is Beurling's theorem, see \cite{Po89b}.
We now show that the correspondence $(a) \rightarrow (c)$ is well
defined. Let $\uE$ on $\CH_E \supset \CH_C$ and $\uE^\prime$ on
$\CH_{E^\prime} \supset \CH_C$ be two subisometric liftings of $\uC$
which are unitarily equivalent, i.e., there exists a unitary $u:
\CH_E \rightarrow \CH_{E^\prime}$ such that $u |_{\CH_C} = \eins
|_{\CH_C}$ and $E^\prime_i u = u E_i$ for $i=1,\ldots,d$. Clearly
unitarily equivalent row contractions have unitarily equivalent mids
and we can extend $u$ (in a trivial way) to a unitary $\hat{u}$
between the spaces $\hat{\CH}_E$ and $\hat{\CH}_{E^\prime}$ of the
mids $\uV^E$ and $\uV^{E^\prime}$, so we have
\begin{align*}
\hat{u}: \hat{\CH}_E \rightarrow \hat{\CH}_{E^\prime}
\;\mbox{unitary}, \quad \hat{u} |_{\CH_E} = u, \quad V^{E^\prime}_i
\hat{u} = \hat{u} V^E_i \quad(i=1,\ldots,d)
\end{align*}
Because $\uE,\uE^\prime$ are subisometric we also have unitaries
$W,W^\prime$ such that
\begin{align*}
W: \hat{\CH}_E \rightarrow \hat{\CH}_C, \quad V^C_i W = W V^E_i,
\quad W |_{\CH_C} = \eins |_{\CH_C}
\\
W^\prime: \hat{\CH}_{E^\prime} \rightarrow \hat{\CH}_C, \quad V^C_i
W^\prime = W^\prime V^{E^\prime}_i, \quad W^\prime |_{\CH_C} = \eins
|_{\CH_C} \nonumber
\end{align*}
If we now define
\[
u_C := W^\prime \hat{u} W^*: \hat{\CH}_C \rightarrow \hat{\CH}_C
\]
then it follows that $u_C$ commutes with the $V^C_i$ for
$i=1,\ldots,d$. To see that, ``chase" the following commuting
diagram
\begin{align}\label{1.21}
\xymatrix@!{
 & \hat{\CH}_C \ar[rr]^{u_C} \ar'[d][dd]_{V^C_i}
   & & \hat{\CH}_C \ar[dd] \ar'[d][dd]_{V^C_i}
\\
 \hat{\CH}_E \ar[ur]^{W}
\ar[rr]^{\quad\quad\quad\quad\hat{u}} \ar'[d]^{V^E_i}[dd] \ar[dd]
 & & \hat{\CH}_{E'} \ar[ur]^{W^\prime} \ar'[d]^{V^{E^\prime}_i}[dd] \ar[dd] &
\\
 & \hat{\CH}_C \ar'[r]_{u_C}[rr]
   & & \hat{\CH}_C
\\
   \hat{\CH}_E \ar[rr]_{\hat{u}} \ar[ur]_W
 & & \hat{\CH}_{E'} \ar[ur]_{W^\prime}
}
\end{align}
Further, because $W, W^\prime$ and $\hat{u}$ all fix $\CH_C$
pointwise the same is true for $u_C$, so we have also $u_C |_{\CH_C}
= \eins |_{\CH_C}$. But by minimality of $\uV^C$ we know that
$\hat{\CH}_C$ is the closed linear span of vectors of the form
$V^C_\alpha h$ with $\alpha \in \tilde{\Lambda},\,h \in \CH_C$ and
from
\[
u_C V^C_\alpha h = V^C_\alpha u_C h = V^C_\alpha h
\]
we infer that $u_C = \eins$. Hence $W = (u_C)^* W^\prime \hat{u} =
W^\prime \hat{u}$. Clearly $\hat{u}$ maps $e_0 \otimes \CD_E \subset
\hat{\CH}_E$ onto $e_0 \otimes \CD_{E^\prime} \subset
\hat{\CH}_{E^\prime}$, so if we define the unitary $v := \hat{u}
|_{\CD_E}: \CD_E \rightarrow \CD_{E^\prime}$ and use that $\Theta =
W |_{\CD_E}$ and $\Theta^\prime = W^\prime |_{\CD_{E^\prime}}$ we
see that $\Theta = \Theta^\prime v$, i.e., the characteristic
functions are equivalent.

Conversely suppose that we are given a multi-analytic inner operator
with symbol $\Theta: \CD \rightarrow \Gamma \otimes \CD_C$, as in
(c). By $(b) \leftrightarrow (c)$ (Beurling's theorem) we have an
invariant subspace $\CN$ which is associated to a subisometric
lifting $\uE$ of $\uC$ and $\CD = \CD_E$, see the discussion
preceding the theorem. It remains to show that if $\Theta =
\Theta^\prime v$ with a unitary $v: \CD_E \rightarrow
\CD_{E^\prime}$ for two subisometric liftings $\uE$ and $\uE^\prime$
then $\uE$ and $\uE^\prime$ are unitarily equivalent. Let
$W,W^\prime$ be the corresponding unitaries from the subisometric
lifting property. Then
\begin{align*}
W^\prime \CH_{E^\prime} \;=\;  \CH_C \oplus ( \Gamma \otimes \CD_C )
\;\ominus\; W^\prime ( \Gamma \otimes \CD_{E^\prime} ) \\
= \CH_C \oplus ( \Gamma \otimes \CD_C )
\;\ominus\; M_{C,E^\prime} \big( \Gamma \otimes v \CD_E \big) \\
= \CH_C \oplus ( \Gamma \otimes \CD_C ) \;\ominus\; M_{C,E} ( \Gamma
\otimes \CD_E ) \;=\; W \CH_E,
\end{align*}
and we can define
\begin{align*}
U := (W^\prime)^* \, W |_{\CH_E}: \CH_E \rightarrow \CH_{E^\prime}.
\end{align*}

Because for $h \in \CH_C,$ $W h = h = W^\prime h$ we have $U h = h$.
In general for $h\in \CH_E$ and $i=1,\ldots,d$ (with $p_E,
p_{E^\prime}$ orthogonal projections onto $\CH_E, \CH_{E^\prime}$)
\begin{align*}
U E_i\, h = (W^\prime)^*\, W \,E_i \,h = (W^\prime)^* \,W \,p_E
\,V^E_i \,h = p_{E^\prime}\, (W^\prime)^*\, W\, V^E_i \,h
\\
= p_{E^\prime}\, (W^\prime)^*\, V^C_i\, W\, h = p_{E^\prime}\,
V^{E^\prime}_i\, (W^\prime)^* \,W\, h = E^\prime_i \,U\, h,
\end{align*}
i.e., $\uE$ and $\uE^\prime$ are unitarily equivalent.
\end{proof}

There is an interesting variant of the classification if we not only
give $\uC$ but also $\uA$, i.e., if we consider liftings of $\uC$ by
$\uA$.

\begin{thm}
Let $\uA$ and $\uC$ be row contractions, $\uA$ $*$-stable. There is
a one-to-one correspondence between
\begin{itemize}
\item[(a)]
unitary equivalence classes of subisometric liftings of $\uC$ by
$\uA$
\item[(b)]
equivalence classes of isometries $\gamma: \CD_{*,A} \rightarrow
\CD_C$, two isometries considered equivalent if they have the same
range
\end{itemize}
\end{thm}

\begin{proof}
The details of this correspondence have already been discussed in
connection with Proposition 4.1.2. It is shown there how to construct
an isometry $\gamma: \CD_{*,A} \rightarrow \CD_C$ if a subisometric
lifting of $\uC$ by $\uA$ is given, and conversely, how to use such
an isometry to find a subisometric lifting. The equivalence in (b)
is chosen in such a way that two isometries are equivalent if and
only if the associated invariant subspaces are the same, compare
(4.12) and (4.14). Hence the result follows from Theorem 4.1.6.
\end{proof}

\begin{cor}
Let $\uA$ and $\uC$ be row contractions, $\uA$ $*$-stable. A
subisometric lifting of $\uC$ by $\uA$ exists if and only if
\begin{align*}
dim \CD_{*,A} \leq dim \CD_C,
\end{align*}
where $dim$ stands for the cardinality of an orthonormal basis. In
the case $dim \CD_{*,A} = dim \CD_C$ (minimal subisometric dilation
in the terminology of \cite{DF84}) the lifting is unique up to
unitary equivalence.
\end{cor}

\section{Coisometric liftings}

The theory of subisometric liftings turns out to be especially
relevant in the case of coisometric row contractions and coisometric
liftings. We start with definitions and elementary properties.

A row contraction $\uC$ on $\CH_1$ is called {\it coisometric} if
$\uC \,\uC^* = \eins$, i.e., $\sum^d_{i=1} C_i C^*_i = \eins$. It is
easy to check that a lifting $\uE$ on ${\cal H} = {\cal H}_C \oplus
{\cal H}_A$ with block matrices
\begin{align*}
E_i =\left(
\begin{array}{cc}
C_i  & 0 \\
B_i  & A_i \\
\end{array}
\right)
\end{align*}
(for all $i=1,\ldots,d$) is coisometric if and only if
$\underline{C}$ is coisometric and
\begin{align}\label{2.1}
\uB\, \uC^* = 0,\quad \mbox{i.e.},\; \sum^d_{i=1} B_i C^*_i = 0,
\\
\uA\, \uA^* + \uB\, \uB^* = \eins,\quad \mbox{i.e.},\; \sum^d_{i=1}
A_i A^*_i + \sum^d_{i=1} B_i B^*_i = \eins.
\end{align}

\begin{thm}
Let $\uA$ and $\uC$ be row contractions, $\uC$ coisometric. Then
there is a one-to-one correspondence between
\begin{itemize}
\item[(a)]
coisometric liftings $\uE$ of $\uC$ by $\uA$
\item[(b)]
isometries $\gamma: \CD_{*,A} \rightarrow \CD_C$
\end{itemize}

Explicitly, if $E_i = \left(
\begin{array}{cc}
C_i  & 0 \\
B_i  & A_i \\
\end{array}
\right)$ for $i=1,\ldots,d$ provides a coisometric lifting $\uE$ of
$\uC$ by $\uA$ then $\gamma: D_{*,A} h \mapsto \uB^* h \subset
\CD_C$ (for $h \in \CH_A$) is isometric.
\\

Conversely, if $\gamma: \CD_{*,A} \rightarrow \CD_C$ is isometric
then with $\uB^* := \gamma D_{*,A}$ we obtain a coisometric lifting
$\uE$ by $E_i = \left(
\begin{array}{cc}
C_i  & 0 \\
B_i  & A_i \\
\end{array}
\right)$ for $i=1,\ldots,d$.
\end{thm}

\begin{proof}
Because $\uC$ is coisometric, $D_C=\eins -\uC^* \uC$ is the
orthogonal projection onto the kernel of $\uC$.

Let $\uE$ be a coisometric lifting of $\uC$ by $\uA$. Then from
(4.25) we have $\uC \uB^* = (\uB \uC^*)^* =0 $ and hence
$range(\uB^*) \subset \CD_C$.

Further for $h \in \CH_A$, using (4.26)
\begin{align*}
\| D_{*,A} h \|^2 = \langle (\eins - \uA\uA^*)h, h \rangle =\langle
\uB\uB^*h, h \rangle = \| \uB^*h\|^2
\end{align*}
So there exist an isometry $\gamma : \CD_{*,A} \to range(\uB^*)
\subset \CD_C$ with $\gamma D_{*,A} h = \uB^*h$ for all $h \in
\CH_A$.

Conversely, let $\gamma: \CD_{*,A} \rightarrow \CD_C$ be an isometry
and define $\uB^* := \gamma D_{*,A}$. From $\uC |_{\CD_C} = 0$ we
obtain $\uC\, \uB^* = 0$ or $\uB \,\uC^*=0 $, which is (4.25).
Further
\begin{align*}
\uB \,\uB^* = D_{*,A} \gamma^* \gamma D_{*,A} = D^2_{*,A} = \eins -
\uA\, \uA^*,
\end{align*}
hence $\uA\, \uA^* + \uB\, \uB^* = \eins$, which is (4.26). Hence
with $E_i = \left(
\begin{array}{cc}
C_i  & 0 \\
B_i  & A_i \\
\end{array}
\right),$ for $i=1,\ldots,d,$ we obtain a coisometric lifting $\uE$
of $\uC$ by $\uA$.

Finally if $\gamma, \gamma^\prime: \CD_{*,A} \rightarrow \CD_C$ are
two isometries and $\gamma \not= \gamma^\prime$ then $\uB^* \not=
(\uB^\prime)^*$ for $\uB^*,(\uB^\prime)^*$ defined by $\gamma,
\gamma^\prime$ as above. Hence the correspondence is one-to-one.
\end{proof}

\begin{cor}
Let $\uA$ and $\uC$ be row contractions, $\uC$ coisometric. A
coisometric lifting $\uE$ of $\uC$ by $\uA$ exists if and only if
\begin{align*}
dim \CD_{*,A} \leq dim \CD_C,
\end{align*}
where $dim$ stands for the cardinality of an orthonormal basis..
\end{cor}

Theorem 4.2.1 gives a kind of free parametrization of the coisometric
liftings. Let us consider an elementary example.
\begin{align} \label{2.3}
\uc = (c_1,\ldots,c_d) \in \C^d, \quad \|\uc\|^2 = \sum^d_{i=1}
|c_i|^2 = 1 \quad {\mbox (unit\, sphere)}
\\
\ua = (a_1,\ldots,a_d) \in \C^d, \quad \|\ua\|^2 = \sum^d_{i=1}
|a_i|^2 \leq 1 \qquad {\mbox (unit\, ball)} \nonumber
\end{align}
Then we get a left lower corner $\ub = (b_1,\ldots,b_d)$ for a
coisometric lifting if $\langle \ub,\uc \rangle = 0$ and $\|\ua\|^2
+ \|\ub\|^2 = 1$, according to (4.25) and (4.26). Obviously the set of
solutions for $\ub$ is the (complex) sphere with radius $r =
\sqrt{1-\|\ua\|^2}$ in the subspace orthogonal to $\uc$. If
$\|\ua\|=1$ the solution is unique. We can check that the
parametrization using isometries $\gamma: \CD_{*,A} \rightarrow
\CD_C$ as in Theorem 4.2.1 yields the same result.
\\

Theorem 4.2.1 and Corollary 4.2.2 are even true if $\uA$ is not
$*$-stable. If $\uA$ is $*$-stable then we should compare these
results with those in Section 4.1. Note in particular that the formula
$\uB^* = \gamma D_{*,A}$ in Theorem 4.2.1 and the formula $\uB^* =
D^*_C \,\gamma\, D_{*,A}$ (1.10) are compatible because, as noted
above, for $\uC$ coisometric the operator $D^*_C$ is nothing but the
embedding of $\CD_C$ into $\oplus^d_{i=1} \CH_C$ which is implicit
in the formulation chosen in Theorem 4.2.1. Further comparison yields
the following result which shows that subisometric liftings occur
very naturally in the coisometric setting.

\begin{prop}
Let $\uC$ be a coisometric row contraction. A lifting of $\uC$ is a
coisometric lifting by a $*$-stable $\uA$ if and only if it is
subisometric.
\end{prop}

\begin{proof}
Using Theorem 4.2.1 we can replace the condition ``coisometric" for
the lifting by the existence of an isometry $\gamma: \CD_{*,A}
\rightarrow \CD_C$ such that $\uB^* = \gamma D_{*,A} = D^*_C \gamma
D_{*,A}$. Now Proposition 4.2.3 is a direct consequence of Proposition
4.1.2.
\end{proof}

In particular, for coisometric liftings by a $*$-stable $\uA$ there
exists an associated invariant subspace and a characteristic
function. In the special case $dim\, \CH_C = 1$ this characteristic
function was introduced in \cite{DG07a} under the name ``extended
characteristic function". For general $\CH_C$, in view of Theorem
4.1.6, it is better to call it the characteristic function of the
lifting (with $\uC$ given), as we have done in Definition 4.1.5.

\section[Reduced liftings]{Characteristic functions of reduced liftings}
In this section we generalize the theory of characteristic functions
for sub\-isometric liftings from Section 4.1 and establish a setting
that also includes the setting of Section 4.2.

Let $\uC$ be a row contraction on $\CH_C$ and $\uE$ on $\CH_E =
\CH_C \oplus \CH_A$ be a (contractive) lifting so that for all
$i=1,\ldots,d$
\begin{align*}
E_i =\left(
\begin{array}{cc}
C_i  & 0 \\
B_i  & A_i \\
\end{array}
\right)
\end{align*}
Then as in \eqref{1.3} we have a mid $\uV^E$ on $\CH_E \oplus \big(
\Gamma \otimes \CD_E \big)$. Clearly $\uV^E$ is an isometric lifting
of $\uC$, so the space of the mid $\uV^C$ can be embedded as a
subspace reducing the $V^E_i$. Let us encode this by introducing the
restriction $\uY$ on the orthogonal complement $\CK$ and a unitary
$W$ by
\begin{align}\label{3.1}
W: \CH_E \oplus \big( \Gamma \otimes \CD_E \big) \rightarrow
\CH_C \oplus \big( \Gamma \otimes \CD_C \big) \oplus \CK \\
\tilde{V}^E_i W = W V^E_i, \quad W |_{\CH_C} = \eins |_{\CH_C} \quad
\mbox{with} \quad \tilde{\uV}^E = \uV^C \oplus \uY \nonumber
\end{align}
By omitting $\CH_C$ we also have a unitary (also denoted by $W$)
\begin{align}\label{3.2}
W: \CH_A \oplus \big( \Gamma \otimes \CD_E \big) \rightarrow \big(
\Gamma \otimes \CD_C \big) \oplus \CK
\end{align}
and an isometric embedding $\CH_{A_*} := W \CH_A \subset \big(
\Gamma \otimes\CD_C \big) \oplus \CK$. Further we obtain
\begin{align}\label{3.3}
\uB^* = p_{\CH_C} (\uV^E)^* |_{\CH_A} = p_{\CH_C} \big[(\uV^C)^*
\oplus \uY^* \big] W |_{\CH_A} = D^*_C\, p_{e_0 \otimes \CD_C} W
|_{\CH_A}
\end{align}
where we used formula \eqref{1.3} for $\uV^C$.
\\

To proceed we need a few facts about the mid $\uV^A$ on
$\tilde{\CH}_A$ of $\uA$. We write its Wold decomposition as
\begin{align}\label{3.4}
\tilde{\CH}_A = \big( \Gamma \otimes \CD_{*,A} \big) \oplus \CR_A \\
V^A_i = (L_i \otimes \eins) \oplus R^A_i, \quad i=1,\ldots,d,
\nonumber
\end{align}
where $\CR_A$ and $\uR^A$ stand for the residual part (cf.
\cite{Po89a}). The embedding of $\CH_A$ into $\tilde{\CH}_A$ can be
written as
\begin{align}\label{3.5}
\CH_A \ni h \mapsto \big( \sum_{\alpha \in \tilde{\Lambda}} e_\alpha
\otimes D_{*,A} A^*_\alpha h \big) \oplus h_{\CR}
\end{align}
Here $h_{\CR}$ belongs to the residual part $\CR_A$. Compare
\cite{BDZ06} for a derivation of this decomposition via
Stinespring's theorem. In fact, it is not difficult to check that a
formula like \eqref{3.5} always reproduces the Wold decomposition
above, compare also \cite{FF90} for similar arguments for $d=1$.
Note that the residual part vanishes if and only if $\uA$ is
$*$-stable, so in this case we are back in the setting of Section 4.1.

Further we need the decomposition $\CH_A = \CH^1_A \oplus \CH^0_A$
with $\CH^1_A$ the largest subspace invariant for the $A^*_i$ and
such that the restriction of $\uA^*$ is isometric, i.e.,
\begin{align}\label{3.6}
\CH^1_A := \{ h \in \CH_A: \sum_{|\alpha|=n} \| A^*_\alpha h \|^2 =
\|h\|^2 \;\mbox{for all}\; n \in \N \}
\end{align}
Then it is easy to check that $\CH^1_A = \CH_A \cap \CR_A$ (cf.
\cite{Po89a}, Proposition 2.9), but the position of $\CH^0_A$ may be
complicated with respect to the decomposition \eqref{3.5} because
$\uA$ restricted to $\CH^0_A$ may not be $*$-stable and in this case
$\CH^0_A$ is not contained in $\Gamma \otimes \CD_{*,A}$. In fact,
if $0 \not= h \in \CH^0_A$ we only have
\begin{align}\label{3.7}
\sum_{|\alpha|=n} \| A^*_\alpha h \|^2 < \|h\|^2 \;\mbox{for some}\;
n \in \N
\end{align}
which (by definition) means that $\uA |_{\CH^0_A}$ is {\em
completely non-coisometric (c.n.c.)}, cf. \cite{Po89a}.

Now we look at $\uA$ and its mid $\uV^A$ embedded into the larger
structure obtained from the lifting $\uE$. Clearly $\uV^E$
restricted to $\CH_A \oplus \big( \Gamma \otimes \CD_E \big)$ is an
isometric dilation of $\uA$, so $\CH_A \oplus \big( \Gamma \otimes
\CD_E \big)$ contains $\tilde{\CH}_A$ as a $V^E_i$-reducing subspace
($i=1,\ldots,d$) which we still denote by $\tilde{\CH}_A$. Using
\eqref{3.2} we see that $\big( \Gamma \otimes \CD_C \big) \oplus
\CK$ contains the $(L_i \otimes \eins) \oplus Y_i$-reducing subspace
$W \tilde{\CH}_A$ and the restriction of $(\underline{L \otimes
\eins}) \oplus \uY$ is a mid of $\uA$ (transferred to $W \CH_A$).
Denoting the restriction of $W$ to $\tilde{\CH}_A$ also by $W$ we
have (for $i=1,\ldots,d$)
\begin{align}\label{3.8}
W \big[ (L_i \otimes \eins) \oplus R^A_i \big] = W V^A_i = \big[
(L_i \otimes \eins) \oplus Y_i \big] W
\end{align}
Where is $\CH_{A_*} = W \CH_A$? Clearly
\begin{align}\label{3.9}
W \CH^1_A = W (\CH_A \cap \CR_A) \subset W \CR_A \subset \CK,
\end{align}
where the last inclusion follows from \eqref{3.8} and the fact that
$\underline{L \otimes \eins}$ is $*$-stable. The position of $W
\CH^0_A$ may be more complicated.

To organize the relevant data we use \eqref{3.4} together with the
embedding of $\tilde{\CH}_A$ into $\CH_A \oplus \big( \Gamma \otimes
\CD_E \big)$ and \eqref{3.2} to define
\begin{align}\label{3.10}
M: \Gamma \otimes \CD_{*,A} \rightarrow \Gamma \otimes \CD_C, \\
M = P_{\Gamma \otimes \CD_C} W |_{\Gamma \otimes \CD_{*,A}}
\nonumber
\end{align}
which is a multi-analytic operator. Then for $h \in \CH_A$
\begin{align*}
P_{e_0 \otimes \CD_C} W h = P_{e_0 \otimes \CD_C} M P_{\Gamma
\otimes \CD_{*,A}} h = P_{e_0 \otimes \CD_C} M P_{e_0 \otimes
\CD_{*,A}} h
\end{align*}
where for the first equality we used \eqref{3.9} and the second then
follows from the fact that $M$ is multi-analytic. But $P_{e_0
\otimes \CD_{*,A}} h = e_0 \otimes D_{*,A} h$ by \eqref{3.5} and we
conclude that $P_{e_0 \otimes \CD_C} W |_{\CH_A}: \CH_A \rightarrow
\CD_C$ factors through $\CD_{*,A}$ in the sense that there exists a
contraction $\gamma := P_{e_0 \otimes \CD_C} M |_{e_0 \otimes
\CD_{*,A}}: \CD_{*,A} \rightarrow \CD_C$ such that
\begin{align}\label{3.11}
P_{e_0 \otimes \CD_C} W |_{\CH_A} = \gamma D_{*,A}
\end{align}
In fact, $\gamma$ is nothing but the the $0$-th Fourier coefficient
of $M$ in the sense of \cite{Po03}. Combined with \eqref{3.3} we
obtain
\begin{align}\label{3.12}
\uB^* = D^*_C\, \gamma\, D_{*,A}: \CH_A \rightarrow
\bigoplus^d_{i=1} \CH_C
\end{align}
This is one half of the following analogue for row contractions of
Lemma 2.1 in Chap.IV of \cite{FF90}, which already has been
discussed in the introduction, see in particular \eqref{0.3}.
\begin{prop}
$\uE = (E_1,\ldots,E_d)$ on $\CH_E = \CH_C \oplus \CH_A$ with block
matrices
\begin{align*}
E_i =\left(
\begin{array}{cc}
C_i  & 0 \\
B_i  & A_i \\
\end{array}
\right)
\end{align*}
(for $i=1,\ldots,d$) is a row contraction if and only if $\uC$ and
$\uA$ are row contractions and there exists a contraction $\gamma:
\CD_{*,A} \rightarrow \CD_C$ such that \eqref{3.12} holds.
\end{prop}

\begin{proof}
Clearly if $\uE$ is a row contraction then $\uC$ and $\uA$ are row
contractions. Above we have already given a (dilation) proof that if
$\uE$ is contractive then $\uB$ satisfies \eqref{3.12} for a
suitable contraction $\gamma$. To prove the converse, let $\gamma:
\CD_{*,A} \rightarrow \CD_C$ be a contraction and $\uB^*$ given as
in \eqref{3.12}. Then for $x \in \CH_C,\; y \in \CH_A$
\begin{align*}
| \langle x, \uC\, \uB^* y \rangle |^2 = | \langle x, \uC D^*_C
\gamma D_{*,A}\, y \rangle |^2 = | \langle D_C \uC^* x, \gamma
D_{*,A}\, y \rangle |^2
\end{align*}
\begin{align*}
\leq \| D_C \uC^* x \|^2 \| \gamma D_{*,A} y \|^2 \leq \langle x,
(\eins - \uC \uC^*) x \rangle \; \langle y, (\eins - \uA \uA^*) y
\rangle
\end{align*}
which implies (see for example Exercise 3.2 in \cite{Pau03}) that
\begin{align*}
0 \leq \left(
\begin{array}{cc}
\eins - \uC \uC^* & -\uC \uB^*\\
-\uB \uC^*  & \eins - \uA \uA^* \\
\end{array}
\right) = \eins - \uE \uE^*
\end{align*}
hence $\uE$ is a row contraction.
\end{proof}

Let us go back to the lifting $\uE$ of $\uC$ by $\uA$. The following
definition is useful to analyze further the position of $W \CH_A$.
\begin{defn}
$\gamma: \CD_{*,A} \rightarrow \CD_C$ is called resolving if for all
$h \in \CH_A$ we have
\begin{align*}
\big( \gamma D_{*,A} A^*_\alpha h =0 \;\mbox{for all}\; \alpha \in
\tilde{\Lambda} \big) \Rightarrow \big( D_{*,A} A^*_\alpha h =0
\;\mbox{for all}\; \alpha \in \tilde{\Lambda} \big)
\end{align*}
\end{defn}
Clearly if $\gamma: \CD_{*,A} \rightarrow \CD_C$ is injective then
it is resolving. Note that $D_{*,A} A^*_\alpha h =0 \;\mbox{for
all}\; \alpha \in \tilde{\Lambda}$ if and only if $h \in \CH^1_A$,
and so the intuitive meaning of `resolving' is that `\,looking at
$\CH_A$ through $\gamma$\,' still allows to detect whether $h \in
\CH_A$ is in $\CH^1_A$ or not. More precisely, $\gamma$ is resolving
if and only if for all $h \in \CH^0_A = \CH_A \ominus \CH^1_A$ there
exists $\alpha \in \tilde{\Lambda}$ such that $\gamma D_{*,A}
A^*_\alpha h \not=0$. In particular if $\uA$ is c.n.c., i.e.
$\CH^1_A = \{0\}$, then $\gamma$ is resolving if and only if for all
$0 \not= h \in \CH_A$ there exists $\alpha \in \tilde{\Lambda}$ such
that $\gamma D_{*,A} A^*_\alpha h \not=0$.
\begin{lem}
The following assertions are equivalent
\begin{itemize}
\item[(a)]
$\gamma$ is resolving.
\item[(b)]
$W \CH_A \cap \CK \subset W \CH^1_A$
\item[(c)]
$W \CH_A \cap \CK \;=\; W \CH^1_A$
\item[(d)]
$\big( \Gamma \otimes \CD_C \big) \bigvee W \big( \Gamma \otimes
\CD_E \big) \; = \; \big( \Gamma \otimes \CD_C \big) \oplus (\CK
\ominus W \CH^1_A)$
\end{itemize}
\end{lem}
\begin{proof}
(b) says that for $h \in \CH_A \setminus \CH^1_A$ the embedded $Wh$
is not in $\CK$, so not orthogonal to $\Gamma \otimes \CD_C$,
equivalently, there exists $\alpha \in \tilde{\Lambda}$ such that
\begin{align*}
0 \not= P_{e_0 \otimes \CD_C} \big[ (L^*_\alpha \otimes \eins)
\oplus Y^*_\alpha \big] W h = P_{e_0 \otimes \CD_C} W (V^A_\alpha)^*
h = \gamma D_{*,A} A^*_\alpha h
\end{align*}
(where we used the embedding of the mid of $\uA$ and in particular
\eqref{3.11}). By comparison with the comments following Definition
4.3.2 we conclude that (a) and (b) are equivalent. We noted in
\eqref{3.9} that always $W \CH^1_A \subset \CK$, so (b) and (c) are
equivalent.

To get the equivalence of (c) and (d) note that $x \in \big( \Gamma
\otimes \CD_C \big) \oplus \CK$ is orthogonal to $\big( \Gamma
\otimes \CD_C \big)$ and to $W \big( \Gamma \otimes \CD_E \big)$ if
and only if $x \in \CK$ and $x \in W \CH_A$ (compare \eqref{3.2}).
Hence the orthogonal complement of $\big( \Gamma \otimes \CD_C \big)
\bigvee W \big( \Gamma \otimes \CD_E \big)$ in $\big( \Gamma \otimes
\CD_C \big) \oplus \CK$ is in fact $W \CH_A \cap \CK$.
\end{proof}
\begin{defn}
A lifting $\uE$ of $\uC$ by $\uA$ is called reduced if $\uA$ is
c.n.c. (i.e., $\CH^1_A = \{0\}$, see \eqref{3.7}) and $\gamma$ is
resolving.
\end{defn}
We have already seen two important classes of reduced liftings.
\begin{itemize}
\item[1)]
Subisometric liftings. Here $\uA$ is $*$-stable and $\gamma$ is
isometric, see Proposition 4.1.2.
\item[2)]
Coisometric liftings by $\uA$ c.n.c. Here $\gamma$ is isometric by
Theorem 4.2.1.
\end{itemize}
Note that by Proposition 4.2.3 the coisometric liftings by $*$-stable
$\uA$ are exactly the intersection of cases 1) and 2).
\begin{lem}
The following assertions are equivalent
\begin{itemize}
\item[(a)]
$\uE$ is reduced.
\item[(b)]
$\{ h \in \CH_A: \gamma D_{*,A} A^*_\alpha h = 0 \;\mbox{for all}\;
\alpha \in \tilde{\Lambda}\} \;=\; \{0\}$
\item[(c)]
$W \CH_A \cap \CK \;=\; \{0\}$
\end{itemize}
\end{lem}
\begin{proof}
If $\gamma$ is resolving then (by definition) the space given in (b)
is contained in $\CH^1_A$. Hence (a) implies (b). Also, from (b) we
first conclude that $\CH^1_A = \{ h \in \CH_A: D_{*,A} A^*_\alpha h
= 0 \;\mbox{for all}\; \alpha \in \tilde{\Lambda}\} =\; \{0\}$ and
then that $\gamma$ is resolving, so (b) implies (a). If we have (c)
then by Lemma 4.3.3(b) $\gamma$ is resolving and then by Lemma 4.3.3(c)
$\uA$ is c.n.c., so we have (a). Given (a), Lemma 4.3.3(c) implies
(c).
\end{proof}

If $\gamma D_{*,A} A^*_\alpha h = 0 \;\mbox{for all}\; \alpha \in
\tilde{\Lambda}$ then by \eqref{3.12} we conclude that $A^*_\alpha h
\in ker \uB^* = (range \uB)^\perp$. Hence vectors in the space $\{ h
\in \CH_A: \gamma D_{*,A} A^*_\alpha h = 0 \;\mbox{for all}\; \alpha
\in \tilde{\Lambda}\}$ do not contribute in any way to the
interaction between $\CH_A$ and $\CH_C$ via $\uB^*$, and it is no
great loss to concentrate on liftings where this space has been
removed. By Lemma 4.3.5(b), in doing this we obtain exactly the
reduced liftings. This also explains our terminology.
\\

For reduced liftings we can successfully develop a theory of
characteristic functions.
\begin{defn}
Let $\uE$ be a reduced lifting of $\uC$ by $\uA$. We call the
multi-analytic operator
\begin{align}\label{3.13}
M_{C,E}: \Gamma \otimes \CD_E \rightarrow \Gamma \otimes \CD_C, \\
M_{C,E} = P_{\Gamma \otimes \CD_C} W |_{\Gamma \otimes \CD_E}
\nonumber
\end{align}
(or its symbol $\Theta_{C,E}: \CD_E \rightarrow \Gamma \otimes
\CD_C$) the characteristic function of the lifting $\uE$.
\end{defn}


Using the characteristic function we can develop a theory of
functional models for reduced liftings. The idea is similar as in
the case of characteristic functions for c.n.c. row contractions,
see \cite{Po89b}.

Let $\uE$ be a reduced lifting of $\uC$ by $\uA$. From $\uA$ c.n.c.
we obtain $\CH^1_A = \{0\}$ and then Lemma 4.3.3 gives
\begin{align} \label{3.14}
\big( \Gamma \otimes \CD_C \big) \bigvee W \big( \Gamma \otimes
\CD_E \big) \; = \; \big( \Gamma \otimes \CD_C \big) \oplus \CK
\end{align}
With the definition
\begin{align} \label{3.15}
\Delta_{C,E} := (\eins - M^*_{C,E} M_{C,E})^{\frac{1}{2}}: \Gamma
\otimes \CD_E \rightarrow \Gamma \otimes \CD_E
\end{align}
we obtain for $x \in \Gamma \otimes \CD_E$
\begin{align} \label{3.16}
\|P_{\CK} W x \|^2 = \| (\eins - P_{\Gamma \otimes \CD_C}) W x \|^2
= \|x\|^2 - \| P_{\Gamma \otimes \CD_C} W x \|^2
\\
= \|x\|^2 - \| M_{C,E} x \|^2 = \| \Delta_{C,E} x \|^2 \nonumber
\end{align}
This means that we can isometrically identify $\CK$ with
$\overline{\Delta_{C,E}(\Gamma \otimes \CD_E)}$ and with this
identification we have
\begin{align} \label{3.17}
W \CH_A = \big[ (\Gamma \otimes \CD_C) \oplus \CK  \big] \ominus W
(\Gamma \otimes \CD_E)
\end{align}
\begin{align*}
= \big[ (\Gamma \otimes \CD_C) \oplus \overline{\Delta_{C,E}(\Gamma
\otimes \CD_E)} \big] \ominus \{ M_{C,E}\, x \oplus \Delta_{C,E}\,
x: x \in \Gamma \otimes \CD_C \} \nonumber
\end{align*}
which is a kind of functional model.
\begin{thm}
Let $\uC$ be a row contraction. Reduced liftings $\uE$ and $\uE'$ of
$\uC$ are unitarily equivalent if and only if their characteristic
functions $M_{C,E}$ and $M_{C,E'}$ are equivalent.
\end{thm}
Recall that $M_{C,E}$ and $M_{C,E'}$ are equivalent if there exists
a unitary $v: \CD_E \rightarrow \CD_{E'}$ such that their symbols
satisfy $\Theta_{C,E} = \Theta_{C,E'}\, v$. Compared with the
analogous result for subisometric liftings contained in Theorem 4.1.6
the modifications necessary to prove Theorem 4.3.7 are technical and
straightforward, so we omit the proof. The important thing to
recognize is that, if a lifting $\uE$ is reduced, we have the
functional model \eqref{3.17} for it which is built only from $\uC$
and from the characteristic function $M_{C,E}$.
\\

Conversely, if $\uC$ on $\CH_C$ is a row contraction and
\begin{align*}
\tilde{M}: \Gamma \otimes \CD \rightarrow \Gamma \otimes \CD_C
\end{align*}
is an arbitrary contractive multi-analytic function (where $\CD$ is
any Hilbert space), then we can define
\begin{align*}
\Delta :(\eins - \tilde{M}^* \tilde{M})^{\frac{1}{2}}: \; \Gamma
\otimes \CD \rightarrow \Gamma \otimes \CD
\end{align*}
\begin{align*}
\tilde{\CH} := \CH_C \oplus (\Gamma \otimes \CD_C) \oplus
\overline{\Delta (\Gamma \otimes \CD)}
\end{align*}
\begin{align*}
\tilde{W}: \Gamma \otimes \CD \rightarrow (\Gamma \otimes \CD_C)
\oplus \overline{\Delta (\Gamma \otimes \CD)},\; x \mapsto
\tilde{M}\, x \oplus \Delta \, x
\end{align*}
$\tilde{W}$ is isometric and by introducing a copy $\CH_A$ of the
orthogonal complement of $\tilde{W}(\Gamma \otimes \CD)$ we can
extend $\tilde{W}$ to a unitary
\begin{align*}
\tilde{W}: \CH_A \oplus (\Gamma \otimes \CD) \rightarrow (\Gamma
\otimes \CD_C) \oplus \overline{\Delta (\Gamma \otimes \CD)}
\end{align*}
Let $\tilde{\uV} = (\tilde{V}_1,\ldots,\tilde{V}_d)$ be defined on
$\tilde{\CH}$ by $\tilde{V}_i := V^C_i \oplus Y_i$ (for
$i=1,\ldots,d$), where $\uV^C$ is the mid of $\uC$ on $\CH_C \oplus
(\Gamma \otimes \CD_C)$ \eqref{1.3} and $Y_i$ is given by
\begin{align*}
Y_i \Delta \,x := \Delta (L_i \otimes \eins) x \quad (\mbox{where}
\; x \in \Gamma \otimes \CD)
\end{align*}
It is not difficult to check that $\uY$ (and hence also
$\tilde{\uV}$) is a row contraction consisting of isometries with
orthogonal ranges (i.e., a row isometry). Further
\begin{align*}
\tilde{W}(\Gamma \otimes \CD) = \{ \tilde{M}\, x \oplus \Delta\,x,\;
x \in \Gamma \otimes \CD \}
\end{align*}
is invariant for the $\tilde{V}_i$. With $E^*_i := \tilde{V}^*_i
|_{\CH_C \oplus \tilde{W} \CH_A}, \; A^*_i := \tilde{V}^*_i
|_{\tilde{W} \CH_A}$ for $i=1,\ldots,d$ we obtain a contractive
lifting $\uE$ of $\uC$ by $\uA$ which we may call {\it the lifting
associated to the multi-analytic function} $\tilde{M}$. The
following result gives another justification for considering reduced
liftings.
\begin{prop}
The contractive lifting $\uE$ associated to a row contraction $\uC$
and a contractive multi-analytic function $\tilde{M}: \Gamma \otimes
\CD \rightarrow \Gamma \otimes \CD_C$ (where $\CD$ is any Hilbert
space) is reduced.
\end{prop}
\begin{proof}
By Lemma 4.3.5 it is enough to show that any vector $y \in \tilde{W}
\CH_A$ which is orthogonal to $\Gamma \otimes \CD_C$ is the zero
vector. But $y \in \tilde{W} \CH_A$ means that $y$ is orthogonal to
$\tilde{M}\, x \oplus \Delta\,x$ for all $x \in \Gamma \otimes \CD$
and $y$ orthogonal to $\Gamma \otimes \CD_C$ means that $y \in 0
\oplus \overline{\Delta (\Gamma \otimes \CD)}$. Hence indeed $y=0$.
\end{proof}

Proposition 4.3.8 shows that the theory of characteristic functions
cannot be extended beyond reduced liftings. Note that $\tilde{M}$ is
not necessarily the characteristic function of the associated
lifting $\uE$ and we used $\,\tilde{}\,$ to indicate this. It is an
interesting question which intrinsic properties of $\tilde{M}$
guarantee that it is the characteristic function. We leave this as
an open problem.

\section{Properties of the characteristic function}

First we shall compute an explicit expression for the characteristic
function of a reduced lifting. We continue to use the notation of
the previous section and consider a reduced lifting $\uE$ on $\CH_E
= \CH_C \oplus \CH_A$ of $\uC$ on $\CH_C$ by $\uA$ on $\CH_A$. As in
\eqref{3.8} the row isometry $(\underline{L \otimes \eins}) \oplus
\uY$ on $(\Gamma \otimes \CD_C) \oplus \CK$ restricts to a mid of
$\uA$ (transferred to $W \CH_A$). So we have for all $\alpha \in
\tilde{\Lambda}$ and $h \in \CH_A$
\begin{align} \label{4.1}
\big[ (L^*_\alpha \otimes \eins) \oplus Y^*_\alpha \big] W h = W \,
A^*_\alpha h
\end{align}
Using \eqref{3.11} we infer that
\begin{align} \label{4.2}
\gamma D_{*,A} A^*_\alpha h = P_{e_0 \otimes \CD_C} W \, A^*_\alpha
h = P_{e_0 \otimes \CD_C} \big[ (L^*_\alpha \otimes \eins) \oplus
Y^*_\alpha \big] W h = P_{e_\alpha \otimes \CD_C} W h
\end{align}
which yields a Poisson kernel type formula, compare \eqref{1.8}:
\begin{align} \label{4.3}
P_{\Gamma \otimes \CD_C} W h = \sum_{\alpha \in \tilde{\Lambda}}
e_\alpha \otimes \gamma D_{*,A} A^*_\alpha h
\end{align}
To compute the symbol $\Theta_{C,E}$ of the characteristic function
we define $d^i_h := (V^E_i - E_i) h = e_0 \otimes (D_E)_i h$ and use
the identification of $\CD_E$ with the closed linear span of all
$d^i_h$ with $i=1,\ldots,d$ and $h \in \CH_E$, see \eqref{1.3}.
Then, using \eqref{3.1} and the Definition 4.3.6 of $\Theta_{C,E}$, we
obtain
\begin{align} \label{4.4}
\Theta_{C,E} d^i_h = P_{\Gamma \otimes \CD_C} W (V^E_i - E_i) h =
P_{\Gamma \otimes \CD_C} V^C_i P_{\CH_C \oplus (\Gamma \otimes
\CD_C)} W h - P_{\Gamma \otimes \CD_C} W E_i h
\end{align}
We distinguish two cases.

\vspace{0.2cm}

{\bf Case I:} $h \in \CH_C.$
\begin{align*}
P_{\Gamma \otimes \CD_C} V^C_i P_{\CH_C \oplus (\Gamma \otimes
\CD_C)} W h = P_{\Gamma \otimes \CD_C} V^C_i h = [e_0 \otimes
(D_C)_i h ] \quad\quad {\mbox{by} \;\eqref{1.3}}
\end{align*}
\begin{align*}
P_{\Gamma \otimes \CD_C} W E_i h = P_{\Gamma \otimes \CD_C} W (C_i h
\oplus B_i h) = \sum_\alpha e_\alpha \otimes \gamma \DsA A^*_\alpha
B_i h \quad\quad {\mbox{by} \;\eqref{4.3}}
\end{align*}
and thus
\begin{align} \label{4.5}
\Theta_{C,E}  d^i_h = e_0 \otimes [(D_C)_i h - \gamma \DsA B_i h] -
\sum_{|\alpha|\geq 1} e_\alpha \otimes \gamma \DsA A^*_\alpha B_i h
\end{align}

{\bf Case II:} $h \in\, \CH_A$.
\begin{align*}
P_{\Gamma \otimes \CD_C} V^C_i P_{\CH_C \oplus (\Gamma \otimes
\CD_C)} W h = V^C_i P_{\Gamma \otimes \CD_C} W h
\end{align*}
\begin{align*}
= (L_i \otimes \eins) P_{\Gamma \otimes \CD_C} W h = \sum_\alpha e_i
\otimes e_\alpha \otimes \gamma \DsA A^*_\alpha h
\end{align*}
\begin{align*}
P_{\Gamma \otimes \CD_C} W E_i h = P_{\Gamma \otimes \CD_C} W A_i h
= \sum_\beta e_\beta \otimes \gamma \DsA A^*_\beta A_i h
\end{align*}
Note that for $h \in\, \CH_A$ we have $(D_A)_i h = (D_E)_i h$ (which
we identify with $d^i_h$) because $\uE$ is an extension of $\uA$.
With $P_j$ the orthogonal projection onto the $j-$th component we
obtain
\begin{align} \label{4.6}
\Theta_{C,E}\, d^i_h = - e_0 \otimes \gamma \DsA A_i h +
\sum^d_{j=1} e_j \otimes \sum_\alpha e_\alpha \otimes \gamma \DsA
A^*_\alpha (\delta_{ji} \eins -A^*_j A_i) h \nonumber
\\
= - e_0 \otimes \gamma \uA (D_A)_i h + \sum^d_{j=1} e_j \otimes
\sum_\alpha e_\alpha \otimes \gamma \DsA A^*_\alpha P_j D_A (D_A)_i
h \nonumber
\\
= - e_0 \otimes \gamma \sum^d_{j=1} A_j P_j d^i_h + \sum^d_{j=1} e_j
\otimes \sum_\alpha e_\alpha \otimes \gamma \DsA A^*_\alpha P_j D_A
d^i_h
\end{align}
We note that if $\gamma$ is omitted from \eqref{4.6} then we obtain
exactly Popescu's definition of the characteristic function of the
(c.n.c.) row contraction $\uA$ as given in \cite{Po89b}. Hence Case
II is essentially the characteristic function of $\uA$ ,
contractively embedded by $\gamma$. In a special case this has been
observed in \cite{DG07a} and, because this special case was
subisometric and hence $\gamma$ isometric, $\Theta$ was called an
extended characteristic function. \eqref{4.6} generalizes this idea.
\\

Let us now illustrate how the characteristic function factorizes for
iterated liftings. Assume that $\tilde{\uE}$ on $\CH_{\tilde{E}}$ is
a two step lifting of the row contraction  $\uC$ on $\CH_C$, i.e.,
$\uE$ on $\CH_E$ with $E_i= \left(\begin{array}{cc}
C_i & 0   \\
B_i & A_i
\end{array}
\right)$ (for $i=1,\ldots,d$) is a contractive lifting of $\uC$ on
$\CH_C$ by $\uA$ on $\CH_A$ (as before) and $\tilde{\uE}$ on
$\CH_{\tilde{E}}$ with $\tilde{E}_i= \left(\begin{array}{cc}
E_i & 0   \\
*   & \tilde{A}_i
\end{array}
\right)$ (for $i=1,\ldots,d$) is a contractive lifting of $\uE$ on
$\CH_E$ by $\tilde{\uA}$ on $\CH_{\tilde{A}}$. Then $\CH_{\tilde{E}}
= \CH_E \oplus \CH_{\tilde{A}} = \CH_C \oplus \CH_A \oplus
\CH_{\tilde{A}}$ and with respect to this decomposition
\begin{align} \label{4.7}
\tilde{E}_i=\left(\begin{array}{ccc}
C_i & 0   & 0 \\
*   & A_i & 0 \\
*   & *   & \tilde{A}_i
\end{array}
\right)
\end{align}
`$*$' stands for entries which we do not need to name explicitly.

\begin{thm}
If the liftings $\uE$ of $\uC$ and $\tilde{\uE}$ of $\uE$ are
reduced then also the lifting $\tilde{\uE}$ of $\uC$ is reduced, and
the characteristic functions factorize:
\begin{align} \label{4.8}
M_{C,\tilde{E}}=M_{C,E} \; M_{E,\tilde{E}}.
\end{align}
\end{thm}

\begin{proof}
As in \eqref{3.1} we obtain the following unitaries from the given
liftings:
\begin{align*}
W: \CH_E \oplus (\Gamma \otimes \CD_E) \rightarrow \CH_C \oplus
(\Gamma \otimes \CD_C) \oplus \CK
\end{align*}
\begin{align*}
\tilde{W}: \CH_{\tilde{E}} \oplus (\Gamma \otimes \CD_{\tilde{E}})
\rightarrow \CH_E \oplus (\Gamma \otimes \CD_E) \oplus \tilde{\CK}
\end{align*}
satisfying
\begin{align*}
W \, V^E_i = (V^C_i \oplus Y_i)\, W
\end{align*}
\begin{align*}
\tilde{W} \, V^{\tilde{E}}_i = (V^E_i \oplus \tilde{Y}_i)
\,\tilde{W}
\end{align*}
We can define another unitary
\begin{align} \label{4.9}
Z := (W \otimes \eins_{\tilde{\CK}})\,\tilde{W}: \CH_{\tilde{E}}
\oplus (\Gamma \otimes \CD_{\tilde{E}}) \rightarrow \CH_C \oplus
(\Gamma \otimes \CD_C) \oplus \CK \oplus \tilde{\CK}
\end{align}
satisfying
\begin{align} \label{4.10}
Z \, V^{\tilde{E}}_i = (V^C_i \oplus Y_i \oplus \tilde{Y}_i) \,Z
\end{align}
Note further that $W, \tilde{W}$ and hence also $Z$ act identically
on $\CH_C$. By assumption the liftings $\uE$ of $\uC$ and
$\tilde{\uE}$ of $\uE$ are reduced and we have characteristic
functions
\begin{align*}
M_{C,E} = P_{\Gamma \otimes \CD_C} W |_{\Gamma \otimes \CD_E}
\end{align*}
\begin{align*}
M_{E,\tilde{E}} = P_{\Gamma \otimes \CD_E} \tilde{W} |_{\Gamma
\otimes \CD_{\tilde{E}}}
\end{align*}
They can be composed to yield a multi-analytic operator
\begin{align*}
M := M_{C,E} \, M_{E,\tilde{E}}: \; \Gamma \otimes \CD_{\tilde{E}}
\rightarrow \Gamma \otimes \CD_C
\end{align*}
Using \eqref{4.9} it is easily checked that
\begin{align*}
M = P_{\Gamma \otimes \CD_C} Z |_{\Gamma \otimes \CD_{\tilde{E}}}
\end{align*}
We conclude by \eqref{4.10} that the lifting $\tilde{\uE}$ of $\uC$
is associated to $M$ and hence, by Proposition 4.3.8, this lifting is
reduced. In fact, comparing with Definition 4.3.6, we see that $M$ is
the characteristic function, i.e., $M = M_{C,\tilde{E}}$.
\end{proof}

\section{Applications to completely positive maps}

If $\uT = (T_1,\ldots,T_d)$ is a row contraction on a Hilbert space
$\CK$ then we denote by $\Phi_T$ the corresponding (normal)
completely positive map on $\CB(\CK)$ given by
\begin{align} \label{5.1}
\Phi_T(\cdot) = \sum^d_{i=1} T_i \cdot T^*_i
\end{align}
If $d=\infty$ this should be understood as a SOT-limit. See for
example \cite{Pau03} for the general theory of completely positive
maps, we shall only work with the concrete representation
\eqref{5.1}. The fact that $\uT$ is a row contraction implies that
$\Phi_T(\eins) \leq \eins$, i.e., $\Phi_T$ is contractive. It is
unital ($\Phi_T(\eins) = \eins$) if and only if $\uT$ is
coisometric.

If $\uE$ is a contractive lifting of $\uC$ by $\uA$, i.e., $E_i =
\left(\begin{array}{cc}
C_i & 0   \\
B_i & A_i
\end{array}
\right)$ (for $i=1,\ldots,d$) then an elementary computation shows
that
\begin{align} \label{5.2}
\Phi_E \left(\begin{array}{cc}
X_{11} & X_{12}   \\
X_{21} & X_{22}
\end{array}
\right) \quad = \quad \sum^d_{i=1} \left(\begin{array}{cc}
C_i X_{11} C^*_i                    & C_i X_{11} B^*_i + C_i X_{12} A^*_i \\
&\\
   B_i X_{11} C^*_i                                 & B_i X_{11} B^*_i + B_i X_{12} A^*_i \\
 + A_i X_{21} C^*_i &+A_i X_{21} B^*_i + A_i X_{22}
A^*_i
\end{array}
\right)
\end{align}
with $X_{11} \in \CB(\CH_C)$, $X_{12} \in \CB(\CH_A,\CH_C)$, $X_{21}
\in \CB(\CH_C,\CH_A)$, $X_{22} \in \CB(\CH_A)$. We denote by $p_C=
\left(\begin{array}{cc}
\eins & 0   \\
0 & 0
\end{array}
\right)$ and $p_A= \left(\begin{array}{cc}
0 & 0   \\
0 & \eins
\end{array}
\right)$ the orthogonal projections onto $\CH_C$ and $\CH_A$. The
following facts are immediate from \eqref{5.2}.

\begin{align} \label{5.3}
p_C \,(\Phi_E)^n \left(\begin{array}{cc}
X & 0   \\
0 & 0
\end{array}
\right) |_{\CH_C} \quad = \quad (\Phi_C)^n (X)
\end{align}
(for $n\in\N_0$ and $X\in\CB(\CH_C)$)

\begin{align} \label{5.4}
\Phi_E \left(\begin{array}{cc}
0 & 0   \\
0 & Y
\end{array}
\right) \quad = \quad \left(\begin{array}{cc}
0 & 0   \\
0 & \Phi_A(Y)
\end{array}
\right)
\end{align}
(for $Y\in \CB(\CH_A)$). So $\Phi_E$ is a kind of (power) dilation
of $\Phi_C$ \eqref{5.3} and an extension of $\Phi_A$ \eqref{5.4}.

\begin{defn}
If $\CH_E = \CH_C \oplus \CH_A$, $\Phi_E: \CB(\CH_E) \rightarrow
\CB(\CH_E)$, $\Phi_C: \CB(\CH_C) \rightarrow \CB(\CH_C)$, $\Phi_A:
\CB(\CH_A) \rightarrow \CB(\CH_A)$ are contractive normal completely
positive maps such that \eqref{5.3} and \eqref{5.4} are valid then
we say that $\Phi_E$ is a contractive lifting of $\Phi_C$ by
$\Phi_A$.
\end{defn}

We have seen that a contractive lifting of row contractions gives
rise to a contractive lifting of completely positive maps. The
converse is also true: Let us assume \eqref{5.4}. If $\Phi_E (\cdot)
= \sum^d_{i=1} E_i \cdot E^*_i$ and we write $E_i=
\left(\begin{array}{cc}
C_i & D_i   \\
B_i & A_i
\end{array}
\right)$ for the moment, then
\begin{align*}
\Phi_E \left(\begin{array}{cc}
0 & 0   \\
0 & \eins
\end{array}
\right) \quad = \quad \left(\begin{array}{cc} \vspace{0.2cm}
\sum^d_{i=1} D_i D^*_i & * \\
           *           & * \\
\end{array}
\right)
\end{align*}
and \eqref{5.4} implies that all the $D_i$ are zero, i.e., we have a
lifting of row contractions. So actually \eqref{5.4} implies
\eqref{5.3} with some $\Phi_C$.

Note that if $\uE = \uV^C$, the mid of $\uC$, then $\Phi_E$ is a
$*$-homomorphism and \eqref{5.3} shows that the powers of $\Phi_E$
are a homomorphic dilation of the completely positive semigroup
formed by powers of $\Phi_C$. See \cite{BP94,Ar03,Go04} for further
information about this kind of dilation theory.

The discussion above shows that we can use our theory of liftings
for row contractions to study liftings of completely positive maps.
If $\uE$ is a reduced lifting of $\uC$ by $\uA$ then we have a
characteristic function $M_{C,E}$. It is well known (see for example
\cite{Pau03,Go04}) that in the decomposition $\Phi_E {\cdot} =
\sum^d_{i=1} E_i \cdot E^*_i$ the tuple $(E_1,\ldots,E_d)$ is not
uniquely determined and that $\sum^d_{i=1} E'_i \cdot (E'_i)^*$
describes the same map if and only if $\uE'$ is obtained from $\uE$
by multiplication with a unitary $d \times d$-matrix (with complex
entries). This does not change the characteristic function because
the latter is defined as an intertwiner between objects which are
transformed in the same way. Hence it is possible to think of
$M_{C,E}$ also as the characteristic function of a reduced lifting
$\Phi_E$ of $\Phi_C$ by $\Phi_A$. (Of course we call this lifting
reduced if the corresponding lifting of row contractions is
reduced.) Theorem 4.3.7 translates immediately into

\begin{cor}
Given $\Phi_C$, two reduced liftings $\Phi_E$ resp. $\Phi_{E'}$ of
$\Phi_C$ by $\Phi_A$ resp. $\Phi_{A'}$ are conjugate, i.e.
\begin{align*}
\Phi_E = U^* \Phi_{E'} (UXU^*) U
\end{align*}
with a unitary $U: \CH_E \rightarrow \CH_{E'}$ such that $U
|_{\CH_C} = \eins |_{\CH_C}$, if and only if the corresponding
characteristic functions are equivalent.
\end{cor}

\noindent Corollary 4.5.2 generalizes Corollary 6.3 in \cite{DG07a}
where $dim\, \CH_C = 1$.
\\

In the following we confine ourselves mainly to liftings which are
coisometric and subisometric and give some concrete and useful
results about the corresponding completely positive maps.

\begin{lem}
Let $\uE$ be a contractive lifting of a row contraction $\uC$ by a
$*$-stable row contraction $\uA.$ Then for all $X_{12}, X_{21},
X_{22}$
\[
\Phi^n_E \left(
\begin{array}{cc}
0  &  X_{12} \\
X_{21}  &  X_{22} \\
\end{array}
\right) \to 0
\]
as $n \to \infty$ (SOT).
\end{lem}

\begin{proof}
$\Phi^n_E (p_A)$ decreases to zero in the strong operator topology
because of \eqref{5.4} and the assumption that $\uA$ is $*$-stable.
Then also $ \Phi^n_E \left(
\begin{array}{cc}
0  &  0 \\
0  &  X_{22} \\
\end{array}
\right) \to 0, $ first for $0 \leq X_{22} \leq \|X_{22}\| \,p_A$,
then for general $X_{22}$ by writing it as a linear combination of
positive elements. Using the Kadison-Schwarz inequality for
completely positive maps (cf. \cite{Ch74} or \cite{Pau03}, Chapter
3) we obtain
\begin{eqnarray*}
\Phi^n_E \left(
\begin{array}{cc}
0  &  0 \\
X^*_{12}  &  0 \\
\end{array}
\right)\; \Phi^n_E \left(
\begin{array}{cc}
0  &  X_{12} \\
0  &  0 \\
\end{array}
\right) &\leq & \Phi^n_E \left(
\begin{array}{cc}
0  &  0 \\
0 & X^*_{12} X_{12} \\
\end{array}
\right) \to 0
\end{eqnarray*}
and hence $ \Phi^n_E \left(
\begin{array}{cc}
0  &  X_{12} \\
0  &  0 \\
\end{array}
\right) \to 0$. Similarly $ \Phi^n_E \left(
\begin{array}{cc}
0  &  0 \\
X_{21} &  0 \\
\end{array}
\right) \to 0$.
\end{proof}

\begin{thm}
Suppose the row coisometry $\uE$ is a lifting of $\uC$ by $\uA$.
Then the following assertions are equivalent:
\begin{itemize}
\item[(a)]
The lifting is subisometric.
\item[(b)]
$\uA$ is $*$-stable.
\item[(c)]
$(\Phi_E)^n (p_C) \to \eins \quad (n \to \infty, SOT)$
\item[(d)]
There is an order isomorphism between the fixed point sets of
$\Phi_E$ and of $\Phi_C$ given by
\begin{align} \label{5.6}
\kappa: X \mapsto p_C \, X \, p_C
\end{align}
\end{itemize}
In this case, $\kappa$ is isometric on selfadjoint elements. If $x$
is a fixed point of $\Phi_C$ then we can reconstruct the preimage
$\kappa^{-1}(x)$ as the SOT-limit
\begin{align} \label{5.7}
\lim_{n \to \infty} (\Phi_E)^n \left(
\begin{array}{cc}
x  &  0 \\
0  &  0 \\
\end{array}
\right)
\end{align}
\end{thm}

Recall further that by the results of Section 4.2 the liftings in
Theorem 4.5.4 are parametrized by $*$-stable row contractions $\uA$
with $dim\,\CD_{*,A} \leq dim\,\CD_C$ together with isometries
$\gamma: \CD_{*,A} \rightarrow \CD_C$ and that they can be
explicitly constructed from these data. Theorem 4.5.4(d) tells us that
(exactly) for such liftings the maps $\Phi_E$ and $\Phi_C$ have
closely related properties in terms of their fixed points. We can
identify this useful situation by checking the convenient conditions
(b) or (c).

\begin{proof}
By Proposition 4.2.3 a coisometric lifting $\uE$ of $\uC$ by $\uA$ is
subisometric if and only if $\uA$ is $*$-stable. Using \eqref{5.4}
the latter means that
\begin{align*}
(\Phi_E)^n (p_A) \to 0 \quad (n \to \infty,\; SOT),
\end{align*}
which is equivalent to (c) because $\Phi_E$ is unital. Hence $(a)
\Leftrightarrow (b) \Leftrightarrow (c)$.

If $X = \left(\begin{array}{cc}
x & *   \\
* & *
\end{array}
\right)$ is a fixed point of $\Phi_E$ then it is immediate from
\eqref{5.2} that $x$ is a fixed point of $\Phi_C$. Hence $\kappa: X
\mapsto p_C \, X \, p_C$ indeed maps fixed points of $\Phi_E$ to
fixed points of $\Phi_C$. (This is true for all contractive
liftings.) Now assume (a), i.e., the lifting is subisometric. Then
\begin{align*}
X = \Phi_E (X) = \lim_{n\to\infty} (\Phi_E)^n (X) =
\lim_{n\to\infty} (\Phi_E)^n \left(
\begin{array}{cc}
x  &  0 \\
0  &  0 \\
\end{array}
\right),
\end{align*}
where the last equality follows from Lemma 4.5.3. Hence $\kappa$ is
injective.

Let $\uV = (V_1,\ldots,V_d)$ simultaneously serve as mid for $\uC$
and $\uE$. Then Theorem 4.5.1 in \cite{BJKW00} or Lemma 4.6.4 in the
Appendix of this paper show that for every fixed point $x$ of
$\Phi_C$ there exists $A'$ in the commutant of $V_1,\ldots,V_d$ such
that $p_C A' p_C = x$. Define $X := p_E A' p_E$, where $p_E$ is the
orthogonal projection onto $\CH_E$. Then, using the lifting property
$E_i\,p_E = p_E \, V_i$ for $i=1,\ldots,d$ for the mid and the fact
that $\sum^d_{i=1} V_i V^*_i = \eins$ (because $\uE$ is coisometric
also $\uV$ is coisometric), we find that
\begin{align*}
\Phi_E (X) = \sum^d_{i=1} E_i X E^*_i = \sum^d_{i=1} E_i p_E A' p_E
E^*_i
= \sum^d_{i=1} p_E V_i A' V^*_i p_E \\
= p_E A' \sum^d_{i=1} V_i V^*_i p_E = p_E A' p_E = X
\end{align*}
So $X$ is a fixed point of $\Phi_E$ and clearly $\kappa(X)=x$. We
conclude that $\kappa$ is also surjective. The fact that $\kappa$ is
isometric on selfadjoint elements is also a consequence of Lemma
4.6.4.

On the other hand, if the lifting $\uE$ of $\uC$ is not subisometric
then the mid $\uV^C$ of $\uC$ is embedded on a proper reducing
subspace $\hat{\CH}_C$ into the space $\hat{\CH}_E$ of the mid
$\uV^E$ of $\uE$. Then $\eins_{\hat{\CH}_E}$ and $p_{\hat{\CH}_C}$
are two different fixed points of $\Phi_{V^E}$. By Lemma 4.6.4 the map
$\hat{X} \mapsto p_E \hat{X} p_E$ maps them into different fixed
points of $\Phi_E$: $p_E \eins_{\hat{\CH}_E} p_E = p_E \not= p_E
p_{\hat{\CH}_C} p_E$. If $\kappa: X \mapsto p_C \, X \, p_C$ from
the fixed point set of $\Phi_E$ into the fixed point set of $\Phi_C$
were injective then also $p_C p_E p_C \not= p_C p_E p_{\hat{\CH}_C}
p_E p_C$. But both sides are equal to $p_C$. Hence in this case
$\kappa$ is not injective. We have proved $(a) \Leftrightarrow (d)$.
\end{proof}

Recall that a unital completely positive map $\Phi_E$ is called {\it
ergodic} if there are no other fixed points than the multiples of
the identity. By abuse of language we also call $\uE$ ergodic in
this case (as in \cite{DG07a}).

\begin{prop}
Let $\uE$ be a coisometric lifting of $\uC$ by $\uA$. Then $\uE$ is
ergodic if and only if $\uC$ is ergodic and $\uA$ is $*$-stable.
\end{prop}

\begin{proof}
If $\uA$ is $*$-stable then use the equivalence $(b) \Leftrightarrow
(d)$ in Theorem 4.5.4 and infer from $\uC$ ergodic that also $\uE$ is
ergodic. Further note that, because $\uE$ is coisometric, we always
have
\begin{align*}
\Phi_E \left(
\begin{array}{cc}
\eins  &  0 \\
0  &  0 \\
\end{array}
\right) \quad = \quad \left(
\begin{array}{cc}
\eins  &  0 \\
0  &  \uB \uB^* \\
\end{array}
\right) \quad \geq \quad \left(
\begin{array}{cc}
\eins  &  0 \\
0  &  0 \\
\end{array}
\right)
\end{align*}
We say that $p_C \left(
\begin{array}{cc}
\eins  &  0 \\
0  &  0 \\
\end{array}
\right)$ is an increasing projection for $\Phi_E$. Hence $(\Phi_E)^n
(p_C)$ increases to a SOT-limit which clearly is a fixed point of
$\Phi_E$.
\\
Now let $\uE$ be ergodic. Then all fixed points are multiples of
$\left(
\begin{array}{cc}
\eins  &  0 \\
0  &  \eins \\
\end{array}
\right)$ and because the left upper corner of $(\Phi_E)^n (p_C)$ is
always $\eins$ we have $(\Phi_E)^n (p_C) \to \left(
\begin{array}{cc}
\eins  &  0 \\
0  &  \eins \\
\end{array}
\right)$. We have verified Theorem 4.5.4(c) and now Theorem 4.5.4(d) and
(b) show that $\uC$ is ergodic and that $\uA$ is $*$-stable.
\end{proof}

This generalizes Proposition 2.3 in \cite{DG07a} where $\CH_C$ is
one dimensional and hence $\uC$ ergodic is automatically fulfilled.
\\

The following provides an interesting example for the liftings
considered above. Let $\Phi_E: \CB(\CH_E) \rightarrow \CB(\CH_E)$ be
any (normal) unital completely positive map and let $\psi$ be a
normal invariant state, i.e., $\psi \circ \Phi_E = \psi$. Define
$\CH_C$ to be the support of $\psi$ (cf.\,\cite{Ta01}) and let
$\CH_A$ be the orthogonal complement, so $\CH_E = \CH_C \oplus
\CH_A$. Then $\uE = (E_1,\ldots,E_d)$ is a coisometric lifting of
$\uC = (C_1,\ldots,C_d)$ if we define $C^*_i := E^*_i |_{\CH_C}$ for
$i=1,\ldots,d$. In fact $p_C E^*_i p_C = E^*_i p_C$ for all $i$ by
Lemma 6.1 of \cite{BJKW00}. Note that the compression $\Phi_C$ has a
faithful normal invariant state, the restriction $\psi_C$ of $\psi$
to $\CB(\CH_C)$. Conversely we can start with $\Phi_C$ and a
faithful invariant state $\psi_C$ and construct liftings $\Phi_E$ .
They have normal invariant states given by $\psi(X) := \psi_C (p_C X
p_C)$. From Proposition 4.2.3 and Theorem 4.5.4 we conclude

\begin{cor}
Let $\Phi_C: \CB(\CH_C) \rightarrow \CB(\CH_C)$ be a (normal) unital
completely positive map with a faithful normal invariant state
$\psi_C$. Then we have a one-to-one correspondence between
\begin{itemize}
\item[(a)]
(normal) unital completely positive maps $\Phi_E: \CB(\CH_E)
\rightarrow \CB(\CH_E)$ with normal invariant state $\psi$ such that
the support of $\psi$ is $\CH_C$ and $\psi |_{\CB(\CH_C)} = \psi_C$,
compression of $\Phi_E$ is $\Phi_C$ and $(\Phi_E)^n (p_C) \to \eins
\quad (n \to \infty,\, SOT)$
\item[(b)]
$*$-stable $\uA$ with $dim\,\CD_{*,A} \leq dim\,\CD_C$ together with
isometries $\gamma: \CD_{*,A} \rightarrow \CD_C$
\end{itemize}
There exist order isomorphisms $\kappa_E: X \mapsto p_C X p_C$
between the fixed point sets of these maps $\Phi_E$ and the fixed
point set of $\Phi_C$.
\end{cor}

In the special case when $\psi$ is an invariant vector state
$\langle \xi, \cdot \,\xi \rangle$ of $\Phi_E$ we have the result
that $\Phi_E$ is ergodic if and only if $(\Phi_E)^n (p_\xi) \to
\eins \; (n \to \infty,\, SOT)$, where $p_\xi$ is the orthogonal
projection onto $\C \xi$, cf. \cite{Go04}, A.5.2. Hence we obtain a
classification of such maps. Here $\CD_C$ is $(d-1)$-dimensional.
This case has been further investigated in \cite{DG07a}.

Corollary 4.5.6 is useful because many techniques only apply to
completely positive maps with faithful invariant states, cf.
\cite{Ku03}. It enables us to transfer information from the faithful
to the non-faithful setting. For example, it is known that in the
case of a faithful normal invariant state the fixed point set is an
algebra (cf. \cite{Ch74,FNW94,BJKW00}). Now $\kappa$ is an order
isomorphism but it is not in general multiplicative. In fact, there
are examples of completely positive maps with a normal invariant
non-faithful state where the fixed point set is not an algebra (cf.
\cite{Ar69,Ar72,BJKW00}. If Corollary 4.5.6 applies we can think of
it as an (order isomorphic) deformation of an algebra.

\section{Appendix}

In Section 4.5 we needed a commutant lifting theorem (Theorem 5.1 of
\cite{BJKW00}) which says that the fixed point set of a normal
unital completely positive map is in one-to-one correspondence with
the commutant of the Cuntz algebra representation generated by the
mid. Below we give a variant of the proof which is based on a
Radon-Nikodym result for completely positive maps by W.Arveson. This
is a good way to think about it and it supports the understanding of
the other arguments in the main text.

\begin{lem} \cite{Ar69}, Theorem 1.4.2 \\
If $\Psi$ is a completely positive map from a $C^*$-algebra $\CB$ to
$\CB(\CH)$, with $\CH$ a Hilbert space, then there exists an affine
order isomorphism of the partially ordered set of operators $\{ A'
\in \pi(\CB)': 0 \leq A' \leq \eins \}$ onto $[0,\Psi]$. Here $\pi$
is the minimal Stinespring representation of $\CB$ associated to
$\Psi$ and $[0,\Psi]$ is the order interval containing all
completely positive maps $\Phi: \CB \rightarrow \CB(\CH)$ with $0
\leq \Phi \leq \Psi$. The order relation for completely positive
maps used here is $\Phi \leq \Psi$ if $\Psi - \Phi$ is completely
positive.
\\
Explicitly, if $\Psi(x) = W^* \pi(x) W$ is the minimal Stinespring
representation of $\Psi$ then $A' \in \pi(\CB)'$ corresponds to
$\Phi = W^* A' \pi(x) W$.
\end{lem}

\begin{lem} \cite{BJKW00}, Corollary 2.4; \cite{Po03}, Theorem 2.1 \\
If $0 \leq D \leq \eins$ is a fixed point of the (normal unital
completely positive) map $\Phi_R(\cdot) = \sum^d_1 R_i \cdot R^*_i$
on $\CB(\CH)$ then there exists a completely positive map $\Psi_D:
\Od \rightarrow \CB(\CH), \; V_\alpha V^*_\beta \mapsto R_\alpha D
R^*_\beta$. Here $\alpha,\beta \in \tilde{\Lambda}$ and $\Od$ is the
Cuntz algebra generated by the $V_i$, where $\uV = (V_1,\ldots,V_d)$
is a mid of $\uR = (R_1,\ldots,R_d)$.
\end{lem}

Using notation from the previous lemmas we get

\begin{lem}
There exists an affine order isomorphism $D \mapsto \Psi_D$ between
\[
\{ 0 \leq D \leq \eins: D \;\mbox{is a fixed point of}\;
\Phi_R(\cdot) = \sum^d_1 R_i \cdot R^*_i \;\mbox{on}\; \CB(\CH) \}
\]
and $[0,\Psi_\eins]$, where $\Psi_\eins$ is the completely positive
map described in Lemma 4.6.2 with $D = \eins$, i.e., $\Psi_\eins: \Od
\rightarrow \CB(\CH), \; V_\alpha V^*_\beta \mapsto R_\alpha
R^*_\beta$.
\end{lem}

\begin{proof}
From $\Psi_\eins = \Psi_D + \Psi_{\eins - D}$ we see that $\Psi_D
\in [0,\Psi_\eins]$ for all fixed points $ 0 \leq D \leq \eins$ of
$\Phi_R$. On the other hand, if $\Phi \in [0,\Psi_\eins]$ then by
Lemma 4.6.1 with $\CB = \Od$ there exists $A' \in \pi(\CB)'$ with $0
\leq A' \leq \eins$ such that $\Phi(x) = W^* A' \pi(x) W$, where
$\Psi_\eins(x) = W^* \pi(x) W$ is a minimal Stinespring
representation. Using that $(V_1,\ldots,V_d)$ is a mid of $\uR =
(R_1,\ldots,R_d)$ it is easily checked that $\Psi_\eins(x) = p
\pi(x) p$ is such a minimal Stimespring representation if $\pi$ is
the Cuntz algebra representation generated by $(V_1,\ldots,V_d)$ and
$p$ is the projection onto the space $\CH$. (In $p \pi(x) p$ the $p$
on the right hand side should be interpreted as the embedding of
$\CH$ into the dilation space.)
\\
Hence if $x = V_\alpha V^*_\beta$ then we obtain
\[
\Phi(V_\alpha V^*_\beta) = p A' V_\alpha V^*_\beta p = p V_\alpha A'
V^*_\beta p = p V_\alpha p A' p V^*_\beta p = R_\alpha p A' p
R^*_\beta.
\]
We conclude that $\Phi = \Psi_D$ with $D := p A' p$. Clearly $0 \leq
D \leq \eins$ and $D$ is a fixed point of $\Phi_R$ (because $\uV$ is
a coisometric lifting of $\uR$, i.e., $\sum^d_{i=1} V_i V^*_i =
\eins$ and $R_i p = p V_i$ for all $i$). The correspondence is
bijective ($\Psi_D(\eins) = D$) and it clearly respects the order.
\end{proof}

\begin{lem} \cite{BJKW00}, Theorem 5.1 \\
There is an affine order isomorphism between
$\{ 0 \leq D \leq \eins: D \;\mbox{is a fixed point}\\
\mbox{of}\;
 \Phi_R(\cdot) = \sum^d_1 R_i \cdot R^*_i
\;\mbox{on}\; \CB(\CH) \}$ and $\{ A' \in \pi(\Od)': 0 \leq A' \leq
\eins \}$, where $\pi$ is the Cuntz algebra representation generated
by the mid $\uV = (V_1,\ldots,V_d)$ of $\uR = (R_1,\ldots,R_d)$. It
is given by $A' \mapsto p A' p$, where $p$ is the projection onto
the space $\CH$. The isomorphism is isometric on the selfadjoint
parts.
\end{lem}

\begin{proof}
For the first part we only have to add to the arguments in the proof
of Lemma 4.6.3 the reminder that by Lemma 4.6.1 the correspondence
between $\{ A' \in \pi(\Od)': 0 \leq A' \leq \eins \}$ and $[0,
\Psi_\eins]$ is a bijection. As pointed out in \cite{BJKW00},
Section 4, it is isometric on the selfadjoint parts because $\eins$
is mapped to $\eins$ (identities on different Hilbert spaces) and
for selfadjoint elements $y$ we have $\|y\| = \inf \{ \alpha > 0:
-\alpha \eins \leq y \leq \alpha \eins \}$.
\end{proof}

\chapter{Constrained Liftings}

In the preceding chapter we have mainly exploited the dilation
theory of  (noncommutative) row contractions. There are some more
important types of dilation theories. They are for row contractions
satisfying some constraints. Formally:

\noindent {\bf Definition ~~}
{\em Assume $\uT$ to be a $d$-tuple of bounded operators on a common
Hilbert space $\CH$ and $\{p_\eta (\uz)\}_{\eta \in J}$ a finite set
of polynomials in $d$ noncommuting variables with index set $J.$
Then $\uT$ is said to be $J$-constrained if
\[ p_{\eta} (\uT)=0  \mbox{~~~for~~} \eta \in J.\]}

The corresponding dilations are called constrained dilations
(defined in Section 5.1). In \cite{BBD04}, \cite{De07a} and
\cite{BDZ06} the question how the mid is related to the constrained
dilation was addressed. Here we carry out this study further in
Section 5.1, especially when the row contraction is defined
 on a finite dimensional Hilbert space. There is a generalization
of Beurling's Theorem for Fock spaces by Popescu \cite{Po05} in
constrained case. A complete invariant for constrained liftings by
c.n.c. is obtained in Section 5.2 motivated by this Beurling type
result. These invariants are at times more interesting than the ones
in preceding chapters. For instance in the commuting case this
invariant is a $H^\infty$ function in the sense of multivariable
complex analysis. Similar analysis to classify constrained c.n.c.
row contractions can be found in
\cite{Ar03},\cite{Po07},\cite{BT07},\cite{BES05}, etc.

Beginning with a row contraction  $\uR=(R_1,\ldots,R_d)$  on a Hilbert space $\CL,$ define
\begin{eqnarray*} {\mathcal C}(\uR ):&=&
\{ \CM \subset \CL: \CM ~~\mbox {is an invariant subspace
 for each}~ R^*_i, \\
&&(p_\eta (\uR))^*h=0,\; \forall h\in \CM,\; \forall \eta \in J \}.
\end{eqnarray*}
${\mathcal C}(\uR)$ is a complete lattice with respect to arbitrary
intersections and  span closures of arbitrary unions. The maximal
element is called {\em maximal $J$-constrained subspace} and is
denoted by $\CL ^J(\uR)$ or $\CL^J.$ The row contraction $\uR
^J=(R_1^J, \ldots , R_d^J)$ obtained by compressing $\uR$ to $\CL^J$
is called {\em maximal $J$-constrained piece}. Clearly the maximal
$J$-constrained piece is $J$-constrained row contraction. The block
form of $\uR$ in terms of the maximal constrained piece is: $R_i =
\left(
\begin{array}{cc}
R^J_i  & 0 \\
\tilde{R}_i  & R^N_i \\
\end{array} \right)$
where $\uR^N$ is the compression of $\uR$ to the orthogonal complement of $\CL^J.$
\\

The three special  sets of polynomials $\{p_\eta\}_{\eta \in J}$
inducing constraints are:
\begin{enumerate}
\item $p_{i,j}(\uz)=z_i z_j -z_j z_i$ for $(i,j) \in \{1,2, \ldots, d\}^2$ are associated with commuting case.
\item $p_{i,j}(\uz)=z_i z_j - q_{ji}z_j z_i$ for $(i,j) \in \{1,2, \ldots, d\}^2$ where $|q_{ij}|=1$ and $q_{ij}=
q^{-1}_{ji}.$
\item $p_{i,j}(\uz)=z_i z_j - a_{ij}z_i z_j$ for $(i,j) \in \{1,2, \ldots, d\}^2$ where $A=(a_{ij})_{d \times d}$
is a $0-1$-matrix and each row and column has at least one non zero
entry.
\end{enumerate}
If the given set of polynomials is a combination of the above three,
we call such cases as $J^\prime$-constrained. The $J$ in the above
stated commuting case is denoted by $J^s.$ For
$J^\prime$-constrained row contractions
of the above mentioned type 1, 2 and 3, the constrained dilations
are  called standard commuting dilation, q-commuting dilation and
Cuntz-Krieger dilation respectively. Note that the type 2 includes the type 1 case.

We recall the definition of Cuntz-Krieger  algebras $\CO_A$. Let
$A=(a_{ij})_{d\times d}$ be a square $0-1$-matrix i.e. $a_{ij} \in
\{ 0 , 1\}$ and such that each row and column has at least one
non-zero entry.
\\

\noindent {\bf Definition ~~} {\em  $\CO_A$ is the universal
$C^*$-algebra generated by $d$ partial isometries $s_1,
\cdots,s_d$ with orthogonal ranges satisfying
\begin{equation}
\begin{array}{c}
s^*_is_i= \sum^d_{j=1} a_{ij} s_js^*_j\\
\\
\eins=\sum^d_{i=1}s_is^*_i.
\end{array}
\end{equation}
}

\section{Minimal constrained dilation}

In this section we assume $d < \infty.$  Lemma 5.1.3 is the only
exception where this assumption is not needed. Consider the maximal
$J^\prime$-constrained piece $\uL^{J^\prime}$ of the creation
operators $\uL$ on the Fock space $\Gamma (\C^d)$ or $\Gamma.$ It is
not difficult to show that the unital $C^*$-algebras generated by
$\uL^{J^\prime}$ satisfy
\[C^*(\uL^{J^\prime})= \overline{\mbox{span}} \{ L^{J^\prime}_\alpha
(L^{J^\prime}_\beta)^*: \alpha, \beta \in \tilde{\Lambda} \}.\]
Using this and the Poisson kernel (2.6) one can show further for any
given row contraction $\uT=(T_1,\ldots,T_d)$ on $\CH$ that there
exist unique unital completely positive map $\psi:C^*(\uL) \to
B(\CH)$  satisfying:
\[\psi(L_\alpha (L_\beta)^*)= T_\alpha (T_\beta)^*.\]
If $\uT$ is $J^\prime$-constrained, then there also exist unique
unital completely positive map $\phi: C^*(\uL^{J^\prime}) \to
B(\CH)$ with
\begin{equation}
\phi(L^{J^\prime}_\alpha (L^{J^\prime}_\beta)^*)= T_\alpha (T_\beta)^*,
\hspace{.5cm} \alpha, \beta \in \tilde{\Lambda}.
\end{equation}
We can use a minimal Stinespring dilation $\pi_1: C^*(\uL^{J^\prime}) \to
B(\CH_1)$ of $\psi$ such that
$$
\phi (X)=P_{\CH}\pi_1(X)|_{\CH}~~~ \forall X \in C^*(\uL^{J^\prime})
$$
and $\overline{\mbox{span}}\{\pi_1(X)h: X \in C^*(\uL^{J^\prime}),
h\in \CH\}=\CH_1$. The tuple
$\tilde{\uS}=(\tilde{S}_1,\cdots,\tilde{S}_d)$ where
$\tilde{S}_i=\pi_1(L^{J^\prime}_i),$ is the {\em minimal
$J^\prime$-constrained dilation} of $\uT$ which is unique up to
unitary equivalence.

We recall a result from \cite{BBD04}, \cite{De07a} and \cite{BDZ06}
(cf. Appendix of this chapter):
\begin{thm}
Let $\uT$ be a $J^\prime$-constrained row contraction on $\CH.$ Then
$\uV^{J^\prime}$ (i.e., the maximal $J^\prime$-constrained piece of $\uV$) is
the minimal $J^\prime$-constrained dilation of $\uT.$
\end{thm}
From the easy observation that the maximal constrained subspace of
$\uT$ is a $\uV$-coinvariant subspace we have $\uV$ to be an
isometric lifting of $\uT^J.$ Therefore using the previous Theorem
we get that the compression of $\uV$ to a
$\uV^{J^\prime}$-coinvariant subspace of $\hat{\CH}^{J^\prime}$ is
the constrained dilation of $\uT^{J^\prime}.$ It is natural to ask
if $\uV^{J^\prime}$ is the minimal constrained dilation of
$\uT^{J^\prime}.$ For $*$-stable row contraction a necessary and
sufficient condition appears in Theorem 9 of \cite{BBD04}. A version
for coisometric case is obtained here in Theorem 5.1.4 and another
for finite dimensional Hilbert spaces in Theorem 5.1.7.

We will need the following lemma. Let $\CG$ be the (non-selfadjoint unital)
weak operator topology-closed algebra generated by the
$V_i \in B(\hat{\CH})$ of the mid $\uV$ of $\uT.$

\begin{lem}
(Lemma 3.4 of \cite{DKS01})  $\CG[\CL]$ reduces $\CG$ if $\CL$ is a
$T^*_i$-invariant subspace.
\end{lem}

\begin{lem}
Suppose  $\uT,$ given on a finite dimensional $\CH,$ is coisometric.
Let $\CM$ be a subspace of (the dilation space) $\hat{\CH}$ which is
invariant for both $\CG$ and $\CG^*$ (i.e., reducing). Denote
$\CM\cap\CH$ by $\CH_\CM.$ Then $\CM=\CG[\CH_\CM].$
\end{lem}

\begin{proof} Note that $\CH_\CM$ is invariant with respect to $T^*_i$ for $ i=1, \ldots, d$
because $\CH$ and $\CM$ are $V^*_i$-invariant and $V^*_i|_\CH =
T^*_i$. Lemma 5.1.2 shows that $\CG[\CH_\CM]$ reduces $\CG$. Because
$\CH_\CM \subset \CM$ and $\CM$ is $\CG$-invariant we also have
$\CG[\CH_\CM] \subset \CM$. Let us assume that $\CH^\prime=\CM\ominus
\CG[\CH_\CM]$ is non-zero. But Corollary 4.2 of \cite{DKS01} gives
that any non-zero $\CG^*$-invariant subspace intersects $\CH$
nontrivially. Hence $\CH^\prime$ has a non-trivial intersection with
$\CH.$ This is a contradiction as $(\CM \ominus \CG[\CH_\CM])
\,\cap\, \CH = \CH_\CM \,\cap\, \CG[\CH_\CM]^\perp = \{0\}$.
Therefore $\CM=\CG[\CH_\CM].$ \end{proof}

Suppose $\pi $ is a representation on a Hilbert space $\CL$ of the
Cuntz algebra $\CO_d$ with generators $g_1,\ldots,g_d$.
The representation $\pi $ is said to be {\em spherical\/ } if $~\overline
{\mbox {span}}~\{ \pi(g_\alpha)\, h:
h\in \CL ^s(\pi (g_1), \ldots , \pi (g_d)),\, \alpha \in \tilde {\Lambda } \}
=\CL.$

Let us now assume that $\CH$ is finite dimensional. Theorem 19 of
\cite{BBD04} states that the mid $\uV$ on $\hat{\CH}$ of $\uT$ on
$\CH$ can be decomposed as $\uV^0 \oplus \uV^1$ with respect to the
decomposition of $\hat{\CH}$ as $\hat{\CH}^0 \oplus \hat{\CH}^1$
into reducing subspaces where $\uV^0$ is associated to a spherical
representation of $\CO_d$ and $\uV^1$ has trivial maximal commuting
piece. Because $\CH$ is finite dimensional, the already mentioned
direct integral decomposition (\cite{BBD04}, Theorem 18) now tells
us that $\hat{\CH}^0$ can be further decomposed into irreducible
subspaces as $\hat{\CH}^0_1 \oplus \ldots\oplus\hat{\CH}^0_k$ for
some $k\in \N$. Let $\CH_j:=\CH \cap \hat{\CH}^0_j.$ We observe that
$\CH_j, j=1, 2, \ldots, k$, are non-zero disjoint $T^*_i$-invariant
subspaces for $ i=1, \ldots, d$ and $\CG[\CH_j]=\hat{\CH}^0_j$ for
$j=1, \ldots, k$ from Lemma 5.1.3. It follows also that the
compressions of $\uT$ to the $\CH_j$'s are coisometric. But as the
restriction of $\uV$ to $\hat{\CH}^0_j$ is associated to an
irreducible and spherical representation, the related maximal
commuting subspace is one dimensional (cf. \cite{BBD04}, Theorem 18
and 19) and hence is a minimal $\CG^*$-invariant subspace for each
$j$. By Lemma 5.8 of \cite{DKS01} such a minimal $\CG^*$-invariant
subspace is unique, and since the $\CH_j$'s are $\CG^*$-invariant
subspaces, it follows that the maximal commuting subspace of $\uV$
on $\hat{\CH}^0_j$ is contained in $\CH_j.$

Consider the case when the maximal commuting subspace of the mid $\uV$
of a row contraction $\uT$ on the Hilbert space $\CH$ is contained in $\CH.$
Proposition 7 of \cite{BBD04} yields that $\hat{\CH}^{J^s} \cap \CH= \CH^{J^s}.$
So the maximal
commuting piece of $\uT$ is also the maximal commuting piece of $\uV$ and
therefore the standard commuting dilation of itself.

\begin{thm}
Suppose the dimension of $\CH$ is finite and $\uT$ is a coisometric row contraction on it.
Then the maximal commuting subspace of $\uV$ is contained in $\CH$
and coincides with the maximal commuting subspace of $\uT$.
\end{thm}

\begin{proof} Let $\uV$ on $\hat{\CH}$ be decomposed as above. From the arguments above,
we obtain that the maximal commuting subspaces of the compressions
of $\uV$ on $\hat{\CH}^0_j$ are contained in $\CH_j.$ The linear
span of all these subspaces is in fact $\hat{\CH}^{J^s} (\uV)$ and
hence is also contained in $\CH.$ The argument for the second
assertion has already been given above also. \end{proof}

For a lifting $\uE$ on $\CH_E=\CH_C \oplus \CH_A$ of $\uC$ by $\uA,$
it is evident that when $\uE$ is $J$-constrained, then both $\uC$
and $\uA$ are $J$-constrained. In case $\uE$ is a subisometric
lifting of a row contraction $\uC,$ then the unitary equivalence of
mids  $\uV^C$ and $\uV^E$ implies unitary equivalence of the maximal
constrained pieces of $\uV^C$ and $\uV^E.$ In addition when $\uE$ is
$J^\prime$-constrained, by Theorem 5.1.1  we obtain unitary equivalence of
their minimal constrained dilations.

\begin{cor} Let $\uE,$ on a finite dimensional Hilbert space, be a coisometric
lifting of $\uC$ by a $*$-stable row contraction $\uA.$ Then the maximal commuting pieces of $\uE$ and $\uC$
coincide.
\end{cor}
\begin{proof} Here $\uE$ is a subisometric lifting of $\uC.$ So maximal commuting pieces
of $\uV^C$ and $\uV^E$ are unitarily equivalent. By the previous
Theorem this means that the maximal commuting pieces of $\uE$ and
$\uC$ coincide. \end{proof}

\begin{thm}
Suppose $\uT$ is coisometric  and $\uT^N$ is $*$-stable. Then the minimal $J^\prime$-constrained dilation
of a $\uT^{J^\prime}$ is the maximal constrained piece of the mid $\uT.$
\end{thm}
\begin{proof} As $\uT^N$ is $*$-stable, by Proposition 4.2.3  $\uT$ is a subisometric lifting of $\uT^{J^\prime}.$
This causes the maximal $J^\prime$-constrained pieces of mids of $\uT$ and
$\uT^{J^\prime}$ to be unitarily equivalent. Moreover the maximal
constrained piece of mid of $\uT^{J^\prime}$ is the minimal constrained
dilation of $\uT^{J^\prime}.$  Hence the Theorem follows. \end{proof}

The maximal commuting piece of any coisometric row contraction $\uT$ on $\CH$ consists of subnormal operators. This is
because
the maximal commuting piece of a coisometric row contraction is commuting and coisometric, and so it has a standard
commuting dilation consisting of normal operators (cf. \cite{Ar98}, Corollary 1 in Section 8).
 Consequently if $\CH$ is finite dimensional, the maximal commuting piece of $\uT$ always
consists of normal operators.

\begin{cor} A commuting coisometric row contraction on a finite dimensional Hilbert space cannot be a
lifting of another commuting row contraction by a $*$-stable tuple.
\end{cor}
\begin{proof} Assume that $\uE$ is a commuting coisometric lifting of a
commuting row contraction $\uC$ on a finite dimensional Hilbert
space by $*$-stable $\uA.$ Then again by Proposition 4.2.3  $\uE$ is
a subisometric lifting of $\uC.$ Consequently, standard commuting
dilations of $\uE$ and $\uC$ are unitarily equivalent. For normal
coisometric row contractions their standard commuting dilations
coincide with the row contractions as shown in Theorem 15
\cite{BBD04}. As $\uC$ consists of normal operators, its standard
commuting dilation is equal to itself. In other words $\uC$ is equal
to the standard commuting dilation of $\uE$ along with the fact that
$\uE$ is a compression of the standard commuting dilation of $\uE.$
This yields that $\uC=\uE.$ \end{proof}

\begin{prop} Let $C$ be a (single) coisometry on a finite dimensional Hilbert space.
Let $E$ be a  coisometric lifting of $C$ by $A.$ Then $A$ is
coisometric.
\end{prop}
\begin{proof}
 Let
\[E =\left(
\begin{array}{cc}
C  & 0 \\
B  & A \\
\end{array} \right).\]
$C$ is infact a unitary. We have $C^*C=\eins$ and $BC^*=0.$ This
implies $B$ is zero and hence $A$ is coisometric. \end{proof}

The following is a passing by remark on the ergodic case treated in
Chapter 3.
\begin{rem}
If the maximal commuting subspace of a coisometric row contraction $\uT$
is of dimension greater than one then $\uT$ is not ergodic.
\end{rem}

\begin{proof} If $\uT$ is ergodic then by Theorem 5.1 of \cite{BJKW00} its
mid $\uV$ on $\hat{\CH}$ is associated to an irreducible
representation of $\CO_d$ and so Theorem 19 of \cite{BBD04} tells us
that $\uV$ has trivial or one dimensional maximal commuting piece.
But as $\hat{\CH}^{J^s} \cap \CH= \CH^{J^s}$ by Proposition 7 of
\cite{BBD04}, we finally get that $\uT$ has trivial or one
dimensional maximal commuting piece. \end{proof}

\section{Constrained characteristic function}
As before let  $\uC$  be a row contraction on $\CH_C$ and $\uE$ on
$\CH_E = \CH_C \oplus \CH_A$ be a contractive lifting of $\uC$ with
$B_i=P_{\CH_A} E_iP_{\CH_C}.$

For a reduced lifting if we define $M:= P_{\CK^\bot} W P_{\CH_E}:
\CH_E \to \CH_C \oplus \Gamma (\C^d) \otimes \CD_C$ then $M$ acts as
identity on $ \CH_C$ and
\[Mh := \sum_\alpha e_\alpha \otimes  \gamma \DsA A^*_\alpha h \mbox{~~~for~~}h \in \CH_A\]
where $\gamma: \CD_{*,A} \to \CD_C$ is a contraction as defined
before. The characteristic function for this reduced lifting is
defined as
\[M_{C,E}:= P_{\Gamma(\C^d) \otimes \CD_C}W |_{\Gamma(\C^d) \otimes \CD_E}.\]

In order to have a similar invariant for $J$-constrained row
contractions, say $\uT$ on $\CH,$ we consider a dilation $\uS^T$ of
 $\uT$ obtained by a compression of
  Popescu's realization of mid to $\CH \oplus (\Gamma_J (\C^d) \otimes
\CD_T).$ Explicitly for $i= 1,\ldots, d$
\[ S^T_i (h \oplus d) = T_i h \oplus [ e_0 \otimes (D_T)_i  h + (L^J_i \otimes \eins) (d)]
\mbox{~~~for~~} h \in \CH, d \in \Gamma_J (\C^d) \otimes \CD_T.
\] That $\uS^T$ is a dilation is immediate and we call it
{\em pseudo $J$-constrained dilation} of $\uT.$ Note that $\uS^T$
does not satisfy constrained relations. In case of $J^\prime$ this
dilation contains the minimal constrained dilation as a corner
(w.r.t. a coinvariant subspace). Indeed Theorem 5.1.1 tells us that
the maximal $J^\prime$-constrained pieces of mids  are minimal
constrained dilations. The constrained pieces are clearly on spaces
$\CH \oplus \CN_J$ where $\CN_J$ are some proper coinvariant
subspaces of $(\Gamma_J (\C^d) \otimes \CD_T)$ with respect to the
$\uL^J \otimes \eins.$
\\

A {\em constrained lifting} of a row contraction is a lifting which
as a tuple of operators is $J$-constrained. For a constrained
lifting $\uE$ of $\uC,$ the tuple $\uA$ is also $J$-constrained and
therefore by \cite{Po05} $M \CH_A$ is contained in $\Gamma_J (\C^d)
\otimes \CD_C.$ We set $M_0:=M|_{\CH_A}$ and $M^J_0:= P_{\Gamma_J
(\C^d) \otimes \CD_C} M_0.$ Let us denote the restriction of $W$ to
$\CH_E \oplus (\Gamma_J (\C^d) \otimes \CD_E)$ by $\W.$ Because
$M_{C,E}$ is multianalytic and
\[ \Gamma_J (\C^d)=\overline{\mbox{span}}\{ L^*_\alpha (L_i L_j -L_j L_i)^* L^*_\beta f =0 : f \in \Gamma (\C^d);
\alpha, \beta \in \tilde{\Lambda}\} \]
(cf. \cite{BDZ06}), we get for $k \in \Gamma_J (\C^d) \otimes \CD_E$
\begin{eqnarray*} [(L^*_\alpha (L_i L_j -L_j L_i)^* L^*_\beta) \otimes \eins] M_{C,E} k &=& M_{C,E} [(L^*_\alpha
(L_i L_j -L_j L_i)^* L^*_\beta) \otimes \eins] k \\
&=&  0 \mbox{~~ for } \alpha, \beta  \in \tilde{\Lambda}.
\end{eqnarray*}
Thus the range of $P_{\Gamma (\C^d) \otimes \CD_C}\W |_{\Gamma_J
(\C^d) \otimes \CD_E}$ is $\Gamma_J (\C^d)\otimes \CD_C.$ We define
the {\em constrained characteristic function} as
$M_{J,C,E}:=P_{\Gamma_J (\C^d) \otimes \CD_C}\W |_{\Gamma_J (\C^d)
\otimes \CD_E}$ which is the same as $M_{C,E}|_{\Gamma_J (\C^d)
\otimes \CD_E}.$ Let $d^i_h := D_E (e_i \otimes h)$ for $h \in
\CH_E.$ From equations (4.49) and (4.50) it is clear that the symbol
$\theta_{J,C,E}$ of $M_{J,C,E}$ for $h \in\, \CH_C$ will be given by
\begin{eqnarray*}
\theta_{J,C,E}  d^i_h &:=& e_0 \otimes [\CD_C( e_i
\otimes h) - \gamma  \DsA  B_i h] - \sum_{|\alpha|\geq 1} L^J_\alpha
e_0 \otimes \gamma \DsA A^*_\alpha B_i h
\end{eqnarray*}
and for $h \in\, \CH_A$ by
\begin{eqnarray*}
\theta_{J,C,E} d^i_h &:=& - e_0 \otimes \gamma \DsA A_i h +
\sum_\alpha L^J_i L^J_\alpha e_0 \otimes \gamma \DsA
A^*_\alpha (\eins - A^*_i A_i) h \\
&& + \sum_{j \neq i} \sum_\alpha L^J_j L^J_\alpha e_0 \otimes
\gamma \DsA A^*_\alpha (-A^*_j A_i) h.
\end{eqnarray*}
It is evident that $\uS^E$ is unitarily equivalent to $\uS^C \oplus
\uZ$ where $\uZ$ is defined on some Hilbert space $\CK_1$ and
\[ \W S^E_i = (S^C_i \oplus Z_i) \W. \]

Let $\CH^1_A$ be
\[\CH^1_A:=\{ h \in \CH_A: \sum_{|\alpha|=n} \|A^*_\alpha h\|^2 = \|h\|^2 \mbox{~for all~} n \in \N \}\]
as in the last chapter. The equation (4.36)  gives us
\[ W \CH^1_A \hspace{.5cm}
\bot \hspace{.5cm} \CH_C \oplus (\Gamma (\C^d) \otimes \CD_C).\]
Therefore
\[ \W \CH^1_A \hspace{.5cm} \bot \hspace{.5cm} \CH_C
\oplus (\Gamma_J (\C^d) \otimes \CD_C) \] and hence $\W \CH^1_A
\subset \CK_1.$ We denote the space on which $S^E_i$s are defined by
$\bCK_J.$ Assuming $\gamma$ to be resolving and repeating the
arguments from the last chapter in this setting, we have
\[ \W \bCK_J \ominus (\CH_C \oplus \CH^1_A) = (\Gamma_J (\C^d)\otimes \CD_C) \bigvee \W (\Gamma_J (\C^d)\otimes \CD_E). \]

In case of reduced constrained liftings, $\CH^1_A$ is trivial and
taking $\Delta_{J,C,E} := (\eins -
M^*_{J,C,E}M_{J,C,E})^{\frac{1}{2}},$  we can realise $\uZ$ on
$\overline{ \Delta_{J,C,E} (\Gamma_J (\C^d) \otimes \CD_E)}.$ Hence
for such liftings $\uS^C \oplus \uZ$ is a row contraction on
\[ \W \bCK_J := \CH_C \oplus (\Gamma_J (\C^d) \otimes \CD_C) \oplus \overline{\Delta_{J,C,E} (\Gamma_J (\C^d)
\otimes \CD_E)} \] and
\[ \W \CH_E := \W \bCK_J \ominus (M_{J,C,E} \oplus \Delta_{J,C,E}) (\Gamma_J (\C^d) \otimes \CD_E). \]
is coinvariant for $\uS^C \oplus \uZ$ and their restriction to $\W
\CH_E$ gives a copy of $\uE.$

\begin{thm}
For a constrained row contraction, two constrained reduced liftings
are unitarily equivalent if and only if its constrained
characteristic functions coincide.
\end{thm}

The proof follows from Theorem 4.3.7 after noting that unitary
equivalence of two liftings say $\uE$ and $\hat{\uE}$ of given row
contraction $\uC$ will also mean unitary equivalence of $\uL \otimes
\eins_{\CD_E}$ and $\uL \otimes \eins_{\CD_{\hat{E}}}.$ This implies
unitary equivalence of the pseudo constrained dilations.
\\

Suppose $\uE$ on $\CH_E=\CH_C \oplus \CH_A$ be such a
$J$-constrained lifting of $\uC$ by a c.n.c. row contraction $\uA.$
We follow the notations of the last Section. It is immediate that
\[ \W (M^J_0)^* v = P_{\W \CH_E} (v \oplus 0), \mbox{~~~~~~for~~} v \in \Gamma (\C^d) \otimes \CD_C.\]
Consequently
\begin{equation}
\| (M^J_0)^* v \| = \| P_{\W \CH_E} (v \oplus 0) \|, \mbox{~~~~~~for~~} v \in \Gamma (\C^d) \otimes \CD_C.
\end{equation}

The map $\phi: \Gamma_J (\C^d) \otimes \CD_E \to \W \bCK_J$  (in fact $\W \bCK_J \ominus \CH_C$) defined by
\[  u \mapsto  M_{J,C,E} u \oplus \Delta_{J,C,E} u  \]
is an isometry and
\[ \phi^* (v \oplus 0) = M^*_{J,C,E} v \]
for $v \in \Gamma_J (\C^d) \otimes \CD_C.$ Suppose $P_{\W \CH_E} \in
B(\bCK_J)$ is the orthogonal projection  onto $\W \CH_E.$ Then for
$v \in \Gamma_J (\C^d) \otimes \CD_C$
\begin{equation}
 \| v\|^2  = \|P_{\W \CH_E} (v \oplus 0) \|^2 + \| \phi \phi^* (v \oplus 0)\|^2
= \|P_{\W \CH_E} (v \oplus 0) \|^2 + \|M^*_{J,C,E} v \|^2.
\end{equation}

Comparing equations (5.3) and (5.4) we get
\begin{equation}
M^J_0 (M^J_0)^*  + M_{J,C,E} M^*_{J,C,E} = \eins_{\Gamma_J (\C^d)\otimes \CD_C}.
\end{equation}
When $J=\{ 0 \}$ we obtain
\begin{equation}
M_0 (M_0)^*  + M_{C,E} M^*_{C,E} = \eins_{\Gamma (\C^d)\otimes \CD_C}.
\end{equation}

\begin{rem}
The above results clearly holds for any contractive liftings which
are resolving, i.e., the liftings need not necessarily be reduced.
\end{rem}

\section{Appendix}

In Lemma 5.1.1 we refer to a result from \cite{BBD04}, \cite{De07a}
and \cite{BDZ06}. We illustrate how one can prove it in a special
case in this Appendix.

\begin{thm}
Let $\uT=(T_1,\ldots,T_d)$ be a  row contraction on a Hilbert space $\CH$ where $T_i$'s are mutually commuting. Then
  the maximal commuting piece of the mid  of $\uT $ is a realization of the
standard commuting dilation of $\uT$.
\end{thm}

It is known that the maximal $J^s$-constrained subspace (defined in
Section 5.1) of $\Gamma (\C^d)$ w.r.t. the standard tuple $\uL$ of
creation operators is the symmetric Fock space $\Gamma_s(\C^d)$. Let
us denote by $\uS$ the  constrained piece $\uL^{J^s}.$ To prove the
above theorem we consider the map $\phi: C(\uS) \to B(\CH)$ from the
Equation (5.2) associated with the standard commuting dilation of
$\uT $  and the dilation $\pi_1$ of $\phi$ on a Hilbert space
$\CH_1.$ As $\uS$ is also a row contraction, out of the discussion
at the beginning of Section 5.1 it follows that there is a unique
unital completely positive map $\eta : C^*(\uL)\to C^*(\uS)$,
satisfying
$$\eta (L_\alpha (L_\beta)^*)= S_\alpha (S_\beta )^*, ~~~~ \alpha , \beta \in
{\tilde \Lambda}.$$
Consider the minimal Stinespring dilation
of the composed map $\pi _1\circ \eta : C^*(\uL)\to \CB(\CH _1)$. We obtain a Hilbert space
 $\CH _2 $ containing  $\CH _1$ and a unital $*$-homomorphism
$\pi _2 : C^*(\uL)\to \CB(\CH _2),$ such that
$$\pi _1\circ \eta (X)=P_{\CH _1 }\pi _2(X)|_{\CH _1}, ~~~~ X\in C^*(\uL),$$
and $\overline {\mbox {span}}~\{ \pi _2(X)h: X\in C^*(\uL), h\in \CH _1\}=
{\CH _2}.$
In the following commuting diagram

\hskip1.5in \begin{picture}(200,115)

\put(0,20){$C^*(\uL)$}

\put(40,21){$\longrightarrow$}

\put(70,20){$C^*(\uS)$}

\put(110,21){$\longrightarrow$}

\put(140,20){$\mathcal{B}(\CH)$}

\put(140,60){$\mathcal{B}(\CH_1)$}

\put(140,100){$\mathcal{B} (\CH_2)$}

\put(40,30){\vector(4,3){80}}

\put(110,30){\vector(4,3){20}}

\put(150,40){$\downarrow $}

\put(150,80){$\downarrow $}

\put(45,10){$\eta$}

\put(115,10){$\phi$}

\put(110,40){$\pi_1$}

\put(70,70){$\pi_2$}

\end{picture}

\noindent  horizontal arrows are unital completely positive maps, down arrows are
compressions and diagonal arrows are minimal Stinespring dilations.

Taking ${\hat  \uV}=({\hat V_1}, \ldots , {\hat V_d})
= (\pi _2(L_1), \ldots , \pi _2(L_d))$, we first show that ${\tilde {\uS }}
=(\pi _1(S_1), \ldots , \pi _1(S_d))$ is the
maximal commuting piece of $\hat {\uV}$. Then we prove that
 ${\hat \uV}$ is the mid of $\uT.$  The last statement follows if we can show that $\pi _2$ is
a minimal dilation of $\phi \circ\eta$  because
the minimal Stinespring dilation is unique up to unitary equivalence.
 For this we would need the a lemma (cf. Theorem 15 in \cite{BBD04}).
\begin{defn}
A $d$-tuple $\uT =(T_1, \ldots , T_d) $ of operators on a Hilbert
space $\CH $ is called a spherical unitary if it is commuting, each
$T_i$ is normal, and $T_1T_1^*+\ldots +T_dT_d^*=\eins.$
\end{defn}

\begin{lem}
Let $\uT$ be a spherical unitary on a Hilbert space $\CH.$ Then
the maximal commuting piece of the mid of
$\uT $ is $\uT $.
\end{lem}

\noindent {\em Proof of Theorem 5.3.1:} As $C^*(\uS)$ contains the
ideal of all compact operators by standard $C^*$-algebra theory we
have a direct sum decomposition of $\pi _1$ as follows. Take
${\mathcal H}_1 = {\mathcal H}_{1C} \oplus {\mathcal H}_{1N}$\\
where ${\mathcal H}_{1C} = \overline{\mbox{span}}\{\pi _1(X)h : h\in {\mathcal
H},
 X \in
C^*(\underline{S})$ and $X$ is compact$\}$ and
 ${\mathcal H}_{1N} = {\mathcal H}_1 \ominus {\mathcal H}_{1C}$,
Clearly ${\mathcal H}_{1C}$ is a reducing subspace for $\pi _1$.
Therefore \[\pi _1(X)=\begin{pmatrix}
  \pi_{1C}(X)&  \\
     & \pi_{1N}(X)
\end{pmatrix}     \]
that is, $\pi _1=\pi _{1C}\oplus \pi _{1N}$ where $\pi_{1C}(X) = P_{{\mathcal H}_{1C}} \pi _1(X) P_{{\mathcal H}_{1C}}$,
$\pi_{1N}(X) = P_{{\mathcal H}_{1N}} \pi _1(X) P_{{\mathcal H}_{1N}}$.
 $\pi _{1C}(X)$ is just the identity representation
with some multiplicity as remarked in \cite{Ar98}. In other words
$\CH _{1C}$ can be factored as $\CH _{1C}= \Gamma _s (\C^d) \otimes
\CD_{*,T},$ such that $\pi _{1C}(X)=X\otimes \eins$. Also $\pi
_{1N}(X)=0$ for compact $X$. Therefore, taking $Z_i=\pi _{1N}(S_i)$,
$\underline {Z}=(Z_1, \ldots , Z_d)$ is a spherical unitary.

Now  $\pi _1\circ \eta =(\pi _{1C}\circ \eta )\oplus (\pi _{1N}\circ \eta)$
and
the minimal Stinespring dilation of a direct sum of two completely
positive maps is the direct sum of minimal Stinespring dilations. So $\CH _2
$  decomposes as $\CH _2= \CH _{2C}\oplus \CH _{2N},$ where
$\CH _{2C}, \CH _{2N}$ are orthogonal reducing subspaces of $\pi _2$,
such that $\pi _2$ also decomposes, say $\pi _2=\pi _{2C}\oplus
\pi _{2N},$ with
$$\pi _{1C}\circ \eta (X)= P_{\CH _{1C}}\pi _{2C}(X)|_{\CH _{1C}}, ~~
\pi _{1N}\circ \eta (X)= P_{\CH _{1N}}\pi _{2N}(X)|_{\CH _{1N}},$$
for $
X\in C^*(\uL)$
 with
$\CH _{2C}= ~\overline {\mbox {span}}~\{ \pi _{2C}(X)h: X\in
C^*(\uL), h\in \CH _{1C}\}$ and $\CH _{2N}= ~\overline {\mbox
{span}}~\{ \pi _{2N}(X)h: X\in C^*(\uL), h\in \CH _{1N}\}$.  It is
also not difficult to see that $\CH _{2C}= ~\overline {\mbox
{span}}~\{ \pi _{2C}(X)h: X\in C^*(\uL), X ~~\mbox {compact},~  h\in
\CH _{1C}\}$ and hence $\CH _{2C}$ factors as $\CH _{2C}= \Gamma
(\C^d) \otimes \CD_{*,T}$ with $\pi _{2C}(V_i)=V_i\otimes \eins.$
Also $(\pi _{2N}(L_1), \ldots , \pi _{2N}(L_d))$ is a mid of
spherical unitary $(Z_1, \ldots , Z_d)$.
Hence by Lemma 5.3.3 and some easy observations we get that $(\pi
_1(S_1), \ldots , \pi _1(S_d))$ acting on $\CH _1$ is the maximal
commuting piece of $(\pi _2(L_1), \ldots , \pi _2(L_d))$.

Next we show  that $\pi _2$ is the minimal Stinespring
dilation of $\phi \circ \eta .$ Assume that this is not true.
Then we get a reducing subspace $\CH _{20}$ for $\pi _2$ by
taking
$\CH _{20} =~\overline {\mbox {span}}~\{\pi _2(X)h: X\in C^*(\uL), h\in \CH
\} .$
Take $\CH _{21}=\CH _2\ominus \CH _{20}$ and correspondingly
decompose $\pi _2$ as $\pi _2= \pi _{20}\oplus \pi _{21}$,
\[\pi _2(X)=\begin{pmatrix}
  \pi_{20}(X)&  \\
     & \pi_{21}(X)
\end{pmatrix}     \]

 Note
that we already have $\CH \subseteq \CH _{20}.$  We claim that $\CH
_2\subseteq \CH _{20}.$ Firstly, as $\CH _1$ is the space where the
maximal commuting piece of $(\pi _2(L_1), \ldots ,$ $ \pi _2(L_d))=
(\pi _{20}(L_1)\oplus \pi _{21}(L_1), \ldots , \pi _{20}(L_d)\oplus
\pi _{21}(L_d))$ acts, $\CH _1$ decomposes as $\CH _1=\CH
_{10}\oplus \CH _{11}$ for some subspaces $\CH _{10} \subseteq \CH
_{20},$ and $\CH _{11}\subseteq  \CH _{21}$.  So for $X\in
C^*(\uL)$, $P_{\CH _1}\pi _2(X)P_{\CH _1}$, has the form (see the
diagram)
$$ P_{\CH _1}\pi _2(X)P_{\CH _1}=
\begin{pmatrix}
  \pi _{10}\circ \eta (X)&0& & \\
  0 & 0& & \\
   & & \pi_{11}\circ \eta(X)& 0 \\
   & & 0 & 0
\end{pmatrix}
$$
where $\pi _{10}, \pi _{11}$ are compressions of $\pi _1$ to
$\CH _{10}$, $\CH _{11}$ respectively. As the mapping $\eta $ from $C^*(\uL)$
to $C^*(\uS)$ is clearly surjective, it follows that $\CH _{10},
\CH _{11}$ are reducing subspaces for $\pi _{1}.$ Now as $\CH $
is contained in $\CH _{20}$, in view of minimality of $\pi _1$
as a Stinespring dilation, $\CH _1\subseteq \CH _{20}$. But then
the minimality of $\pi _2$ shows that $\CH _2\subseteq \CH _{20}$.
So finally we get $\CH _2=\CH _{20}.$
\qed

\chapter{Modules}

In this chapter we focus on two module structures arising in
operator algebras. By continuing the analysis from the previous
chapters we derive some interesting results for modules of both
types.

For commuting row contractions Arveson introduced the concept of
Hilbert module in \cite{Ar03} where he also gave some geometric
invariants for these modules. These notions where extended to the
noncommutative case by Kribs (\cite{Kr01}) and Popescu (\cite{Po02})
independently. We obtain how these invariants of Hilbert modules can
be expressed in terms of our characteristic functions in Section 6.1
(cf. \cite{Po05}).  The subsequent section is devoted to computing
some examples. In the final section we are concerned with the second
module structure called Hilbert $C^*$-module. We generalize part of
our theory developed in chapter 4 for liftings of covariant
representations of such modules.

\section{Invariants of Hilbert modules}

Let $\uT=(T_1,\ldots,T_d)$ be a mutually commuting $d$-tuple of
operators on $\CH.$ A {\em Hilbert module} structure on $\CH$ over
the algebra of polynomials in $d$-variables $\CP=\C[z_1,\ldots,z_d]$
is obtained by setting
\[ g . h := g(T_1, \ldots, T_d ) h, \mbox{~~for ~~} g  \in \CP, h \in \CH_E.  \]
This Hilbert module is contractive in the sense that:
\[ \| z_1. h_1 + \ldots + z_d . h_d \|^2 \leq \|h_1\|^2 + \ldots + \|h_n\|^2, \hspace{1cm} h_1,\ldots, h_d \in \CH.\]
We caution the reader that this notion of Hilbert module is
different from the ones discussed in Sections 2.4 and 6.3.
\\

We will further assume that $D_{*,T}$ has finite rank. Consider the function $T: \C^d \to B(\CH)$ given by
\[T(z)= \overline{z}_1 T_d + \ldots + \overline{z}_d T_d.\]
If $z$ belongs to the open unit ball $B_d$ of $\C^d,$ then it
follows immediately that $\|T\| < 1$ and $\eins - T(z)$ is
invertible. For every $z \in B_d$ let us form an operator $F(z)\in
B(D_{*,T}(\CH))$ by putting
\[F(z)h= D_{*,T} ((\eins - T(z))^*)^{-1} (\eins - T(z))^{-1} D_{*,T}h,  \mbox{~~~~~} h \in \CD_{*,T}. \]
 Arveson in \cite{Ar03} shows that with respect to the natural rotation-invariant probability
 measure $\sigma$ on $\partial B_d,$ the limit
\[K_0 (z)= \lim_{r \uparrow 1} (1 - r^2) \mbox{~trace~} F(rz)= 2. \lim_{r \uparrow 1} \mbox{~trace~} F(rz) \]
exists for $\sigma$-almost every $z \in \partial B_d$ and he defined
the {\em curvature invariant} of the Hilbert module as the scalar
\[ curv_s(T) := \int_{\partial B_d} K_0 (z) d \sigma (z).\]

Another invariant for Hilbert modules is the Euler characteristic
$\chi_s(T).$  For defining it we need to consider the submodule of
$\CH$:
\[\CM_T=\mbox{span} \{g.h: g \in \CP, h \in \CD_{*,T} \}.\]
If $r=\dim \CD_{*,T}$ and $k_1,\ldots,k_r$ is a basis for
$\CD_{*,T},$ then
\[\CM_T = \{ f_1. k_1 + \ldots + f_r. k_r : f_i \in \CP\}.\]
Thus $\CM_T$ is finitely generated. It follows from Hilbert's syzygy
theorem (cf. Theorem 182 of \cite{Ka70}) that $\CM_T$ has finite
free resolution, i.e., there is an exact sequence of $\CP$-modules
\[0 \to F_n \to \ldots \to F_2 \to F_1 \to \CM_T \to 0\]
where $F_k$ is a free module of finite rank $\beta_k.$ The numbers
$\beta_k$ are called {\em Betti numbers.} The {\em Euler
characteristic} is defined as
\[\chi_s (T)= \sum^n_{k=1} (-1)^{k+1} \beta_k \]
and is independent of the choice of the finite free resolution.

In the above quoted article of Arveson it is further shown that:
\[ curv_s T := (d-1)! \lim_{n \to \infty} \frac{ \mbox{trace} [K^*_{J^s}(T) (Q_{\leq n} \otimes \eins_{\CD_{*,T}} )
K_{J^s}(T) ]}{d^n}, \]
\[ \chi_s (T) :=  d! \lim_{n \to \infty} \frac{ \mbox{rank} [K^*_{J^s}(T) (Q_{\leq n} \otimes \eins_{\CD_{*,T}} )
K_{J^s}(T) ]} {n^d},  \] where  $Q_{\leq n}$ is the orthogonal
projection of the {\em symmetric Fock space} $\Gamma_s (\C^d)$ onto
$ P_{\Gamma_s (\C^d)} (\C \oplus (\C^d)^{\otimes 2} \oplus \ldots
\oplus (\C^d)^{\otimes n})$ and $K_{J^s} (T):=P_{\Gamma_s (\C^d)
\otimes \CD_{*,T}} K(T) $ is obtained by projecting the Poisson
kernel $K(T)$ on $\Gamma_s (\C^d) \otimes \CD_{*,T}.$ (The Poisson
kernel is denoted by $C$ in Equation (2.6). We are using a different
notation here for convenience).
\\

Take $\F^+_d$ to be the free semigroup with $d$ generators
$f_1,\ldots,f_d$ and denote the corresponding complex free semigroup
algebra  with $\C \F^+_d.$ In an analogous way as above, to any row
contraction (not necessarily commuting) $\uT=(T_1,\ldots,T_d)$ on a
Hilbert space $\CH,$ a contractive { \em Hilbert module} over $\C
\F^+ _d$ (cf. \cite{Po02}) can be associated through
\[ g . h := g(T_1, \ldots, T_d ) h, \mbox{~~for ~~} g  \in \C \F^+_d, h \in \CH.  \]
The contractivity of the module now means
\[ \| f_1. h_1 + \ldots + f_d . h_d \|^2 \leq \|h_1\|^2 + \ldots + \|h_d\|^2, \mbox{~~for ~~} h_1,\ldots, h_d \in \CH.\]

 We recall the definition of the
curvature invariant $curv T$  and Euler characteristic $\chi (T)$ of
Hilbert modules introduced by Kribs and Popescu  (cf.  \cite{Kr01},
\cite{Po02}):
\[ curv T := \lim_{n \to \infty} \frac{ \mbox{trace} [K^*(T) (P_{\leq n} \otimes \eins_{\CD_{*,T}} ) K(T) ]}{d^n}, \]
\[ \chi (T) :=  \lim_{n \to \infty} \frac{ \mbox{rank} [K^*(T) (P_{\leq n} \otimes \eins_{\CD_{*,T}} ) K(T) ]}
{1 + d + \ldots + d^{n-1}},  \] where $K(T)~(~=M_0~)$ is the Poisson
kernel of $\uT$ as before and $P_{\leq n}$ is the orthogonal
projection of $\Gamma (\C^d)$ onto $ \C \oplus  (\C^d)^{\otimes 2}
\oplus \ldots \oplus (\C^d)^{\otimes n}.$

These invariants are shown in this section to be related to
constrained reduced liftings  when $\gamma$ is an isometry. When
$\uE$ on $\CH_E=\CH_C \oplus \CH_A$ is such a reduced lifting of
$\uC$ by $\uA,$ we can get three Hilbert modules namely those
associated with $\uC,\uA$ and $\uE.$ Assume in the sequel that rank
$D_C$ is finite.

\begin{thm}
\[ curv A = \mbox{rank} D_C- \lim_{n \to \infty} \frac{\mbox{trace} [M_{C,E} M^*_{C,E}
(P_{\leq n} \otimes \eins )]}{d^n}, \]
\[ \chi (A) =  \lim_{n \to \infty} \frac{ \mbox{rank} [(\eins - M_{C,E} M^*_{C,E}) (P_{\leq n}
\otimes \eins )] }{1 + d + \ldots + d^{n-1}}.  \]
\end{thm}
\begin{proof} An easy simplification shows that for a lifting $\uE$
\begin{equation}
 curv A =  \lim_{n \to \infty} \frac{ \mbox{trace} [M_0^* (P_{\leq n} \otimes \eins_{\CD_C} )M_0 ]}{d^n},
\end{equation}
\begin{equation}
 \chi (A)  =   \lim_{n \to \infty} \frac{ \mbox{rank} [M_0^* (P_{\leq n} \otimes \eins_{\CD_C} )M_0 ]}
{1 + d + \ldots + d^{n-1}}.
\end{equation}
Evidently
\[ \mbox{trace} M_0^* (P_{\leq n} \otimes \eins )M_0 = \mbox{trace} M_0 M_0^* (P_{\leq n} \otimes \eins ).\]
Since rank $D_C < \infty$ and $M_0$ is injective on the range of
$M_0^* (P_{\leq n} \otimes \eins)$ we have
\[ \mbox{rank} M_0^* (P_{\leq n} \otimes \eins )M_0 = \mbox{rank} M_0^* (P_{\leq n} \otimes \eins )=
 \mbox{rank} M_0 M_0^* (P_{\leq n} \otimes \eins ).\]
Now from equations (5.6), (6.1) and  (6.2) the claim follows.
\end{proof}

Let us consider the case when $\uE$ is commuting. The  symmetric
Fock space $\Gamma_s (\C^d)$  can be identified with the space
$H^2,$ of all analytic functions on the open unit ball $B_d,$ the
reproducing kernel Hilbert space with the kernel $K_d:B_d \times B_d
\to \C$ defined by
\[ K_d (z,w):= \frac{1}{1- \langle z,w \rangle_{\C^d}}  \mbox{~~~for~~} z,w \in B_d\]
(see \cite{Ar03}). In this picture  $\uL^{J^s}$ corresponds to the
tuple $(M_{z_1},\ldots, M_{z_d})$ of multiplication operators  by
the coordinate functions. The constrained characteristic function
gets identified with a multiplication operator $M_{J^s,C,E}: H^2
\otimes \CD_E \to H^2 \otimes \CD_C$ with its symbol
$\theta_{J^s,C,E}$ as a $B(\CD_E, \CD_C)$-valued bounded analytic
function on $B_d.$ The $curv_s A$ simplifies using similar arguments
as for the proof of Theorem 6.1.1 to
\[ curv_s A = (d-1)! \lim_{n \to \infty} \frac{ \mbox{trace} [(M^{J^s}_0)^* (Q_{\leq n} \otimes \eins_{\CD_{*,A}} )
 M^{J^s}_0 ]}{d^n}, \]
From Theorem A of \cite{Ar03} and equation (5.5) we get
\[ curv_s A = \int_{\partial B_d} \lim_{r \to 1} \mbox{trace} [\eins - \theta_{J^s,C,E}(r \zeta)
\theta^*_{J^s,C,E}(r \zeta)] d \sigma (\zeta)\] and like in Theorem
6.1.1 we also have the final simplifications
\[ curv_s A = \mbox{rank} D_C- (d-1)!\lim_{n \to \infty} \frac{\mbox{trace} [\theta_{J^s,C,E} \theta^*_{J^s,C,E}
(Q_{\leq n} \otimes \eins )]}{d^n}, \]
\[ \chi_s (A) =  d! \lim_{n \to \infty} \frac{ \mbox{rank} [(\eins - \theta_{J^s,C,E} \theta^*_{J^s,C,E}) (Q_{\leq n}
\otimes \eins )] }{n^d}.  \]
\\

\section{Examples}

\noindent {\bf Example 1:}\\
Here we consider a coisometric lifting  for $d=1.$
Let
\[
E=\left(
\begin{array}{cc}
C  & 0 \\
B  & A \\
\end{array}
\right),
\]
on $\CH_E=  l^2(\N) \oplus l^2(\Z),$ be a lifting of $C$ on
$\CH_C=l^2(\N).$ We denote the standard basis for $l^2(\N)$ and
$l^2(\Z)$ by $\{e_i\}^\infty_{i=1}$ and $\{g_i\}^{\infty}_{-
\infty}$ respectively. Denote the shift operator on $l^2(\N)$ by
$S.$ Take operator $C=S^*.$  Further, fix  $0 \leq \lambda \leq 1$
and define $B$ and $A$ as follows:
\[ B( \sum^\infty_{i=1} a_i e_i ):=\sqrt{1- \lambda^2} \hspace{.2cm} a_1 g_1,\]
\[A( \sum^\infty_{-\infty} a_i g_i)=  \sum^\infty_{-\infty} c_i g_i,\]
where $c_1=\lambda a_0$ and $c_{i+1}=a_i$ for $i \neq 0.$

It is easy to verify that $CC^*=\eins,BC^*=0$ and
\[BB^*+ AA^*=\eins.\]
This implies  that $E$ is coisometric. Here $\eins - C^* C$ is a
projection onto $\C e_1.$ In this case $\CD_{*,A}=\C g_1$ and
$\CD_C= \C e_1.$ The isometry $\gamma: \CD_{*,A} \to \CD_C$ is given
by
\[\gamma (a_1 g_1) = a_1 e_1.\]

The characteristic function $M_{C,E}$ of the lifting in this case
maps $l^2(\N) \otimes \CD_E$ to $l^2(\N) \otimes \CD_C$ with symbol
$\theta_{C,E}.$ Let $\omega$ denote the unit vector $(1,0,0,
\ldots)$ of $l^2(\N)$ in $l^2(\N) \otimes \CD_C.$ Thus for $h=\sum^\infty_{i=1} a_i e_i  \in \CH_C =
l^2(\N)$
\[ \theta_{C,E} d_h= ( \lambda^2 a_1 e_1) \omega, \]
and for $h= \sum^\infty_{-\infty} a_i g_j \in \CH_A =  l^2(\Z)$
\[ \theta_{C,E} d_h= - (\lambda \sqrt{1- \lambda^2} \hspace{.2cm} a_0 e_1) \omega.\]
Note that here $M_{C,E}$ and $M_{J^s,C,E}$ will be the same.

We make use of the formula  from Theorems C and D of [Ar03] to
calculate
\[ curv_s A= \lim_{n \to \infty} \frac{\mbox{trace}(\eins - A^{n+1} (A^*)^{n+1})}{n}= \lim_{n \to \infty} \frac{(n+1)(1-\lambda^2)}{n}
=1-\lambda^2\]
\[ \chi_s (A)=  \lim_{n \to \infty}\frac{\mbox{rank}(\eins - A^{n+1} (A^*)^{n+1})}{n}= \lim_{n \to \infty} \frac{n+1}{n}=1.\]

\noindent {\bf Example 2:}\\
Assume $\uE=(E_1,E_2)$ on $\CH_E=l^2(\N) \oplus \C f$ for a unit
vector $f,$ given by
\[
E_i=\left(
\begin{array}{cc}
C_i  & 0 \\
B_i  & A_i \\
\end{array}
\right)
\]
for $i=1,2.$ The subspaces $\CH_C=l^2(\N)$ and $\CH_A=\C f$ give the
decomposition of $\CH_E$ relative to the above block matrix form.
Denote the shift operator on $l^2(\N)$ by $S.$ Take
\[C_1=C_2= \frac{1}{\sqrt{2}} S^*.\]
We fix a real number $0 < t < 1$ and then take
\[B_1=B_2=\frac{1}{\sqrt{2}}(\sqrt{1-t^2},0,0,\ldots),\]
\[A_1=A_2=\frac{t}{\sqrt{2}}.\]
Clearly $\uE$ is a coisometric lifting of $\uA.$ Because $\uA$ is
$*$-stable, this lifting is also subisometric. The defect space
$\CD_{*,A}= \C f.$ The isometry $\gamma: \CD_{*,A} \to \CD_C$ turns
out to be
\[ \gamma (k f)= \frac{k}{\sqrt{2}}\left(
\begin{array}{c}
(1,0,0,\ldots)^T \\
(1,0,0,\ldots)^T   \\
\end{array}
\right)\] for $k \in \C.$ Superscript $T$ denotes transpose. Finally
the constrained characteristic function is given for
$h=(h_1,h_2,\ldots) \in \CH_C$ by
\begin{eqnarray*}
\theta_{J^s,C,E}d^1_h &=& e_0 \otimes \left(
\begin{array}{c}
((\frac{t^2+1}{2})h_1,-\frac{h_2}{2} ,-\frac{h_3}{2} ,\ldots)^T \\
 ((\frac{t^2-1}{2})h_1,-\frac{h_2}{2} ,-\frac{h_3}{2} ,\ldots)^T \\
\end{array}
\right)\\
&&- \sum_{|\alpha| \geq 1} L^{J^s}_\alpha e_0 \otimes
\frac{(1-t^2)t^{|\alpha|}h_1}{2}  \left(
\begin{array}{c}
(1,0,0,\ldots)^T \\
(1,0,0,\ldots)^T   \\
\end{array}
\right)\\
\theta_{J^s,C,E}d^2_h &=& e_0 \otimes \left(
\begin{array}{c}
((\frac{t^2-1}{2})h_1,-\frac{h_2}{2} ,-\frac{h_3}{2} ,\ldots)^T \\
 ((\frac{t^2+1}{2})h_1,-\frac{h_2}{2} ,-\frac{h_3}{2} ,\ldots)^T \\
\end{array}
\right)\\
&&- \sum_{|\alpha| \geq 1} L^{J^s}_\alpha e_0 \otimes
\frac{(1-t^2)t^{|\alpha|}h_1}{2}  \left(
\begin{array}{c}
(1,0,0,\ldots)^T \\
(1,0,0,\ldots)^T   \\
\end{array}
\right)
\end{eqnarray*}
and for $k\in \CH_A, \hspace{.2cm} i=1,2 \hspace{.2cm}$ by
\begin{eqnarray*}
\theta_{J^s,C,E}d^i_k &=& - e_0 \otimes
 \frac{t(\sqrt{1-t^2})k}{2} \left(
\begin{array}{c}
(1,0,0,\ldots)^T \\
(1,0,0,\ldots)^T   \\
\end{array}
\right)\\
&&+ \sum_{\alpha} L^{J^s}_i L^{J^s}_\alpha e_0 \otimes
\frac{t^{|\alpha|}(\sqrt{1-t^2})(1-\frac{t^2}{2})k}{\sqrt{2}} \left(
\begin{array}{c}
(1,0,0,\ldots)^T \\
(1,0,0,\ldots)^T   \\
\end{array}
\right)\\
&&+ \sum_{j \neq i} \sum_{\alpha} L^{J^s}_j L^{J^s}_\alpha e_0
\otimes \frac{t^{|\alpha|+2}(\sqrt{1-t^2})k}{2 \sqrt{2}} \left(
\begin{array}{c}
(1,0,0,\ldots)^T \\
(1,0,0,\ldots)^T   \\
\end{array}
\right)
\end{eqnarray*}

We calculate the invariants of the corresponding Hilbert module like
in the previous example and as expected they are zero.
\[ curv_s A= 2!  \ \lim_{n \to \infty} \frac{\mbox{trace}(\eins - \sum_{|\alpha|=n+1} A_\alpha A^*_\alpha)}{n^2}= 2! \lim_{n \to \infty}
\frac{1 - t^{2(n+1)}}{n^2} =0\]
\[ \chi_s (A)=  2! \lim_{n \to \infty}\frac{\mbox{rank}(\eins - \sum_{|\alpha|=n+1} A_\alpha A^*_\alpha)}{n^2}=
2! \lim_{n \to \infty} \frac{1}{n^2} =0.\]

\noindent {\bf Example 3:}\\
The following example is for the noncommutative case. Let $\uE$ on
$\CH_E=\C \oplus \Gamma (\C^2)$ be a coisometric lifting of
$\uC=(1,0)$ on $\CH_C=\C.$ We take $B_1=0$ and $B_2(k)=k e_0,$ and
$\uA=(L_1,L_2),$ i.e., the tuple of creation operators. Thus $D_{*,A}$
is the projection onto the vacuum $\C e_0.$ Clearly for $h_1, h_2
\in \CH_C$ we have $D_C(h_1, h_2)=h_2.$ It follows that
$\gamma : \Gamma (\C^2) \to \C$ takes $e_0$ to 1 and maps $\Gamma
(\C^2) \ominus \C e_0$ to $0.$ The characteristic function is zero.
Finally using Equations (2.16) and (4.1) of \cite{Po02} we get

\[ curv A=    \lim_{m \to \infty} \frac{\mbox{trace}(\eins - \sum_{|\alpha|=m} L_\alpha L^*_\alpha)}{2^m}=  \lim_{m \to \infty}
\frac{1 +2+2^2+...+2^{m-1} }{2^m} =1\]
\[ \chi (A)=   \lim_{m \to \infty}\frac{\mbox{rank}(\eins - \sum_{|\alpha|=m} L_\alpha L^*_\alpha)}{2^m}=
 \lim_{m \to \infty} \frac{1 +2+2^2+...+2^{m-1} }{2^m} =1.\]

Instead if one chooses $\uA$ to be
$S_{\bot}$ constructed in  Theorem 3.5, 3.8 and 4.9 of \cite{Po02},
then one can realise any $t$ in $[0,1]$ as curvature invariant and Euler
characteristic.

\section[Covariant representations, Hardy algebras]
{Liftings of covariant representations of $W^*$-correspondences and
Hardy algebras}

First we prefer to remark that any row contraction
$\uT=(T_1,\ldots,T_d)$ on a Hilbert space $\CH$ can be encoded as
the covariant representation $(T,\sigma)$ of the
$W^*$-correspondence $\C^d$ over the von Neumann algebra $\C$ on
$\CH.$ In this picture if $\{e_1,\ldots,e_d\}$ denote the standard
basis of $\C^d,$ then $T_i=T(e_i)$ for all $i.$ In the current
section we will generalize the theory of chapter 4 for covariant
representations of $W^*$-correspondence. The  constrained dilation
theory corresponding to Cuntz-Kreiger constraints of the last
chapter fits into this scheme but  the $q$-commuting dilation theory
does not. The reader may need to refer back to some notions defined
in Section 2.4.

Let $(C,\sigma_C)$ be  a contractive covariant representation of
$\CE$ on $\CH_C.$ Then a contractive covariant representation
$(E,\sigma_E)$  of $\CE$ on a Hilbert space $\CH_E \supset \CH_C$ is
called a {\em contractive lifting} of $(C,\sigma_C)$ if

\begin{itemize}
\item[(a)] $\CH_C$ reduces $\sigma_E$ and for all $a \in \CM$
\[ \sigma_E (a)|_{\CH_C} = P_{\CH_C} \sigma_E (a)|_{\CH_C} = \sigma_C (a).\]

\item[(b)] $\CH_C^\bot$ is invariant w.r.t. $E(\xi)$ for all $\xi \in \CE.$

\item[(c)] $P_{\CH_C} E(\xi)|_{\CH_C}=C(\xi)$ for all $\xi \in \CE.$

\end{itemize}

Set $A(\xi)=P_{\CH_C^\bot} E(\xi)|_{\CH_C^\bot}$ and $\sigma_A (a):=
\sigma_E (a) |_{\CH_C^\bot}$ for all $\xi \in \CE, a \in \CM.$
Observe that $(A,\sigma_A)$ is also a covariant representation of
$\CE.$ This definition of contractive lifting is equivalent to
assuming that $\tilde{E}$ is contractive and has the form
\[  \tilde{E}= \left(
\begin{array}{cc}
\tilde{C}  & 0 \\
\tilde{B}  & \tilde{A} \\
\end{array}
\right).\]

Note that if $(E,\sigma_E)$ is completely contractive then
$(C,\sigma_C)$ and $(A,\sigma_A)$ are also completely contractive.
This follows easily by passing to $\tilde{E}$ and using Lemma 2.4.3.

\begin{defn}
Let $(T,\sigma)$ be a completely contractive covariant (c.c.c. for
short) representation of $\CE$ on $\CH.$ An isometric dilation $(V,
\pi)$ of $(T,\sigma)$ is an isometric covariant representation of
$\CE$ on $\tilde{\CH} \supset \CH$ such that $(V,\pi)$ is a lifting
of $(T,\sigma).$ A {\em minimal isometric dilation} (mid) of
$(T,\sigma)$  is an isometric dilation $(V,\pi)$ on $\hat{\CH}$ for
which (as before)
\[\hat{\CH} = \overline{\mbox{span}} \{ V(\xi_1)\ldots V(\xi_n) h: h \in \CH , \xi_i \in \CE \}.\]
\end{defn}

Further one defines the {\em full Fock module} over $\CM$ as
\[ \CF (\CE) = \oplus^n_{i=0} \CE^{\otimes ^n} \]
where $\CE^{\otimes ^0}= \CM.$ We will write $\CF$ for $\CF(\CE)$ in
short. For a $\xi \in \CE$ the creation operator $L_\xi$ on $\CF
(\CE)$ is given by $L_\xi \eta = \xi \otimes \eta.$ We have an
induced homomorphism $\varphi^n$ from $\CM$ to $\CL (\CE^{\otimes
^n})$ which for each $a \in \CM$ is given by
\[  \varphi^n (a) (\xi_1 \otimes \xi_2 \otimes \ldots \otimes \xi_n)=(\varphi(a) \xi_1) \otimes
\xi_2 \otimes \ldots \otimes \xi_n.\] Take the operator
\[ \varphi_\infty (a) = diag (a, \varphi(a), \varphi^2 (a), \ldots)\]
on $\CF.$

A mid of a c.c.c. representation always exist and is unique up to
unitary equivalence (cf. \cite{MS05}). We give a brief sketch of the
proof: Given a c.c.c. representation $(T,\sigma)$ of a
$W^*$-correspondence $\CE$ on a Hilbert space $\CH,$ we consider the
associated $\tilde{T}:\CE \otimes \CH \to \CH.$ We set
$D_{*,T}:=(\eins - \tilde{T}\tilde{T}^*)^{\frac{1}{2}}$ (in
$B(\CH)$) and $ D_T:=(\eins - \tilde{T}^*\tilde{T})^{\frac{1}{2}}$
(in  $B(\CE \otimes_\sigma \CH)$). Let $\CD_{*,T}:=
\overline{\mbox{range~} D_{*,T}}$  and
$\CD_{T}=\overline{\mbox{range~} D_T}.$ The space of mid $(V,\pi)$
is
\[ \hat{\CH}=\CH \oplus  \CF  \otimes_{\sigma_1} \CD_T\]
where $\sigma_1$ is defined to be the restriction to $\CD_T$ of
$\varphi (.) \otimes \eins.$ Finally a representation $\pi$ of $\CE$
on $\hat{\CH}$  given by
\[\pi = \sigma \oplus \sigma_1^{\CF} \circ \varphi_\infty \]
with $(\sigma_1^{\CF} \circ \varphi_\infty)(a)= \varphi_\infty(a) \otimes \eins_{\CD_T}$ for $a \in \CM,$
and a linear map $V:\CE \to B(\hat{\CH})$  given in operator matrix
form by:

\[ V(\xi)= \left(
\begin{array}{cccc}

T(\xi) & 0& 0& \ldots \\
D_T(\xi \otimes .) & 0& 0 & \ldots \\
0 & \eins & 0 &  \\
0 & 0 & \eins \\
\vdots &&& \ddots \\

\end{array}
\right)\] together constitute a mid $(V,\pi)$  of $(T,\sigma).$ \qed

Note that $\tilde{V}$ will be the mid of $\tilde{T}.$ Moreover if
\[\tilde{T} \tilde{T}^*= \eins,\]
then $(T,\sigma)$ is said to be {\em coisometric.} It is also known
that the mid  $(V,\pi)$ is coisometric if and only if $(T,\sigma)$
is coisometric.

Let us define $\tilde{T}^n= \tilde{T}(\eins \otimes
\tilde{T}^{n-1})$ for $n\geq 1.$ We say $(T,\sigma)$ is {\em
$*$-stable} if $\lim_{n \to \infty} \tilde{T}^n (\tilde{T}^n)^*=0$
in SOT.

The algebra defined below is intrinsically related to $H^\infty$ which will become
apparent when we go through the examples listed after it. Because the characteristic functions and Poisson kernel
of Sz-Nagy and Foias are $H^\infty$ functions, this new algebra will be crucial for extending our theory of chapter
4 for $W^*$-correspondences.

\begin{defn}
The ultraweakly closed subalgebra of $\CL(\CF)$ generated by the
$L_\xi$'s and $\varphi_\infty (a)$'s is called the Hardy algebra of
$\CE$ and is denoted by $H^\infty (\CE).$
\end{defn}

From works of Muhly and Solel (\cite{MS05}) it is known that there
are 1-1 correspondences:
\begin{itemize}
\item[(a)] between completely contractive covariant representations $(T,\sigma)$ and contractive
$\tilde{T}$'s.
\item[(b)] between completely contractive $(T,\sigma)$ and its integrated form $\sigma \times T: H^\infty (\CE) \to B(\CH),$
if $\|\tilde{T}\| < 1.$ Here $\sigma \times T$ maps
\[ \varphi_\infty (a) \mapsto \sigma (a), L_\xi \mapsto T(\xi).\]
\end{itemize}

\begin{lem}
\[ \tilde{T} (\varphi (a) \otimes \eins)  = \sigma (a) \tilde{T}  \mbox{~~~~~~for~~} a \in  \CM.  \]
Therefore $\tilde{T}^* \tilde{T}$ commutes with $(\varphi
(\CM)\otimes \eins)$ and $\tilde{T} \tilde{T}^*$ commutes with
$\sigma (\CM).$
\end{lem}
\begin{proof} For $\xi \otimes h \in \CH^E$ and $a \in \CM$
we have
\begin{align*}
\tilde{T} (\varphi (a) \otimes \eins) (\xi \otimes h) = \tilde{T}
(\varphi
(a) \xi \otimes h)\\
= T (\varphi (a) \xi ) h =\sigma (a) T (\xi) h\\
=\sigma (a) \tilde{T}(\xi \otimes h).
\end{align*}
\end{proof}

We remark that each element $\mu \in \CE^\sigma$ there is a representation $(T,\sigma)$
of $\CE$ such that $\tilde{T}=\mu^*.$
Using the above bijective relations
we define for each $G \in H^\infty (\CE)$ a function $\hat{G}$ given
by
\[ \hat{G}(\mu^*):= (\sigma \times T) (G) ~~~~ \forall \mu=\tilde{T}^* \in \D(\CE^\sigma).\]
This $\hat{G}$ is called the {\em Fourier transform} of $G.$ Thus
elements of the Hardy algebra $H^\infty (\CE)$ can be realised as
functions on unit ball $\D(\CE^\sigma)$ analogous to classical
$H^\infty$ functions.
\\

Consider the following special case:

{\bf Example 1:}  When $\CM=\CE=\C$ and $\sigma$ is the identity
representation of $\CM$ on $\CH=\C,$ then  $\D (\CE^\sigma)$ is the open
unit disc in the complex plane. Any $G\in \CH^\infty (\CE)$ is
basically an infinite, lower-triangular, Toeplitz matrix on
$l^2(\N)$:

\[ G= \left(
\begin{array}{ccccc}

a_0 & 0& 0& 0&\ldots \\
a_1 & a_0& 0& 0&\ldots \\
a_2 & a_1 & a_0& 0&  \\
a_3 & a_2 & a_1& a_0&\\
\vdots & \vdots &&& \ddots \\

\end{array}
\right).\] The Fourier transform $\hat{G}:\D(\CE^\sigma)^* \to
B(\CH)$ is
\[\overline{\mu} \mapsto \sum^\infty_{i=0} a_i \overline{\mu}^i \]
with $\sum^\infty_{i=0} a_i \overline{\mu}^i$ acting via
multiplication on $\CH=\C.$

{\bf Example 2:} With the same $\CM=\CE=\C$ as above but $\CH$ a
Hilbert space instead, we take in this example the representation
$\sigma$ of $\C$ on $\CH$ as multiplication, i.e,
\[\sigma (c) h = c h, \hspace{1cm} c \in \C, h \in \CH. \]
Then $\CE^\sigma \cong B(\CH)$ and $\D(\CE^\sigma)$  is the set of
all strict contractions on $\CH.$ $\CH^\infty (\CE)$ is still the
set of all lower-triangular Toeplitz matrices as before. If $T \in
\D(\CE^\sigma),$ then
\[\hat{G}(T^*) =  \sum^\infty_{i=0} a_i (T^*)^i. \]

{\bf Example 3:} Next, consider the case when $\CM=\C,\CE=\C^n$ and
take $\sigma$ to be the representation of $\CM$ on $\CH$ same as in
Example 2. The $\D(\CE^\sigma)$ is the set of all row contractions
$\uT$ with norm $\|\uT \uT^*\|$ less than $1.$
\\

We can also define the Poisson kernel in module context. For every
$\mu \in \D(\CE^\sigma)$ we set an operator $\mu^{(n)}: \CH \to
\CE^{\otimes^n} \otimes_\sigma \CH$ by
\[ \mu^{(n)}:=(\eins_{\CE^{\otimes^{n-1}}} \otimes \mu)(\eins_{\CE^{\otimes^{n-2}}}\otimes \mu) \ldots
(\eins_{\CE} \otimes \mu)  \mu\] Now with it we associate the
operator called {\em Poisson kernel}
\[ K(\mu):=(\eins_\CF \otimes (\eins - \mu^* \mu )^\frac{1}{2})[\mu^{(0)},\mu^{(1)},\mu^{(2)},\ldots]^T.\]
 which maps $\CH$ to $\CF \otimes_\sigma \CH.$

Characteristic functions of covariant representations of $W^*$-
correspondences have been studied by Muhly and Solel (cf.
\cite{MS05}). Here we are interested in the corresponding theory for
liftings of covariant representations. We consider two special cases
of liftings as in the last chapter and then investigate the general
case.

Let $(V^E,\pi_E)$ and $(V^C,\pi_C)$ denote the mids of $(E,\sigma_E)$ and $(C,\sigma_C)$ respectively.
From the definition of lifting it is immediate that the space of the mid $V^C$ can be embedded as a subspace
reducing $V^E.$

\begin{defn}
A lifting $(E,\sigma_E)$ of a completely contractive covariant representation $(C,\sigma_C)$ on $\CH_E \supset \CH_C$
is called subisometric if the corresponding
mids $V^E$ and $V^C$ are unitarily equivalent, i.e., there exists a unitary $W: \hat{\CH}_E \to \hat{\CH}_C$
such that $W|_{\CH_C}=\eins|_{\CH_C},$ $W V^E(\xi) = V^C(\xi) W$ for all $\xi \in \CE$ and $W \pi_E (a) = \pi_C (a) W$ for all
$a \in \CM.$
\end{defn}

\begin{rem}
Alternatively, a subisometric lifting means the existence of a
unitary $W$ (same as above) such that
\[\tilde{V}^C (1 \otimes W) =W \tilde{V}^E.\]
\end{rem}

\begin{prop}
Let $(C,\sigma_C)$ be a completely contractive covariant (c.c.c. for
short) representation of $W^*$-correspondence  $\CE$ on $\CH_C.$ A
completely contractive lifting $(E,\sigma_E)$ on $\CH_E=\CH_C \oplus
\CH_A$ of  $(C,\sigma_C)$ with
\[ E(\xi)= \left(
\begin{array}{cc}
C(\xi) & 0\\
B(\xi) & A(\xi) \\
\end{array}
\right), \mbox{~~~~~} \xi \in \CE,
\]
is subisometric if and only if $(A,\sigma_A)$ is $*$-stable and $\tilde{B}=D_{*,A} \gamma^* D_C$ with an isometry
$\gamma:\CD_{*,A} \to \CD_C.$
\end{prop}

The {\em characteristic function} of subisometric lifting as before
is defined as
\[ M_{C,E}:= W|_{\CF \otimes \CD_E}.\]

\begin{thm}
Let $(C,\sigma_C)$ be a c.c.c. representation of $\CE$ on a Hilbert
space $\CH_C.$ Two subisometric liftings $(E,\sigma_E)$ and
$(E^\prime,\sigma_E)$ of $(C,\sigma_C)$ are unitarily equivalent if
and only if the the corresponding characteristic functions $M_{C,E}$
and $M_{C,E^\prime}$ are unitarily equivalent.
\end{thm}
\begin{proof} For the proof of the necessary part we assume that
 the liftings $(E,\sigma_E)$ and $(E^\prime,\sigma_{E^\prime})$ are c.c.c. representations on $\CH_E$
 and  $\CH_{E^\prime},$ and $U:\CH_E \to \CH_{E^\prime}$ is a unitary such that
 $U|_{\CH_C}=\eins_{\CH_C}$ and
\[
UE(\xi)=E^\prime(\xi)U, \hspace{1cm} U \sigma_E (a) = \sigma_{E^\prime}(a)
U \hspace{1cm} \forall \xi \in \CE, a \in \CM,
\]
The mids of unitarily equivalent row contractions are unitarily
equivalent. Hence we extend $U$ canonically to a unitary
$\hat{U}:\hat{\CH}_E \to \hat{\CH}_{E^\prime}$ defined between the
spaces of mids $(V^E,\pi_E)$ and $(V^{E^\prime},\pi_{E^\prime})$
with $\hat{U}|_{\CH_E} =U,$ and we get
\[
\hat{U} V^E(\xi)=V^{E^\prime}(\xi)\hat{U}, \hspace{1cm} \hat{U}
\pi_E (a) = \pi_{E^\prime}(a) \hat{U} \hspace{1cm} \forall \xi \in \CE, a
\in \CM
\]

Because $(E,\sigma_E)$ and $(E^\prime,\sigma_{E^\prime})$ are
subisometric we also have unitaries $W: \hat{\CH}_E \rightarrow
\hat{\CH}_C$ and $W^\prime: \hat{\CH}_{E^\prime} \rightarrow
\hat{\CH}_C$ respectively with:
\begin{align*}
 \quad V^C(\xi) W = W V^E(\xi), \quad W |_{\CH_C} = \eins |_{\CH_C}
\\
 \quad V^C(\xi) W^\prime = W^\prime V^{E^\prime}(\xi), \quad W^\prime
|_{\CH_C} = \eins |_{\CH_C} \nonumber
\end{align*}
Let us take
\[
U_C := W^\prime \hat{U} W^*: \hat{\CH}_C \rightarrow \hat{\CH}_C.
\]
Chasing a commuting diagram similar to diagram 4.24 and arguing on
the lines of Theorem 4.1.6  we find that the characteristic
functions are equivalent.

Conversely we show that if $\Theta = \Theta^\prime V$ with a unitary
$V: \CD_E \rightarrow \CD_{E^\prime}$ then the two subisometric
liftings $(E,\sigma_E)$ and $(E^\prime,\sigma_{E^\prime})$ are
unitarily equivalent. Let $W$ and $W^\prime$ be the corresponding
unitaries from the subisometric lifting property. Then
\begin{align*}
W^\prime \CH_{E^\prime} \;=\;  \CH_C \oplus ( \CF \otimes \CD_C )
\;\ominus\; W^\prime ( \CF \otimes \CD_{E^\prime} ) \\
= \CH_C \oplus ( \CF \otimes \CD_C )
\;\ominus\; M_{C,E^\prime} \big( \CF \otimes V \CD_E \big) \\
= \CH_C \oplus ( \CF \otimes \CD_C ) \;\ominus\; M_{C,E} ( \CF
\otimes \CD_E ) \;=\; W \CH_E,
\end{align*}
and we can define
\begin{align*}
U := (W^\prime)^* \, W |_{\CH_E}: \CH_E \rightarrow \CH_{E^\prime}.
\end{align*}

Because for $h \in \CH_C,$ $W h = h = W^\prime h$ we have $U h = h$.
In general for $h\in \CH_E$ and $\xi \in \CE$ (with $p_E,
p_{E^\prime}$ orthogonal projections onto $\CH_E, \CH_{E^\prime}$)
\begin{align*}
U E(\xi)\, h = (W^\prime)^*\, W \,E(\xi) \,h = (W^\prime)^* \,W
\,p_E \,V^E(\xi) \,h = p_{E^\prime}\, (W^\prime)^*\, W\, V^E(\xi)
\,h
\\
= p_{E^\prime}\, (W^\prime)^*\, V^C(\xi)\, W\, h = p_{E^\prime}\,
V^{E^\prime}(\xi)\, (W^\prime)^* \,W\, h = E^\prime(\xi) \,U\, h.
\end{align*}

Identical computation gives
\[ U \sigma_E (a) = \sigma_{E^\prime} (a) U \mbox{~~~~for~} a \in \CM.\]

Hence $\uE$ and $\uE^\prime$ are unitarily equivalent.
\end{proof}

Consider the coisometric liftings of $(C,\sigma_C)$ by $*$-stable  $(A,\sigma_A).$ Then the unitary equivalence
classes of those  which are also subisometric liftings of $(C,\sigma_C)$ are parametrized by isometries $\gamma: \CD_{*,A}
 \to \CD_C.$

Next we deal with the general case where $(E,\sigma_E)$ is a
contractive lifting of $(C,\sigma_C).$ Because of the structure of
lifting it is immediate that the space of mid $(V^C,\pi_C)$ is
embedded in $(V^E,\pi_E).$ We introduce a c.c.c. representation
$(Y,\pi_Y)$ on the orthogonal complement $\CK$ of the space of mid
$(V^C,\pi_C)$ to encode this. Hence we can get a unitary $W$ such
that
\begin{equation*}
W: \CH_E \oplus \big( \CF \otimes \CD_E \big)
\rightarrow
\CH_C \oplus \big( \CF \otimes \CD_C \big) \oplus \CK,
\end{equation*}
\begin{equation*}
\hat{V}^E(\xi) W = W V^E(\xi), \quad W |_{\CH_C} = \eins |_{\CH_C} \quad
\mbox{with} \quad \hat{V}^E(\xi) = V^C(\xi) \oplus Y(\xi).
\end{equation*}
Recall that $\CF$ denote the full Fock module on $\CE.$ We denote by
the same $W$ its restriction to the complement of $\CH_C$ too, i.e.,
\begin{align*}
W: \CH_A \oplus \big( \CF \otimes \CD_E \big)
\rightarrow
\big( \CF \otimes \CD_C \big) \oplus \CK
\end{align*}
With this we get
\begin{align*}
B(\xi)^* = p_{\CH_C} V^E(\xi)^* |_{W \CH_A} = p_{\CH_C} \big[V^C(\xi)^* \oplus Y(\xi)^* \big] W |_{\CH_A}\\
 = (D_C(\xi \otimes .))^*\, p_{\CD_C} W |_{\CH_A}.
\end{align*}
Doing computations on the lines of those appearing in Section 4.3 we
obtain
\[
P_{\CD_C} W |_{\CH_A} = \gamma D_{*,A}.
\]
We have shown that if $(E,\sigma_E)$ is a c.c.c. lifting of
$(C,\sigma_C)$ by $(A,\sigma_A)$ as above then
\begin{align}
\tilde{B}=D_{*,A} \gamma^* D_C.
\end{align}

For the converse we start with two c.c.c. representations
$(C,\sigma_C)$ and $(A,\sigma_A),$ a contraction $\gamma: \CD_{*,A}
\to \CD_C$ and $\tilde{B}$ as in equation (5.9). Then for $x \in
\CH_C,\; y \in \CH_A$
\begin{align*}
| \langle x, \tilde{C}\, \tilde{B}^* y \rangle |^2 = | \langle x, \tilde{C} D^*_C \gamma D_{*,A}\, y \rangle |^2
= | \langle D_C \tilde{C}^* x, \gamma D_{*,A}\, y \rangle |^2
\end{align*}
\begin{align*}
\leq \| D_C \tilde{C}^* x \|^2 \| \gamma D_{*,A} y \|^2
\leq \langle x, (\eins - \tilde{C} \tilde{C}^*) x \rangle \; \langle y, (\eins - \tilde{A} \tilde{A}^*) y \rangle.
\end{align*}
As in Exercise 3.2 in \cite{Pau03} it follows that
\begin{align*}
0 \leq
\left(
\begin{array}{cc}
\eins - \tilde{C} \tilde{C}^* & -\tilde{C} \tilde{B}^*\\
-\tilde{B} \tilde{C}^*  & \eins - \tilde{A} \tilde{A}^* \\
\end{array}
\right)
= \eins - \tilde{E} \tilde{E}^*.
\end{align*}
Thus $\tilde{E}$ is a contraction. We summarize this in the following lemma:

\begin{lem}
Let $(E,\sigma_E)$ be a lifting of $(C,\sigma_C)$ by $(A,\sigma_A).$ Then $(E,\sigma_E)$ is completely contractive
if and only if
$(C,\sigma_C)$ and $(A,\sigma_A)$ are  completely contractive and there exists a contraction $\gamma$
\[ \tilde{B}=D_{*,A} \gamma^* D_C.\]
\end{lem}

 It turns out that
\[ (\CF \otimes \CD_C) \vee W(\CF \otimes \CD_E) = (\CF \otimes \CD_C) \oplus \CK, \]

if
\begin{enumerate}

\item $\gamma$ is {\em resolving}, i.e., for $h \in \CH_A$
\begin{align*}
\big( \gamma D_{*,A} (A(\xi))^* h =0 \;\mbox{for all}\; \xi \in \CE \big)
\Rightarrow \\
\big( D_{*,A} (A(\xi))^* h =0 \;\mbox{for all}\; \xi \in \CE \big), \mbox{~and~}
\end{align*}

\item $(A,\sigma_A)$ is {\em c.n.c.}, i.e., $\CH^1_A:=\{h \in \CH_A: \|(\tilde{A}^n)^* h\|=\|h\| \mbox{~for all~}
n\in \N \}=\{0\}.$

\end{enumerate}
We call such liftings {\em reduced liftings}.

In this case the {characteristic function} of lifting is defined as
\[ M_{C,E}= P_{\CF \otimes \CD_C} W|_{\CF \otimes \CD_E}.\]
It can be shown as in the Chapter 4 that for any  completely
contractive covariant representation $(C,\sigma_C),$ these
characteristic functions are complete invariants for reduced
liftings of $(C,\sigma_C)$ up to unitary equivalence.
\\

We set $\mu = \tilde{A}^*$ and $d^\xi_h:=(V^E(\xi)- E(\xi))h= D_E
(\xi \otimes h)$ for $\xi \in \CE$ and $h \in \CH_E.$ The expanded
form of the characteristic function is the following:
\\

{\bf Case I:} $h \in \CH_C.$
\[
\Theta_{C,E}  (d^\xi_h) =
[D_C(\xi \otimes h) - \gamma \DsA B(\xi) h] - \sum^\infty _{j=1}  \gamma D_{*,A}
(L^*_\xi \tilde{A}^*)^j B(\xi) h.
\]

Alternatively
\[
\Theta_{C,E}  (d^\xi_h) =
[D_C(\xi \otimes h) - \gamma \DsA B(\xi) h] - \sum^\infty _{j=1}  (\eins \otimes \gamma D_{*,A})
(\eins \otimes \mu^j) (\eins \otimes B(\xi) h).
\]

{\bf Case II:} $h \in \CH_A.$

\[
 \theta_{C,E} (d^\xi_h)= - \gamma \tilde{A} d^\xi_h + \sum^\infty _{j=0} \gamma D_{*,A} (L^*_\xi \tilde{A}^*)^j
L^*_\xi D_A d^\xi_h.
\]

Alternatively
\[
\theta_{C,E} (d^\xi_h)= - \gamma \mu^* d^\xi_h + \sum^\infty _{j=0} (\eins \otimes \gamma D_{*,A})(\eins \otimes
\mu^j) D_A d^\xi_h.
\]

Let us briefly mention one  potential good application of this
theory to analytic crossed products of the type $\CM \rtimes \Z_+$
(cf. Section 6 of \cite{MS05}). Muhly and Solel showed in this last
quoted work that one can associate a contraction $t$ to any
($\sigma$-weakly continuous) representation of this crossed product.
When $t$ is c.n.c., its  Sz. Nagy- Foias characteristic function is
unitarily equivalent to the characteristic function of the covariant
representation associated to this representation of $\CM \rtimes
\Z_+.$ This theory needs to be explored for liftings of covariant
representations.

\section{Appendix}

We list some facts about Hilbert $C^*$-modules and some associated important
$C^*$-algebras. Given a Hilbert $C^*$-bimodule $\CE$ on a $C^*$-algebra $\CA$ with
the associated left action of $\CA$ denoted by $\varphi,$ the {\em
Toeplitz $C^*$-algebra} $\CT (\CE)$ of  $\CE$ is the
$C^*$-subalgebra generated by $\{L_\xi\}_{\xi \in \CE}$ and
$\{\varphi_\infty(a)\}_{a\in \CA}$ of $\CL(\CE).$

\begin{thm}
(Pimsner) If $\CE$ is a Hilbert $C^*$-bimodule on $\CA,$ then there
is a bijective map from the set of all isometric covariant
representations $(V,\sigma)$ from $\CE$ on $\CH$ to the set of all
$C^*$-representation from $\CT (\CE)$ on $\CH,$ determined  by
\[L_\xi \mapsto V(\xi), \hspace{1cm} \varphi_\infty (a) \mapsto \sigma (a).\]
\end{thm}

An ideal $\CI$ in a $C^*$-algebra $\CC$ is called {\em essential} if there is no nonzero ideal of $\CC$ that has
zero intersection with $\CI.$ For any algebra $\CA$ there always exists a unique (up to isomorphism) maximal $C^*$-algebra
that contains $\CA$ as an essential ideal. This maximal ideal is called the {\em multiplier algebra} of $\CA,$ denoted by $M(\CA).$
Set $\CB$ as the $C^*$-subalgebra of $\CL(\CF(\CE))$ generated by $\CL(\sum^N_{n=0} \CE^{\otimes^n})$ for all $N \in \N.$

\begin{defn}
A {\em Cuntz-Pimsner algebra} $\CO(\CE)$ is the image of $\CT(\CE)$ (under canonical embedding) in $M(\CB)/ \CB.$
\end{defn}

The Cuntz-Pimsner algebra generalizes Cuntz-Krieger algebra and
cross-product of any $C^*$-algebra by $\Z.$ The following are some
examples of Cuntz-Pimsner algebras:
\begin{itemize}
\item[(a)] For $\CA=\CE=\C$ we can identify the Fock module $\CF(\CE)$ with $l^2(\N).$ Hence
$\CT(\CE)= C^*(L)$ where $L$ is the unilateral shift and
\[\CO(\CE) \cong C(\T)\]
where $\T$ denotes the unit circle in the complex plane.

\item[(b)] Now take $\CA$ to be any unital $C^*$-algebra with an automorphism $\varphi.$ Then $\CA$ has a Hilbert module
structure $\CE$ with the left action  $\tilde{\varphi}:\CA \to
\CL(\CE)$ is
\[a \mapsto \varphi(a)\]
Here $\varphi(a)$s are realised as elements of $\CL(\CE)$ acting as
multiplication operator. We get $\CO(\CE) \cong \CA \times_\varphi \Z.$

\item[(c)] The Cuntz algebra $\CO_d$ can be realised as $\CO(\CE)$ for $\CA=\C, \CE=\C^d$ and
\[ \varphi (k) \xi = (k. \eins) \xi \mbox{~~~~~~for~~} k \in \CA, \xi \in \CE.\]

\end{itemize}

From a separable infinite dimensional Hilbert space $\CH$ and a
$C^*$-algebra $\CA$ we can construct an important example of Hilbert
module $\CH_\CA:=\CH \otimes \CA$  where the tensor product is given
by
\[ <h_1 \otimes a_1, h_2 \otimes a_2>=<h_1, h_2>a^*_1 a_2 \mbox{~~~~~~for~~} h_1, h_2 \in \CH; a_1, a_2 \in \CA.\]

Bellow some results are quoted which we need for Theorem 2.4.3. They show the usefulness of
$\CH_\CA.$
\begin{thm} (Kasparov's stabilisation theorem)
If $\CE$ is a countably generated Hilbert $\CA$-module, then $\CE \oplus \CH_\CA \cong \CH_\CA.$
 \end{thm}

\begin{cor}
If $\CE$ is a countably generated Hilbert $\CA$-module, then the $C^*$-algebra of compacts $K(\CE)$ is $\sigma$-unital, i.e., $K(\CE)$
has a countable approximate unit.
\end{cor}

\backmatter


\end{document}